\documentclass{amsart}

\usepackage{amsmath,amsthm, txfonts, graphicx}
\usepackage{curves}
\usepackage[colorlinks=true]{hyperref}
 \hypersetup{urlcolor=blue, citecolor=red}

\newtheorem{theorem}{Theorem}[section]
\newtheorem{lemma}[theorem]{Lemma}
\newtheorem{proposition}[theorem]{Proposition}
\newtheorem{corollary}[theorem]{Corollary}

\numberwithin{equation}{section}

\theoremstyle{definition}
\newtheorem{definition}[theorem]{Definition}
\newtheorem{notation}[theorem]{Notation}
\newtheorem{assumption}[theorem]{Assumption}

\theoremstyle{remark}
\newtheorem{remark}{Remark}
\newtheorem{example}{Example}

\DeclareMathOperator{\im}{Im}%
\DeclareMathOperator{\re}{Re}%
\DeclareMathOperator{\Arg}{Arg}%
\DeclareMathOperator{\dom}{dom}%
\DeclareMathOperator{\diam}{diam}%

\newcommand{\w}{\sigma}
\newcommand{\id}{I}% FIX THIS!!
\newcommand{\DF}{\ensuremath{E}}

\newcommand{\Sect}{A_{\alpha}}
\newcommand{\scale}{M}

\begin{document}

\title[Resolvent estimates and fractal blowups]{Estimates for the resolvent kernel of the Laplacian on p.c.f. self similar fractals and blowups.}
\author{Luke~G. Rogers}

\date{\today}

\maketitle

\section{Introduction}
One of the main features of analysis on post-critically finite self-similar (pcfss) sets is that it is possible to understand the behavior of the Laplacian and its inverse, the Green operator, in terms of the self-similar structure of the set.  Indeed, a major step in the approach to analysis on self-similar fractals via Dirichlet forms was Kigami's proof~\cite{Kigami1993TAMS,Kigami2003JFA} that for a self-similar Dirichlet form the Green kernel can be written explicitly as a series in which each term is a rescaling of a single expression via the self-similar  structure.  In~\cite{IoPeRoHuSt2010TAMS} this result was extended to show that the resolvent kernel of the Laplacian, meaning the kernel of $(z-\Delta)^{-1}$, can also be written as a self-similar series for suitable values of $z\in\mathbb{C}$. Part of the motivation for that work was that it gives a new understanding of functions of the Laplacian (such as the heat operator $e^{t\Delta}$) by writing them as integrals of the resolvent.

The purpose of the present work is to establish estimates that permit the above approach to be carried out.  We study the functions occurring in the series decomposition from~\cite{IoPeRoHuSt2010TAMS} (see Theorem~\ref{maintheoremfromIPRRS} below for this decomposition) and give estimates on their decay.  From this we determine estimates on the resolvent kernel and on kernels of operators defined as integrals of the resolvent kernel.  In particular we recover the sharp upper estimates for the heat kernel (see Theorem~\ref{heatkernelestimates}) that were proved for pcfss sets by Hambly and Kumagai~\cite{HamblKumag1999PLMS} by probabilistic methods (see also~\cite{BarloPerki1988PTRF,Kumagai1993PTRF,FitzHamKum1994CMP} for earlier results of this type on less general classes of sets).  It is important to note that the preceding authors were able to prove not just upper estimates but also lower bounds for the heat kernel, and therefore were able to prove sharpness of their bounds.  Our methods permit sharp bounds for the resolvent kernel on the positive real axis, but we do not know how to obtain these globally in the complex plane or how to obtain lower estimates for the heat kernel from them.  Therefore in this direction our results are not as strong as those obtained in~\cite{HamblKumag1999PLMS}.  However in other directions we obtain more information than that known from heat kernel estimates, and we hope that our approach will complement the existing probabilistic methods.  In particular we are able to obtain resolvent bounds on any ray in $\mathbb{C}$ other than the negative real axis (where the spectrum lies), while standard calculations from heat kernel bounds only give these estimates in a half-plane.

A further consequence of our approach is that we extend (in Theorem~\ref{mainblowuptheorem}) the decomposition from~\cite{IoPeRoHuSt2010TAMS} to the case of blowups, which are non-compact sets with local structure equivalent to that of the underlying self-similar sets.  The blowup of a pcfss set bears the same relation to the original set as the real line bears to the unit interval, see~\cite{Strichartz1998CJM} for details.

The structure of the paper is as follows.  In Section~\ref{preliminaries} we recall some basic features of analysis on pcfss sets, as well as the main result of~\cite{IoPeRoHuSt2010TAMS}, which is the decomposition of the resolvent as a weighted sum of piecewise eigenfunctions.  Section~\ref{partitionssection} then discusses the natural decomposition of the pcfss set according to a self-similar harmonic structure.  The next two sections deal with sharp estimates on the real axis.  In Section~\ref{estimatesofpiecewiseefnssection} we obtain estimates of piecewise eigenfunctions for which the eigenvalues are real.  To show these we decompose into pieces in the manner of Section~\ref{partitionssection}, prove that each piece has one normal derivative that is larger than the others by a factor, and use this to show that the only way to glue pieces together smoothly is by requiring the size to decay by this factor each time we cross a cell of the decomposition.   This implies that piecewise eigenfunctions have sub-Gaussian decay.  From these estimates it is routine, though somewhat long and technical, to estimate the kernel of $(\lambda-\Delta)^{-1}$ for $\lambda$ on the positive real axis and to prove that the decomposition of~\cite{IoPeRoHuSt2010TAMS} is valid on blowups in this setting. The former is in Section~\ref{resolventestimatesection} and the latter in Section~\ref{blowupsection}.

In Section~\ref{plsection} we switch gears and prove some complex analytic estimates related to the Phragmen-Lindel\"{o}f theorem.  Our goal is a method for obtaining resolvent estimates on sectors in the complex plane from our decay estimates on the real axis.  This method is not restricted to the setting of pcfss sets, and may be of interest for proving resolvent estimates in more general settings, such as metric measure spaces.  We find that combining decay estimates on the real axis with some weak estimate on a sector away from the spectrum proves decay on any ray in the sector.  In Section~\ref{offnegaxissection} we show that the required weak estimates can be proved on pcfss sets by a generalization of some arguments from Chapter~4 of~\cite{Kigami2001}.  As a result we obtain our main result, Theorem~\ref{Allestimatesoffnegaxis}, giving decay estimates for piecewise eigenfunctions and the resolvent kernel away from the negative real axis in $\mathbb{C}$.  Section~\ref{otherkernelsection} contains some examples, including a proof of the upper bounds for the heat kernel.

\section{Acknowledgements}
This paper relies heavily on its precursor~\cite{IoPeRoHuSt2010TAMS}; the author thanks his co-authors on that paper, especially Bob Strichartz, whose suggestion that there might be an approach to heat kernel estimates via Kigami's theory and the results of~\cite{IoPeRoHuSt2010TAMS} inspired the present work.

\section{PCFSS Fractals: Energy, Laplacian, Metric, Eigenvalue Counting and Resolvent kernel}\label{preliminaries}

We work on a pcfss fractal $X$ corresponding to an iterated function system $\{F_{1},\dotsc
F_{J}\}$.  In full generality this could be a compact metrizable topological space $X$ equipped
with continuous injective self-maps $F_{j}$ so that there is a continuous surjection from the space
of infinite words over $\{1,\dotsc,J\}^{\mathbb{N}}$ to $X$, as in Chapter~1 of \cite{Kigami2001}.
The reader may prefer, however, to think of the more intuitive situation in which the $F_{j}$ are
Lipschitz contractions on a finite dimensional Euclidean space and $X$ is their fixed point in the
sense of Hutchinson~\cite{Hutchinson1981IUMJ}.  In either case, the boundary of $X$ is the post-critical
set $V_{0}$, which is a finite set with the property that $F_{j}(X)\cap F_{k}(X)\subset
F_{j}(V_{0})\cap F_{k}(V_{0})$ for $k\neq j$.  We assume that $X\setminus V_{0}$ is connected. The
best known example of a set of this type is the Sierpinski Gasket, for which the full analytic
theory we will use is expounded in~\cite{Strichartz2006}.  Note that the Sierpinski Carpet is not
of this type, as it is not post-critically finite.

If $w$ is a finite word on the letters $\{1,\dotsc,J\}$ then $|w|$ is its length. The set $V_{m}$
is $\cup_{|w|=m} F_{w}(V_{0})$ and $V_{\ast}=\cup_{m=0}^{\infty} V_{m}$. Points in
$V_{\ast}\setminus V_{0}$ are called junction points.  For a word $w=w_{1}\dotsm w_{m}$ we write
$F_{w}=F_{w_{1}}\circ F_{w_{2}}\circ\dotsm\circ F_{w_{m}}$.  The set of finite words is $W_{\ast}$,
and the set of infinite words is $\Sigma$. If $\w=\w_{1}\w_{2}\dotsm$ is an infinite word we use
$[\w]_{m}$ to denote the subword $\w_{1}\dotsm \w_{m}$. For a finite word $w$ we use the notation
$wW_{\ast}$ (respectively $w\Sigma$) for finite (respectively infinite) words that begin with $w$.

We assume that we have a regular self-similar harmonic structure on $X$ which provides an irreducible
self-similar Dirichlet form $\DF$ with domain $\dom(\DF)$ and scalings $0<r_{j}<1$ so that
\begin{equation*}
    \DF(u) = \sum_{1}^{J} r_{j}^{-1} \DF(u\circ F_{j}).
    \end{equation*}
Details of the construction and properties of such forms may be found in~\cite{Kigami2001}. We also
fix a probability measure $\mu$ satisfying
\begin{equation}\label{muscaling}
    \mu(A) =\sum_{1}^{J} \mu_{j} \mu(F_{j}^{-1}(A)).
    \end{equation}
We will require that this measure is related to the Dirichlet form by $\mu_{j}=r_{j}^{S}$ for the
unique $S$ so that $\sum_{j} r_{j}^{S}=1$.  This measure is known to be the natural one for many
aspects of the analysis on $X$; it arises in studying the analogue of Weyl-type asymptotics of
Laplacian eigenvalues~\cite{KigamLapid1993CMP}, and in determining heat kernel
estimates~\cite{Kigami2001,Kumagai1993PTRF,FitzHamKum1994CMP,HamblKumag1999PLMS}.  It is the Hausdorff
measure (with Hausdorff dimension $S$) for the resistance metric, which is the metric in which the
distance between points $x$ and $y$ is $R(x,y)$ given by
\begin{equation*}
    R(x,y)^{-1} = \min \bigl\{ \DF(u): u(x)=0,\, u(y)=1 \bigr\}.
    \end{equation*}

\begin{remark}
In most of the theory of analysis on fractals one can take any Bernoulli measure defined
by~\eqref{muscaling} for some $0<\mu_{j}<1$ with $\sum_{j}\mu_{j}=1$.  In particular, this is all
that is required for the construction of the Laplacian resolvent in~\cite{IoPeRoHuSt2010TAMS}.  The author would be curious to know whether it suffices for Theorem~\ref{mainblowuptheorem} or even Theorem~\ref{pathdecayforGonblowupthm} regarding the resolvent on blowups.
\end{remark}

The Dirichlet form and measure give rise to a weak Laplacian $\Delta$ by defining $u\in\dom(\Delta)$ with
$\Delta u=f$\/ if there is a continuous $f$ such that
\begin{equation*}
    \DF(u,v)= -\int f v \, d\mu
    \end{equation*}
for all $v\in\dom_{0}(\DF)$, the functions in $\dom(\DF)$ that vanish on $V_{0}$.  The Laplacian is
self-adjoint and has compact resolvent, with negative eigenvalues $-\lambda_{j}$ of finite
multiplicity that accumulate only at $-\infty$.  The asymptotic distribution of the eigenvalues was
determined by Kigami and Lapidus~\cite{KigamLapid1993CMP}, who proved a more general version of the
following result.
\begin{proposition}[\protect{\cite{KigamLapid1993CMP}}]\label{KigamiLapidusWeylest}
Let $N(x)=\#\{\lambda_{j}:|\lambda_{j}|\leq x\}$.  Then
\begin{equation*}
    0< \liminf_{x\rightarrow\infty} \frac{N(x)}{x^{S/(S+1)}}
    \leq \limsup_{x\rightarrow\infty} \frac{N(x)}{x^{S/(S+1)}} <\infty
    \end{equation*}
\end{proposition}

We also recall from~\cite{Kigami2001} that there is an explicit formula for a continuous Green
kernel that is positive on the interior of $X$, zero on $V_{0}$, and inverts $-\Delta$ with
Dirichlet boundary conditions. Kigami's construction was generalized in~\cite{IoPeRoHuSt2010TAMS}, from which we
will need the following results.
\begin{proposition}\label{propexistofeta}
If $p\in V_{0}$ and $z\in\mathbb{C}$ is such that none of the values $r_{w}\mu_{w}z$ is a Dirichlet eigenvalue of
$\Delta$ then there is a function $\eta_{p}^{(z)}$ with measure-valued Laplacian satisfying
\begin{equation*}
\begin{cases}
    ( z\id - \Delta )\eta_{p}^{(z)} = 0 &\quad \text{on $X\setminus V_{0}$}\\
    \eta_{p}^{(z)} = \delta_{pq} &\quad\text{for $q\in V_{0}$}
    \end{cases}
\end{equation*}
where $\delta_{pq}$ is the Kronecker delta.  We call this a piecewise $z$-eigenfunction.
\end{proposition}

Under the same assumptions on $z$ we define functions that form a natural basis for the
piecewise $z$-eigenfunctions on $1$-cells.  For $p\in V_{1}\setminus V_{0}$, let
\begin{equation}\label{definitionofpsi}
    \psi^{(z)}_{p} = \sum_{j} \eta^{(r_{j}\mu_{j}z)}_{F_{j}^{-1}p} \circ F_{j}^{-1}
    \end{equation}
where by convention the sum ranges only over those terms that are well-defined, in this case those
$j$ so $p\in F_{j}(V_{0})$.  Then $\psi_{p}^{(z)}$ has a measure-valued Laplacian, with $(
z\id - \Delta )\psi^{(z)}_{p} =0$ on all $1$-cells and Dirac masses at the points of
$V_{1}\setminus V_{0}$.  The main result of~\cite{IoPeRoHuSt2010TAMS} is that for suitable values of $z$
the resolvent $G^{(z)}(x,y)$ of the Laplacian with Dirichlet boundary conditions may be
written in terms of these piecewise eigenfunctions.
\begin{proposition}
For each $z\in\mathbb{C}$ such that none of the values $r_{w}\mu_{w}z$ is a Dirichlet eigenvalue of
$\Delta$, let
\begin{align}\label{maintheoremfromIPRRS}
    G^{(z)}(x,y)
    &= \sum_{w\in W_{\ast}} r_{w} \Psi^{(r_{w}\mu_{w}z)}(F_{w}^{-1}x,F_{w}^{-1}y)\\
\intertext{in which}
    \Psi^{(z)}(x,y)
    &= \sum_{p,q\in V_{1}\setminus V_{0}} G^{(z)}_{pq}
    \psi^{(z)}_{p}(x) \psi^{(z)}_{q}(y) \label{PsiintermsofGpq}
    \end{align}
where $G^{(z)}_{pq}$ is the inverse of the symmetric matrix $B^{(z)}_{pq}$ with entries
\begin{equation}\label{defnofBpq}
    B^{(z)}_{pq}
    =\sum_{j:F_{j}(V_{0})\ni q} \partial_{n}^{F_{j}(X)}\psi^{(z)}_{p}(q).
    \end{equation}
Then $G^{(z)}(x,y)$ is continuous and vanishes on $V_{0}$, and it inverts
$(z\id-\Delta)$ in that if $f\in L^{1}(\mu)$ then
\begin{equation*}
 (z\id-\Delta)\: \int_{X} G^{(z)}(x,y)f(y)\,d\mu(y)=f(x).
 \end{equation*}
\end{proposition}

In what follows, we will suppose that $z=\lambda>0$  and study the decay of $\eta^{(\lambda)}$ and its dependence on
$\lambda$, so as to estimate the decay of the terms in $G^{(\lambda)}$.  This will then be used to study the dependence of  $\eta^{(z)}$ and $G^{(z)}$ for those $z\in\mathbb{C}$ that are not on the negative real axis.   As usual
these estimates will include constants depending on the structure of the fractal, for which reason
we define a shorthand notation as follows.

\begin{notation}
We write $a\lesssim b$ or $b\gtrsim a$ if there is a constant $c>0$ that depends only on the
fractal, the harmonic structure or the measure, and for which $a\leq cb$.  If $a\lesssim b$ and
$b\lesssim a$ then we write $a\simeq b$.  The implicit constant, as well as constants explicitly
named, may vary from line to line in a computation.
\end{notation}

\section{Resistance Partitions and Chemical Paths}\label{partitionssection}

The short-time behavior of diffusion on fractal sets may be analyzed by partitioning the fractal such that all pieces have a prescribed resistance (up to a constant factor) and studying the lengths of minimal paths in this decomposition (see, for example, \cite{HamblKumag1999PLMS}), which are often called chemical paths.  We should therefore expect that the behavior of the resolvent $G^{(\lambda)}(x,y)$ for large $\lambda>0$ may be determined in the same manner.  In this section we record some known estimates for chemical paths that will be useful later.

\begin{definition}
A partition of $X$ is a finite set $\Theta\subset W_{\ast}$ with the properties
\begin{enumerate}
    \item If $\theta,\theta'\in\Theta$ then $\theta\Sigma\cap \theta'\Sigma\neq\emptyset$ only if $\theta=\theta'$,
    \item $\cup_{\theta\in\Theta} \theta\Sigma=\Sigma$.
    \end{enumerate}
We will refer to the cells $F_{\theta}(X)$, $\theta\in\Theta$, as the cells of $\Theta$.
\end{definition}

The partitions we use are
\begin{equation*}
    \Theta_{k}
    = \Bigl\{ \theta_{1}\dotsm\theta_{m}:
        r_{\theta_{1}}\dotsm r_{\theta_{m}} \leq e^{-k} < r_{\theta_{1}}\dotsm r_{\theta_{m-1}}
        \Bigr\},
    \end{equation*}
and we write $V(\Theta_{k})$ for the junction points corresponding to words in $\Theta_{k}$
\begin{equation*}
    V(\Theta_{k}) = \bigcup_{\theta\in\Theta_{k}} F_{\theta}(V_{0}).
    \end{equation*}
We will frequently use that when $\theta\in\Theta_{k}$,
\begin{equation}\label{comparabilityrelationsforpartition}
    r_{\theta}\simeq e^{-k},\quad \mu_{\theta}\simeq e^{-kS},\quad r_{\theta}\mu_{\theta}\simeq
    e^{-k(S+1)}.
    \end{equation}

There is a graph structure corresponding to the partition $\Theta_{k}$.  Let $\Gamma_{k}$ be the
graph with vertices $V(\Theta_{k})$ and edges between every $x$ and $y$ for which there is
$\theta\in\Theta_{k}$ such that $x,y\in F_{\theta}(V_{0})$. Points joined by an edge are said to be
neighbors. As a consequence of the definition of $\Theta_{k}$ we have (see Lemma~3.2 of
\cite{HamblKumag1999PLMS} for a proof)
\begin{equation}\label{resistcomparabletopartitionscale}
    R(x,y) \simeq e^{-k} \text{ for any $x,y\in V(\Theta_{k})$ that are neighbors in $\Gamma_{k}$.}
    \end{equation}

A path joining $x$ and $y$ in $\Gamma_{k}$ is a finite sequence $\{x_{i}\}_{i=0}^{I}$ of vertices such
that adjacent vertices are connected by an edge and $x_{0}=x$, $x_{I}=y$.  The length of the path
is $I$ and the graph distance $d_{k}(x,y)$, also called the chemical distance, is the length of the shortest path joining $x$ and $y$
in $\Gamma_{k}$.

It is generally difficult to obtain good estimates relating the graph distance $d_{k}(x,y)$ and the
resistance $R(x,y)$; the interested reader is directed to~\cite{Kigami2009MAMS} for a detailed analysis of this question and its connection to heat kernel estimates.  We will satisfy ourselves with some elementary but crude results.  For example, it is
well known that $R(x,y)$ does not exceed the resistance along a path from $x$ to $y$, so that for $x,y\in V(\Theta_{k})$ we must have
$R(x,y)\lesssim e^{-k}d_{k}(x,y)$.  In particular if $x,y\in V_{0}$ then $d_{k}(x,y)\gtrsim e^{k}$. Conversely, in Lemma~3.3 of~\cite{HamblKumag1999PLMS} it is shown that for any $x,y\in X$, the bound $d_{k}(x,y)\lesssim e^{(S+1)k/2}$ holds.  We will need these estimates in a form that compares $d_{k+k'}(x,y)$ to $d_{k}(x,y)$.
\begin{lemma}\label{boundonratioofchemicaldistancesfordifferentscales}
\begin{equation*}
    e^{k'} \lesssim  \frac{d_{k+k'}(x,y) }{d_{k}(x,y)}  \lesssim e^{(S+1)k'/2}
\end{equation*}
\end{lemma}
\begin{proof}
Observe that the partition by words
\begin{equation*}
    \bigl\{\theta''=\theta\theta':\theta\in\Theta_{k}, \theta'\in\Theta_{k'}\bigr\}
    \end{equation*}
has $r_{\theta''}\leq e^{-(k+k')}$, so is a subpartition of $\Theta_{k+k'}$ and consequently has longer paths.  Since a path in this partition describes a path $\Gamma_{k}$ in which any pair of vertices is joined by a path no longer than the maximal path between points of $V_{0}$ in $\Gamma_{k'}$, which we know is at most $e^{(S+1)k'/2}$, we conclude that $d_{k+k'}(x,y)\lesssim e^{(S+1)k'/2}d_{k}(x,y)$.

For the lower bound we note that there is a constant $c$ depending only on the harmonic structure and such that the partition by words
\begin{equation*}
    \bigl\{\theta''=\theta\theta':\theta\in\Theta_{k}, \theta'\in\Theta_{k'-c}\bigr\}
    \end{equation*}
contains $\Theta_{k+k'}$ and therefore has shorter paths.  Again each such path restricts to a path in $\Gamma_{k}$, and now every pair of $\Gamma_{k}$ vertices is separated by a path of length at least the minimal distance between points of $V_{0}$ in $\Gamma_{k'-c}$, which is $\gtrsim e^{k'-c}\gtrsim e^{k'}$.
\end{proof}

For the special classes of nested and affine nested fractals there are quite precise results about chemical distances in~\cite{Kumagai1993PTRF,FitzHamKum1994CMP}.  In particular the authors construct a geodesic metric on affine nested fractals that is comparable to a
power of the resistance (see Proposition~3.6 and Remark~3.7(2) of~\cite{FitzHamKum1994CMP}),
from which the following is easily obtained.

\begin{proposition}[\protect{\cite{FitzHamKum1994CMP}}]\label{affinenesteedversionofdk}
For an affine nested fractal there is $\gamma\in \Bigl[\frac{1}{S+1}, \frac{1}{2}\Bigr]$ which may be
obtained by solving an explicit optimization problem, such that for all sufficiently large $k$,
\begin{equation*}
    d_{k}(x,y)
    \simeq e^{k\gamma} R(x,y)^{\gamma}.
    \end{equation*}
\end{proposition}
\begin{proof}[Sketch of Proof]
Knowing that $R(x,y)^{\gamma}$ is comparable to a geodesic metric we see that edges of $\Gamma_{k}$
have geodesic length like $e^{-k\gamma}$, so the geodesic distance from $x$ to $y$ is comparable to
$d_{k}(x,y)e^{-k\gamma}$ if $k$ is large.  This must be comparable to $R(x,y)^{\gamma}$, so we have
$d_{k}(x,y) \simeq e^{k\gamma} R(x,y)^{\gamma}$.  The bounds on $\gamma$ are evident from the
preceding discussion.
\end{proof}

The proof in~\cite{FitzHamKum1994CMP} is a generalization of an earlier argument,
Proposition~3.5 of~\cite{Kumagai1993PTRF}, for the case of nested fractals.  We do not wish to give the
definitions of nested or affine nested fractals here, but we recall that they are subsets of
$\mathbb{R}^{n}$, are generated by iterated function systems consisting of Euclidean similarities
that have a high degree of symmetry.  Full details may be found
in~\cite{Kumagai1993PTRF,FitzHamKum1994CMP} and the references therein.

\begin{definition}
For later use we extend the definition of $d_{k}(x,y)$ from $\Gamma_{k}$ to all of $X$ in the
obvious fashion. For $x$ in the interior of a cell $F_{\theta}(X)$ of $\Theta_{k}$ and
$y\in\Gamma_{k}$ let $d_{k}(x,y)=\min\{d_{k}(z,y):z\in F_{\theta}(V_{0})\}$.  For $x$ in the
interior of $F_{\theta}(X)$ and $y$ in the interior of $F_{\theta'}(X)$ set $d_{k}(x,y)=\min
\{d_{k}(z,z'): z\in F_{\theta}(V_{0}), z'\in F_{\theta'}(V_{0})\}$. When we refer to a $d_{k}$ geodesic between points
$x$ and $y$ that are not in $\Theta_{k}$ we mean a geodesic joining the cells containing them.
\end{definition}

%%%%%%%%%%%%%%%%%%%%%%%%%%%%%%%%%%%%%%%%%%%%%%%%%%%%%%%%%%%%%%%%%%%%%%%%%%%%%%%%%%%%%%%%%%%%%%%%%%%%%%%%%%%%%%%%%%%%%%%%%%%%%%%%%%%%%%%%%%%%%%%

\section{Estimates for piecewise eigenfunctions with positive eigenvalue}\label{estimatesofpiecewiseefnssection}

In this section we develop decay estimates for the piecewise $\lambda$-eigenfunctions  $\eta_{p}^{(\lambda)}$ in the case
where $\lambda$ is a positive real number.  We summarize the results in the following theorem.
\begin{theorem}\label{maintheoremoneta}
There is $k(\lambda)$ with $e^{k(\lambda)}\simeq (1+\lambda)^{1/(S+1)}$ and a constant $c>0$ depending only on the fractal and harmonic structure such that for $\lambda\in(0,\infty)$ and $p\in V_{0}$, each piecewise $\lambda$-eigenfunction $\eta^{(\lambda)}_{p}$ satisfies the following bounds.  For all $q\in V_{0}$, $q\neq p$, and all $x\in X$
\begin{gather}
    \exp \Bigl( - c d_{k(\lambda)}(p,x)  \Bigr)
        \lesssim \eta^{(\lambda)}_{p}(x)
        \lesssim \exp \Bigl( - d_{k(\lambda)}(p,x)  \Bigr), \label{finalboundsforeta}\\
    \partial_{n}\eta^{(\lambda)}_{p}(p) \simeq (1+\lambda)^{1/(S+1)}, \label{finalboundfordiagonalnormalderivofeta}\\
    -(1+\lambda)^{1/(S+1)} \exp \Bigl( - d_{k(\lambda)}(p,q)  \Bigr)
        \lesssim \partial_{n}\eta^{(\lambda)}_{p}(q)
        \lesssim -(1+\lambda)^{1/(S+1)} \exp \Bigl( -  c d_{k(\lambda)}(p,q)  \Bigr), \label{finalboundforoffdiagnormalderiveta}
    \end{gather}
except that the lower bound of~\eqref{finalboundsforeta} is not valid on cells $F_{\theta}(X)$, $\theta\in\Theta_{k(\lambda)}$ such that $F_{\theta}(V_{0})$ contains a point of $V_{0}\setminus\{p\}$.
\end{theorem}
\begin{proof}
The number $k(\lambda)$ is introduced in Definition~\ref{defnofklambda}, where it is shown it is comparable to $(1+\lambda)^{1/(S+1)}$.  The upper bound of~\eqref{finalboundsforeta} for sufficiently large $\lambda$ is Corollary~\ref{etaexpdecayinchemicalmetric}, and continuity ensures we may take a suitably large constant multiple to make it true for all $\lambda>0$. Both~\eqref{finalboundfordiagonalnormalderivofeta} and the lower bound of~\eqref{finalboundforoffdiagnormalderiveta} are in Lemma~\ref{normalderivativeestimate}, while the lower bound of~\eqref{finalboundsforeta} is Corollary~\ref{lowerdecayestimateforeta} and the upper bound of~\eqref{finalboundforoffdiagnormalderiveta} is Corollary~\ref{lowerdecayboundcorolforeta}
\end{proof}

\begin{remark}
We note in passing that the estimate~\eqref{finalboundsforeta} is sufficient to complete Strichartz's proof of the smoothness of finite-energy harmonic functions on products of fractals~\cite[Theorem 11.4]{Strichartz2005TAMS}.  Strichartz proves this modulo an assumption on the behavior of the normal derivatives of the heat kernel~\cite[Equation (8.3)]{Strichartz2005TAMS}, but the key estimate in his proof is actually the decay estimate~\cite[Equation (8.25)]{Strichartz2005TAMS}, which is a consequence of the above in the case of nested fractals by Proposition~\ref{affinenesteedversionofdk}.  More generally, for products based on fractals and harmonic structures of the type discussed in this paper, the estimate~\eqref{finalboundsforeta} is sufficient to imply the convergence in~\cite[Equation (8.19)]{Strichartz2005TAMS}, and thus the argument proving~\cite[Theorem 11.4]{Strichartz2005TAMS}.  The author hopes that is a first step toward proving hypoellipticity of solutions of elliptic PDE on products of pcfss fractal sets (~\cite{Strichartz2005TAMS}, Section~11 and~\cite{RogerStricJdM} Section~8).
\end{remark}

The remainder of the section is devoted to the proof of this theorem.  As the working is at times technical it may help the reader to have a concrete example in mind.  The most elementary choice is to let $X$ be the unit interval with Lebesgue measure and the maps $F_{j}$ be the contractions onto the left and right halves of $X$. In this case we obtain the usual Dirichlet energy, the Laplacian is the second derivative, and the normal derivative is the outward-directed first derivative at the endpoints~\cite{Strichartz2006}.  Eigenfunctions and piecewise eigenfunctions are exponentials, and a quick computation shows that $\eta^{(\lambda)}_{1}=\sinh \sqrt{\lambda}x/\sinh\sqrt{\lambda}$.  Our goal is to show that the basic features of this function, as listed in the theorem, are also true for the functions $\eta^{(\lambda)}_{p}$ on any pcfss set.

An observation we use throughout this section is~\eqref{celldecompofeta}, which says that we may decompose these functions into pieces corresponding to cells of a partition $\Theta_{k}$ and obtain a linear combination of pieces; moreover if $\Theta_{k}$ is suitably chosen then each piece looks almost like a copy of a fixed function.  In the case of the unit interval the pieces are literally translates and reflections of a fixed hyperbolic sine, and the decomposition~\eqref{celldecompofeta} can be obtained by iteration of the double angle formula (see Section~2 of~\cite{IoPeRoHuSt2010TAMS}).  The main idea of this section is that in the decomposition each piece has a larger normal derivative near its peak than at the other boundary points, and in order for such pieces to join together in a smooth fashion it is necessary that the peak of each is smaller than that of its neighbor by a constant factor.  Thus the exponential decay may be derived from the shape of a collection of basic functions and counting the cells in paths on $\Theta_{k}$.  In essence, our estimates on the relative sizes of the normal derivatives come from the fact that $\eta^{(\lambda)}_{p}$ is subharmonic, but there is some unavoidable technical work to obtain the correct quantative bounds.  It should be emphasized that most of the work in this section depends on the fact that $\lambda$ is positive.

\begin{assumption}
In this section we require that $\lambda>0$.
\end{assumption}
We will extensively use the decomposition of $\eta_{p}^{(\lambda)}$ as a piecewise $\lambda$-eigenfunction on the
cells of a partition $\Theta_{k}$ as defined in Section~\ref{partitionssection} above.  Using the Laplacian scaling
\begin{equation*}
    \Delta \bigl( \eta_{p}^{(r_{w}\mu_{w}\lambda)}\circ F_{w}^{-1} \bigr)
    = (r_{w}\mu_{w})^{-1}  \bigl( \Delta \eta_{p}^{(r_{w}\mu_{w}\lambda)} \bigr) \circ F_{w}^{-1}
    = \lambda \eta_{p}^{(r_{w}\mu_{w}\lambda)} \circ F_{w}^{-1}
    \end{equation*}
we see that $\eta_{p}^{(\lambda)}$ may be written as a linear combination of the functions
$\eta_{y}^{(r_{\theta}\mu_{\theta}\lambda)}$ for $\theta\in\Theta_{k}$ and $y\in V(\Theta_{k})$.
More precisely,
\begin{equation}\label{celldecompofeta}
    \eta_{p}^{(\lambda)}(x)
    = \sum_{\theta\in\Theta_{k}} \sum_{q\in V_{0}} \eta_{p}^{(\lambda)}(F_{\theta}q)\,
        \bigl( \eta_{q}^{r_{\theta}\mu_{\theta}\lambda}\circ F_{\theta}^{-1}(x)\bigr).
    \end{equation}

Before proceeding we require several preparatory lemmas about the function $\eta^{(\lambda)}_{p}$.
The first is a maximum principle for smooth subharmonic functions. It is well known but
does not seem to appear in the literature, except for that part which is in Proposition~2.11 of \cite{NeeStrTeYu2004JFA}.

\begin{lemma}\label{localminimpliesDeltaunonneg}
Suppose $u\in\dom(\Delta)$ and $\Delta u\geq 0$.  If $u$ attains its global maximum at an interior
point then $u$ is constant.   Moreover $\Delta u\leq0$ at any local maximum point.
\end{lemma}

\begin{proof}
The proof uses the fact that on any cell $F_{w}(X)$ there is a Green kernel $g_{w}$ that is
non-negative on $F_{w}(X)$ and strictly positive on $F_{w}(X)\setminus F_{w}(V_{0})$, and such that
\begin{equation}\label{uintermsofhandgreens}
    u(x)
    = h_{w}(x) - \int_{F_{w}(X)} g_{w}(x,y)\bigl(\Delta u(y)\bigr) \, d\mu
    \end{equation}
where $h_{w}(x)$ is the harmonic function on $F_{w}(X)$ with $h_{w}(x)=u(x)$ for all $x\in
F_{w}(V_{0})$.

Let $h_{m}$ be the piecewise harmonic function at scale $m$ with $h_{m}(x)=u(x)$ for all $x\in
V_{m}$.  Since $\Delta u\geq 0$ we find from~\eqref{uintermsofhandgreens} that $u\leq h_{m}$. This
implies that $\partial_{n}u(x)\geq\partial_{n}h_{m}(x)$ for all $x\in V_{m}$. However the sum of
the normal derivatives of $u$ at any point of $V_{m}\setminus V_{0}$ must vanish because
$u\in\dom(\Delta)$, so  the sum of the normal derivatives of $h_{m}$ must be non-positive. This
gives that the $m$-scale graph Laplacian of $h_{m}$ is non-negative, so if $h_{m}$ achieves its
maximum at a point $x\in V_{m}\setminus V_{0}$ then this maximum is also attained at all neighbors
of $x$ in the $m$-scale graph.  Connectivity then implies $h_{m}$ is constant on $V_{m}$.

If $u$ attains its global maximum at a point $x\in X\setminus V_{0}$ then $u\leq h_{m}$ and $h_{m}$
harmonic implies that $h_{m}$ achieves the same value at a point on the boundary of any cell
containing $x$. For all sufficiently large $m$, this point (which could be $x$) is not in $V_{0}$,
so $h_{m}$ attains its maximum value $u(x)$ at a point of $V_{m}\setminus V_{0}$, and is therefore
constant by our previous reasoning. Applying this for all large $m$ implies $u=u(x)$ on the dense
set $V_{\ast}$, and since $u$ is continuous it must be constant.

Suppose in order to obtain a contradiction that $u$ has a local maximum at $x\in X\setminus V_{0}$
and $\Delta u(x)>0$.  The easy case is when $x\in V_{\ast}$, because  it is then easily seen that
$h_{m}$ has a local maximum at $x$ for all sufficiently large $m$.  This contradicts the above
reasoning showing the graph Laplacian of $h_{m}$ to be non-negative if there is a neighborhood on
which $\Delta u\geq0$. The alternative is that $x\not\in V_{\ast}$.  Then there is an infinite word
$\w\in\Sigma$ so $\cap_{m} F_{[\w]_{m}}(X)=\{x\}$ and these sets form a neighborhood base of $x$.
From~\eqref{uintermsofhandgreens} we see that $h_{[\w]_{m}}(x)>u(x)$ for all $m$ large enough that
$\Delta u>0$ on $F_{[\w]_{m}}(X)$. The maximum principle for harmonic functions then gives
\begin{equation*}
    \max\bigl\{ u(y) : y\in F_{w}(V_{0})\bigr\}
    = \max \bigl\{h_{[\w]_{m}}(y):y\in F_{w}(V_{0}) \bigr\}
    > u(x),
    \end{equation*}
so that every neighborhood of $x$ contains a point at which $u$ exceeds $u(x)$, in contradiction to
$u(x)$ being a local maximum.
\end{proof}

\begin{definition}
Let $\zeta_{p}$ be the function harmonic on $X$ with $\zeta_{p}(q)=0$ for $q\in V_{0}$, $q\neq p$ and $\zeta_{p}(p)=1$.
\end{definition}

\begin{lemma}\label{etapositivesubharmonic}
$0\leq \eta^{(\lambda)}_{p}\leq \zeta_{p}$
\end{lemma}
\begin{proof}
If the first inequality fails then the fact that $\eta^{(\lambda)}_{p}\geq0$ on $V_{0}$ implies
that $\eta^{(\lambda)}_{p}$ has a strictly negative minimum at some interior point $x$. However by
Lemma~\ref{localminimpliesDeltaunonneg} the Laplacian must be non-negative at a minimum point, in
contradiction to $\Delta\eta^{(\lambda)}_{p}(x)=\lambda\eta^{(\lambda)}_{p}(x)<0$.  Having
established the first inequality, it follows that $\eta^{(\lambda)}_{p}$ is subharmonic and thus is
bounded above by the harmonic function with the same boundary values, which is precisely
$\zeta_{p}$.
\end{proof}

\begin{corollary}\label{maxprinciplecorol}
On any connected open set in $X$ the maximum of $\eta_{p}^{(\lambda)}$ is attained at the boundary.
\end{corollary}
\begin{proof}
We have $\Delta \eta^{(\lambda)}_{p}=\lambda\eta^{(\lambda)}_{p}\geq0$ and so may apply
Lemma~\ref{localminimpliesDeltaunonneg}.
\end{proof}

\begin{lemma}\label{etasorderedbylambda}
If $\lambda>\lambda'$ then $\eta^{(\lambda)}_{p}\leq\eta^{(\lambda')}_{p}$.
\end{lemma}
\begin{proof}
Supposing the contrary we find that $\eta^{(\lambda')}_{p}-\eta^{(\lambda)}_{p}$ has a negative
local minimum at a point $x$.  However using that  $\eta^{(\lambda')}_{p}\geq0$ and the hypothesis,
\begin{equation*}
    \Delta\bigl(\eta^{(\lambda')}_{p}-\eta^{(\lambda)}_{p} \bigr)(x)
    = (\lambda'\eta^{(\lambda')}_{p}-\lambda\eta^{(\lambda)}_{p} \bigr)(x)
    \leq \lambda (\eta^{(\lambda')}_{p}-\eta^{(\lambda)}_{p} \bigr)(x)
    <0
    \end{equation*}
which cannot occur at a minimum point by Lemma~\ref{localminimpliesDeltaunonneg}.
\end{proof}

\begin{corollary}\label{signandorderingofnormalderivs}
For any $p,q\in V_{0}$ with $p\neq q$, the values $\partial_{n}\eta_{p}^{(\lambda)}(p)$ are
non-negative and increasing in $\lambda$, while the values $\partial_{n}\eta_{p}^{(\lambda)}(q)$ are
non-positive and increasing in $\lambda$.  In particular for $p\neq q$ the values $\partial_{n}\eta_{p}^{(\lambda)}(q)$ are bounded below by $\min_{p,q\in V_{0}}\partial_{n}\zeta_{p}(q)<0$.
\end{corollary}
\begin{proof}
Positivity of $\partial_{n}\eta_{p}^{(\lambda)}(p)$ is evident from
Lemma~\ref{etapositivesubharmonic} because $\eta_{p}^{(\lambda)}\leq \zeta_{p}<1$ on
$X\setminus\{p\}$.  Then the fact that $\eta^{(\lambda)}_{p}\leq\eta^{(\lambda')}_{p}$ when
$\lambda>\lambda'$ from Lemma~\ref{etasorderedbylambda} implies
$\partial_{n}\eta_{p}^{(\lambda)}(p)\geq\partial_{n}\eta_{p}^{(\lambda')}(p)$.  Similarly
$\partial_{n}\eta_{p}^{(\lambda)}(q)$ is negative because
$\eta_{p}^{(\lambda)}\geq0=\eta_{p}^{(\lambda)}(q)$ from Lemma~\ref{etapositivesubharmonic}, and
$\eta^{(\lambda)}_{p}\leq\eta^{(\lambda')}_{p}$ implies
$\partial_{n}\eta_{p}^{(\lambda)}(q)\geq\partial_{n}\eta_{p}^{(\lambda')}(q)$.   The lower bound
comes from Lemma~\ref{etapositivesubharmonic}, because $\eta_{p}^{(\lambda)}\leq\zeta_{p}$ for all
$\lambda>0$ implies $\partial_{n}\eta_{p}^{(\lambda)}(q)\geq\partial_{n}\zeta_{p}(q)$ for all
$p,q\in V_{0}$.
\end{proof}

With these basic observations in hand we look more closely at the
decomposition~\eqref{celldecompofeta}.

\begin{lemma}\label{estimateofetaintegral}
For all $p\in V_{0}$,
\begin{equation*}
    \int_{X} \eta_{p}^{(\lambda)}\, d\mu \gtrsim (1+\lambda)^{-\frac{S}{S+1}}.
    \end{equation*}
\end{lemma}
\begin{proof}
For $\lambda\leq1$ the result is clear from the positivity of $\eta_{p}^{(\lambda)}$ and the
monotonicity shown in Lemma~\ref{etasorderedbylambda}; in fact we have $\simeq$ rather than an
inequality. For $\lambda\geq 1$, fix $k$ such that $e^{-k(S+1)}\lambda\in (e^{-1},1]$, so that for
$\theta\in\Theta_{k}$ we have $r_{\theta}\mu_{\theta}\lambda\simeq1$. The monotonicity of
Lemma~\ref{etasorderedbylambda} then implies
\begin{equation*}
    \int_{X} \eta^{(r_{\theta}\mu_{\theta}\lambda)} \, d\mu
    \simeq 1.
    \end{equation*}
All terms in the decomposition~\eqref{celldecompofeta} are positive by
Lemma~\ref{etapositivesubharmonic}, so writing $\theta_{p}$ for a word in $\Theta_{k}$ such that
$F_{\theta_{p}}(X)\ni p$, we have
\begin{align*}
    \int_{X} \eta_{p}^{(\lambda)}\, d\mu
    &\geq \int_{F_{\theta_{p}}(X)} \eta_{p}^{(r_{\theta_{p}}\mu_{\theta_{p}}\lambda)}\circ
    F_{\theta_{p}}^{-1}\, d\mu \\
    &= \mu_{\theta_{p}}\int_{X} \eta_{p}^{(r_{\theta_{p}}\mu_{\theta_{p}}\lambda)}\, d\mu.
    \end{align*}
However $\mu_{\theta_{p}}\simeq \lambda^{-S/(S+1)}$ for $\theta\in\Theta_{k}$ and our choice of
$k$, so the proof is complete.
\end{proof}

\begin{lemma}\label{estimateofsumofnormalderivsoncell}
For any word $w\in W_{\ast}$,
\begin{equation*}
\sum_{x\in F_{w}(V_{0})} \partial_{n}\bigl( \eta^{(r_{w}\mu_{w}\lambda)}_{q}\circ
        F_{w}^{-1} \bigr)(x)
        \gtrsim \frac{\lambda}{\bigl((r_{w}\mu_{w})^{-1}+\lambda \bigr)^{\frac{S}{S+1}}}
    \end{equation*}
\end{lemma}
\begin{proof}
Applying the Gauss-Green formula to $\eta^{(r_{w}\mu_{w}\lambda)}_{q}\circ F_{w}^{-1}$ and the
constant function $1$ yields
\begin{align*}
    \sum_{x\in F_{w}(V_{0})} \partial_{n}\bigl( \eta^{(r_{w}\mu_{w}\lambda)}_{q}\circ
        F_{w}^{-1} \bigr)(x)
    &= \int_{F_{w}(X)} \Delta \bigl( \eta^{(r_{w}\mu_{w}\lambda)}_{q}\circ
        F_{w}^{-1} \bigr) \, d\mu\\
    &= \lambda \int_{F_{w}(X)} \eta^{(r_{w}\mu_{w}\lambda)}_{q} \circ
        F_{w}^{-1}  \, d\mu\\
    &= \lambda \mu_{w} \int_{X} \eta^{(r_{w}\mu_{w}\lambda)}_{q}\, d\mu\\
    &\gtrsim \lambda \mu_{w} \bigl( 1+ r_{w}\mu_{w}\lambda \bigr)^{\frac{-S}{S+1}} \\
    &= \frac{\lambda}{\bigl((r_{w}\mu_{w})^{-1}+\lambda \bigr)^{\frac{S}{S+1}}}
    \end{align*}
where we used Lemma~\ref{estimateofetaintegral} and $\mu_{w}=r_{w}^{S}$.
\end{proof}

On the first reading of the following lemma one should think of the case $\Theta'=\Theta_{k}$.  It
will later be used for $\Theta'=\Theta_{k}\setminus\{\theta:F_{\theta}(X)\ni p\}$ and more
complicated sets.  This lemma is the main argument in this section of the paper, in that it uses smoothness of the join between pieces of the decomposition~\eqref{celldecompofeta} to show that $ \eta_{p}^{(\lambda)}$ must decay.

\begin{lemma}\label{interiorvalueestimatefromboundary}
For $\Theta'\subset\Theta_{k}$ and $Y=\cup_{\theta\in\Theta'} F_{\theta}(X)$ we have
\begin{equation*}
    \sum_{y\in V(\Theta')\setminus \partial Y} \eta_{p}^{(\lambda)}(y)
    \leq C \Bigl( e^{k(S+1)}\lambda^{-1} + e^{k}\lambda^{-1/(S+1)} \Bigr)
        \sum_{z\in\partial Y} \eta_{p}^{(\lambda)}(z).
    \end{equation*}
\end{lemma}
\begin{proof}
At any point $x\in V(\Theta_{k})$ the sum of the normal derivatives of $\eta_{p}^{(\lambda)}$ over
the cells meeting at $x$ must be zero because $\eta_{p}^{(\lambda)}\in\dom(\Delta)$.  If we sum
this cancelation over all $V(\Theta_{k})$ points that are interior to $Y$ we may
use~\eqref{celldecompofeta} and rearrange to obtain
\begin{equation}\label{normalderivsumoverinterior}
    \sum_{\theta\in\Theta'} \sum_{q\in V_{0}} \eta_{p}^{(\lambda)}(F_{\theta}(q))
    \biggl( \sum_{x\in F_{\theta}(V_{0})\setminus\partial Y}
    \partial_{n} \bigl( \eta_{q}^{(r_{\theta}\mu_{\theta}\lambda)} \circ F_{\theta}^{-1} \bigr)(x)
    \biggr)
    =0.
    \end{equation}
We estimate the innermost sum in~\eqref{normalderivsumoverinterior} using
Lemma~\ref{estimateofsumofnormalderivsoncell}.  Suppose first that $\theta$ and $q$ are such that
$F_{\theta}q$ is an interior point of $Y$.  Then Corollary~\ref{signandorderingofnormalderivs}
tells us that the terms $\partial_{n} \bigl( \eta_{q}^{(r_{\theta}\mu_{\theta}\lambda)} \circ
F_{\theta}^{-1} \bigr)(x)$ for $x\in\partial Y$ are negative.  Thus the lower bound of
Lemma~\ref{estimateofsumofnormalderivsoncell} is still valid with these points removed.
Substituting into~\eqref{normalderivsumoverinterior} we find
\begin{align}\label{normalderivsestimatewithbdyandinteriorseparate}
    \lefteqn{ -\sum_{\theta\in\Theta'} \sum_{\{q\in V_{0}:F_{\theta}(q)\in\partial Y\}} \eta_{p}^{(\lambda)}(F_{\theta}(q))
    \biggl( \sum_{x\in F_{\theta}(V_{0})\setminus\partial Y}
    \partial_{n} \bigl( \eta_{q}^{(r_{\theta}\mu_{\theta}\lambda)} \circ F_{\theta}^{-1} \bigr)(x)
    \biggr) }\quad\\
    &= \sum_{\theta\in\Theta'} \sum_{\{q\in V_{0}:F_{\theta}(q)\not\in\partial Y\}} \eta_{p}^{(\lambda)}(F_{\theta}(q))
    \biggl( \sum_{x\in F_{\theta}(V_{0})\setminus\partial Y}
    \partial_{n} \bigl( \eta_{q}^{(r_{\theta}\mu_{\theta}\lambda)} \circ F_{\theta}^{-1} \bigr)(x)
    \biggr) \notag\\
    &\gtrsim \sum_{\theta\in\Theta'} \sum_{\{q\in V_{0}:F_{\theta}(q)\not\in\partial Y\}} \eta_{p}^{(\lambda)}(F_{\theta}(q))
    \frac{\lambda}{\bigl((r_{\theta}\mu_{\theta})^{-1}+\lambda \bigr)^{\frac{S}{S+1}}} \notag\\
    &\gtrsim \frac{\lambda}{\bigl(ce^{k(S+1)}+\lambda \bigr)^{\frac{S}{S+1}}}
    \sum_{y\in V(\Theta')\setminus \partial Y} \eta_{p}^{(\lambda)}(y). \notag
    \end{align}
Now in the terms
\begin{equation*}
    \partial_{n} \bigl( \eta_{q}^{(r_{\theta}\mu_{\theta}\lambda)} \circ
        F_{\theta}^{-1} \bigr)(x)
    = r_{\theta}^{-1} \partial_{n}\eta_{q}^{(r_{\theta}\mu_{\theta}\lambda)} \bigl( F_{\theta}^{-1}(x)
    \bigr)
    \end{equation*}
on the left in~\eqref{normalderivsestimatewithbdyandinteriorseparate}, the points $q$ and
$F_{\theta}^{-1}(x)$ cannot coincide, because $x\not\in\partial Y$ and $F_{\theta}q\in\partial Y$.
Again appealing to Corollary~\ref{signandorderingofnormalderivs} we see that such normal
derivatives are negative, increasing in $\lambda$, and bounded below by $-c=\min_{p,q\in
V_{0}}\partial_{n}\zeta_{p}(q)$, which depends only on the harmonic structure of $X$. Putting this
and the estimate $r_{\theta}\simeq e^{-k}$
into~\eqref{normalderivsestimatewithbdyandinteriorseparate} gives a positive constant $C$ depending
on $-c$ and the degree of vertices so that
\begin{align*}
    Ce^{k} \sum_{z\in\partial Y} \eta_{p}^{(\lambda)}(z)
    &\geq -\sum_{\theta\in\Theta'} \sum_{\{q\in V_{0}:F_{\theta}(q)\in\partial Y\}} \eta_{p}^{(\lambda)}(F_{\theta}(q))
    \sum_{x\in F_{\theta}(V_{0})\setminus\partial Y} (-c e^{k})\\
    &\gtrsim \frac{\lambda}{\bigl(ce^{k(S+1)}+\lambda \bigr)^{\frac{S}{S+1}}} \sum_{y\in  V(\Theta')\setminus \partial Y} \eta_{p}^{(\lambda)}(y).
    \end{align*}
We also use that $\bigl(ce^{k(S+1)}+\lambda \bigr)^{\frac{S}{S+1}}\lesssim e^{kS}+
\lambda^{\frac{S}{S+1}}$.
\end{proof}

We may obtain decay estimates of $\eta_{p}^{(\lambda)}$ by iterative use of
Lemma~\ref{interiorvalueestimatefromboundary} for an appropriate choice of $k$ and sets $\Theta'$.

\begin{definition}\label{defnofklambda}
Given $\lambda>0$ we let $k(\lambda)$ be the larger of $0$ and the greatest integer such that
\begin{equation*}
    k(\lambda)\leq \frac{1}{S+1}\log \lambda -\log C - 2
    \end{equation*}
where $C$ is the constant in Lemma~\ref{interiorvalueestimatefromboundary}.  Note that
$e^{k(\lambda)}\simeq (1+\lambda)^{1/(S+1)}$.
\end{definition}

The following result is an immediate consequence of the definition of $k(\lambda)$.
\begin{corollary}\label{simplerinteriorvalueestimatefromboundary}
If $\lambda$ is large enough that $k(\lambda)\geq1$, then for $\Theta'\subset\Theta_{k(\lambda)}$
\begin{equation*}
    \sum_{y\in V(\Theta')\setminus \partial Y} \eta_{p}^{(\lambda)}(y)
    \leq \frac{1}{e} \sum_{z\in\partial Y} \eta_{p}^{(\lambda)}(z).
    \end{equation*}
\end{corollary}

\begin{lemma}\label{expdecayalongpaths}
For $\lambda$ as in Corollary~\ref{simplerinteriorvalueestimatefromboundary} and each $p\in V_{0}$, let $X_{k(\lambda)}(p,i)= \bigl\{ x\in V(\Theta_{k(\lambda)}): d_{k(\lambda)}(p,x)\geq i \bigr\}$. Then
\begin{equation*}
    \sum_{y\in X_{k(\lambda)}(p,i)} \eta_{p}^{(\lambda)}(y)
    \leq e^{-i}
    \end{equation*}
\end{lemma}
\begin{proof}
We induct over $i$.  The base case $i=0$ is simply the fact that $\eta_{p}^{(\lambda)}(p)=1$.
Observe that points $y\in \partial X_{k(\lambda)}(p,i)$ for $i\geq2$ satisfy $d^{(\lambda)}(p,y)\geq i-1$ or
are in $V_{0}\setminus\{p\}$. The latter may be ignored because $\eta_{p}^{(\lambda)}$ is zero on
$V_{0}\setminus\{p\}$, so substituting the inductive estimate into
Lemma~\ref{simplerinteriorvalueestimatefromboundary} gives the result.
\end{proof}

If we restrict $\eta_{p}^{(\lambda)}$ to a cell, Corollary~\ref{maxprinciplecorol} implies the maximum is at the boundary. If a point $x\in X$ has $d_{k(\lambda)}(p,x)=i$ then the boundary points of the cell of $\Theta_{k(\lambda)}$ that contains $x$ are in $X_{k(\lambda)}(p,i)$  Thus  Lemma~\ref{expdecayalongpaths} implies
\begin{corollary}\label{etaexpdecayinchemicalmetric}
For any $x\in X$ and $\lambda$ as in Corollary~\ref{simplerinteriorvalueestimatefromboundary},
\begin{equation*}
    \eta_{p}^{(\lambda)}(x)
    \leq \exp \Bigl( - d_{k(\lambda)}(p,x)  \Bigr).
    \end{equation*}
\end{corollary}

Lemma~\ref{expdecayalongpaths} also allows us to show that Lemma~\ref{estimateofetaintegral} gives the correct value
for the integral as $\lambda\rightarrow\infty$.

\begin{corollary}\label{asymptoticestimateofetaintegral}
For all $p\in V_{0}$,
\begin{equation*}
    \int_{X} \eta_{p}^{(\lambda)}\, d\mu \simeq (1+\lambda)^{-\frac{S}{S+1}}.
    \end{equation*}
\end{corollary}
\begin{proof}
For $\lambda\leq 1$ the result was observed in Lemma~\ref{estimateofetaintegral}. If $\lambda>1$
then for $\theta\in\Theta_{k(\lambda)}$ we reason as in the proof of
Lemma~\ref{estimateofetaintegral} to find
\begin{equation*}
    \int_{F_{\theta}(X)} \eta_{q}^{(r_{\theta}\mu_{\theta}\lambda)}\, d\mu
    \simeq \mu_{\theta}
    \simeq \lambda^{-\frac{S}{S+1}}.
    \end{equation*}
Integrating the decomposition~\eqref{celldecompofeta} then yields
\begin{align*}
    \int_{X} \eta_{p}^{(\lambda)}\, d\mu
    &\simeq \lambda^{-\frac{S}{S+1}} \sum_{\theta\in\Theta(k(\lambda))} \sum_{q\in V_{0}}
    \eta_{p}^{(\lambda)}(F_{\theta}q)\\
    &=\lambda^{-\frac{S}{S+1}} \sum_{i} \sum_{x\in \in X_{k(\lambda)}(p,i)}
    \eta_{p}^{(\lambda)}(x)\\
    &\leq \lambda^{-\frac{S}{S+1}} \sum_{i} e^{-i}
    \end{align*}
where the last inequality is from Lemma~\ref{expdecayalongpaths}.  The reverse inequality is
Lemma~\ref{estimateofetaintegral}.
\end{proof}

Using the result of Corollary~\ref{asymptoticestimateofetaintegral} in the proof of
Lemma~\ref{estimateofsumofnormalderivsoncell} improves it as well.
\begin{corollary}\label{asymptoticestimateofsumofnormalderivsoncell}
For any word $w\in W_{\ast}$,
\begin{equation*}
\sum_{x\in F_{w}(V_{0})} \partial_{n}\bigl( \eta^{(r_{w}\mu_{w}\lambda)}_{q}\circ
        F_{w}^{-1} \bigr)(x)
        \simeq \frac{\lambda}{\bigl((r_{w}\mu_{w})^{-1}+\lambda \bigr)^{\frac{S}{S+1}}}
    \end{equation*}
\end{corollary}

From here it is not difficult to obtain estimates of the normal derivatives of
$\eta_{p}^{(\lambda)}$ at points of $V_{0}$.  Recall from Section~\ref{partitionssection} that
$d_{k}(p,q)$ is $\gtrsim e^{k}$, so $\exp \Bigl(- d_{k(\lambda)}(p,q) \Bigr)$ tends to be small
when $\lambda$ is large.  The following lemma therefore tells us that the normal derivatives
$\partial_{n}\eta_{p}^{(\lambda)}(q)$ are much smaller when $q\neq p$ than when $q= p$.

\begin{lemma}\label{normalderivativeestimate}
For $p,q\in V_{0}$,
\begin{gather*}
    \partial_{n}\eta_{p}^{(\lambda)}(p) \simeq (1+\lambda)^{\frac{1}{S+1}}, \\
    - (1+\lambda)^{\frac{1}{S+1}} \exp \Bigl( - d_{k(\lambda)}(p,q) \Bigr) \lesssim \partial_{n}\eta_{p}^{(\lambda)}(q)\leq 0.
    \end{gather*}
\end{lemma}
\begin{proof}
Suppose $q\neq p$. Let $\theta^{1},\dotsc,\theta^{n}$ be those words from $\Theta_{k(\lambda)}$ for
which $F_{\theta^{i}}(V_{0})\ni q$. Using~\eqref{celldecompofeta} and the lower bound for the
normal derivatives $\partial_{n}\bigl( \eta_{x}^{(r_{\theta^{i}}\mu_{\theta^{i}}\lambda)}
\bigr)(y)$ for $x,y\in V_{0}$ with $y\neq x$ from Corollary~\ref{signandorderingofnormalderivs}, as
well as $r_{\theta^{i}}^{-1}\simeq e^{k(\lambda)}\simeq (1+\lambda)^{1/(S+1)}$,
\begin{align*}
    0\geq\partial_{n}\eta_{p}^{(\lambda)}(q)
    &= \sum_{i=1}^{n} \sum_{x\in V_{0}} \eta_{p}^{(\lambda)} (F_{\theta^{i}}x)
        \partial_{n}\bigl( \eta_{x}^{(r_{\theta^{i}}\mu_{\theta^{i}}\lambda)}\circ F_{\theta^{i}}^{-1}\bigr)
        (q)\\
    &=\sum_{i=1}^{n} \sum_{x\in V_{0}} \eta_{p}^{(\lambda)} (F_{\theta^{i}}x)
        r_{\theta^{i}}^{-1}\partial_{n}\bigl( \eta_{x}^{(r_{\theta^{i}}\mu_{\theta^{i}}\lambda)} \bigr)
        \bigl(F_{\theta^{i}}^{-1}(q)\bigr)\\
    &\gtrsim (1+\lambda)^{\frac{1}{S+1}} \Bigl( \max_{x,y\in V_{0}, x\neq y} -\partial_{n} \zeta_{x}(y) \Bigr)
        \sum_{i=1}^{n} \sum_{x\in V_{0}} \eta_{p}^{(\lambda)} (F_{\theta^{i}}x)\\
    &\geq -C (1+\lambda)^{\frac{1}{S+1}} \sum_{z\in X_{k(\lambda)}(q,1)} \eta_{p}^{(\lambda)}(z)\\
    &\gtrsim - (1+\lambda)^{\frac{1}{S+1}} \exp \Bigl(1 - d_{k(\lambda)}(p,q) \Bigr)
    \end{align*}
where in the final step we applied Lemma~\ref{expdecayalongpaths} and used that points in $z\in
X_{k(\lambda)}(q,1)$ have $d_{k(\lambda)}(p,x)\geq d_{k(\lambda)}(p,q)-1$.  This gives the desired
estimate for $q\neq p$.

As a special case of the above estimate we have
$0\geq\partial_{n}\eta_{p}^{(\lambda)}(q)\gtrsim- \lambda^{\frac{1}{S+1}}$.  Combining this with
Corollary~\ref{asymptoticestimateofsumofnormalderivsoncell} for the empty word $w$ gives the
desired result for $\partial_{n}\eta_{p}^{(\lambda)}(p)$ when $\lambda\geq 1$.  For $\lambda\leq1$
the conclusion is clear from Corollary~\ref{signandorderingofnormalderivs}, because
$\partial_{n}\eta_{p}^{(\lambda)}(p)$ is increasing, so is bounded below by
$\partial_{n}\zeta_{p}(p)$ and above by $\partial_{n}\eta_{p}^{(1)}(p)$.
\end{proof}

To complete our picture of the behavior of $\eta^{(\lambda)}$ we need a lower estimate on its decay and a corresponding upper estimate for the normal derivative $\partial_{n}\eta^{(\lambda)}_{p}(q)$ when $p\neq q$.  After a preliminary lemma, these may be obtained by somewhat simpler reasoning than that used earlier.
\begin{lemma}\label{upperboundfornormalderivofnormalizedeta}\
\begin{enumerate}
    \item The function $\eta^{(\lambda)}_{p}$ is non-zero all points in $V_{\ast}\setminus V_{0}$,
    \item For $p\neq q$ we have $\partial_{n}\eta^{(\lambda)}_{p}(q)<0$,
    \item For $w\in \Theta_{k(\lambda)}$ and $p\neq q$ we have $\partial_{n}\eta^{(r_{w}\mu_{w}\lambda)}_{p}(q)\lesssim -1$.
    \end{enumerate}
\end{lemma}
\begin{proof}
The first step is to prove a weaker version of the second statement.  From the Gauss-Green formula we have
\begin{align*}
    \partial_{n}\eta^{(\lambda)}_{p}(q)-\partial_{n}\zeta_{p}(q)
    &=\partial_{n}\eta^{(\lambda)}_{p}(q)-\partial_{n}\zeta_{q}(p)\\
    &=\sum_{x\in V_{0}} \partial_{n}\eta^{(\lambda)}_{p}(x)\zeta_{q}(x)-\eta^{(\lambda)}_{p}(x)\partial_{n}\zeta_{q}(x)\\
    &=\int_{X} \bigl(\Delta\eta^{(\lambda)}_{p}\bigr)\zeta_{q}\, d\mu\\
    &=\lambda \int_{X} \eta^{(\lambda)}_{p}\zeta_{q}\, d\mu
    \end{align*}
so that $\partial_{n}\eta^{(\lambda)}_{p}(q)\to\partial_{n}\zeta_{p}(q)$ as $\lambda\downarrow0$. With $p\neq q$ we have $\partial_{n}\zeta_{p}(q)<0$, hence we may find some $\tilde{c}<0$ and $\tilde{\lambda}$ such that $\partial_{n}\eta^{(\lambda)}_{p}(q)\leq\tilde{c}$ for all $\lambda\leq\tilde{\lambda}$.

Now at $x\in V_{\ast}\setminus V_{0}$ the fact that $\eta^{(\lambda)}_{p}\in\dom(\Delta)$ requires that the normal derivatives sum to zero. Take $k$ so large that $x\in V(\Theta_{k})$ and for all $\theta\in\Theta_{k}$ we have $r_{\theta}\mu_{\theta}\lambda_{0}\leq\tilde{\lambda}$.  Suppose that $\eta^{(\lambda)}_{p}(x)=0$. Using~\eqref{celldecompofeta} we can write this sum of normal derivatives as
\begin{align}\label{intermedstepforupperboundofnormalderivative}
    0
    &=\sum_{\{\theta\in\Theta_{k}:F_{\theta}(V_{0})\ni x\}} \sum_{p'\in V_{0}} \eta_{p}^{(\lambda)}(F_{\theta}p')\,
        \partial_{n} \bigl( \eta_{p'}^{(r_{\theta}\mu_{\theta}\lambda)}\circ F_{\theta}^{-1}(x)\bigr) \notag\\
    &= \sum_{\{\theta\in\Theta_{k}:F_{\theta}(V_{0})\ni x\}} \sum_{\{p'\in V_{0}:F_{\theta}(p')\neq x\}} \eta_{p}^{(\lambda)}(F_{\theta}p')\,
        r^{-1}_{\theta} \bigl( \partial_{n} \eta_{p'}^{(r_{\theta}\mu_{\theta}\lambda)}\bigr) (F_{\theta}^{-1}(x)) \notag\\
    &\leq \tilde{c} \sum_{\{\theta\in\Theta_{k}:F_{\theta}(V_{0})\ni x\}}\sum_{\{p'\in V_{0}:F_{\theta}(p')\neq x\}}
        r^{-1}_{\theta} \eta_{p}^{(\lambda)}(F_{\theta}p')
    \end{align}
where the first step uses the scaling of the normal derivative and the fact that $\eta^{(\lambda)}_{p}(x)=0$, while the second uses that $\bigl( \partial_{n} \eta_{p'}^{r_{\theta}\mu_{\theta}\lambda}\bigr) (F_{\theta}^{-1}(x))\leq\tilde{c}$  because $r_{\theta}\mu_{\theta}\lambda\leq\tilde{\lambda}$ and $p'\neq F_{\theta}^{-1}(x)$.  However $\tilde{c}<0$ and all $r^{-1}_{\theta} \eta_{p}^{(\lambda)}(F_{\theta}p')\geq0$, so the only way~\eqref{intermedstepforupperboundofnormalderivative} can be true is if these values are all zero, meaning that $\eta_{p}^{(\lambda)}$ vanishes at the boundary points $F_{\theta}(V_{0})$ of each cell meeting at $x$.  Repeating the argument inductively we see after finitely many steps that $\eta_{p}^{(\lambda)}$ must vanish at all points in $V(\Theta_{k})$, which is impossible because $p\in V(\Theta_{k})$ and $\eta_{p}^{(\lambda)}(p)=1$.  This proves the first statement of the lemma.

The second assertion of the lemma now follows fairly easily from the first.  Using the partition $\Theta_{k}$ as above and the decomposition~\eqref{intermedstepforupperboundofnormalderivative} for the normal derivative at $x=q$ we see
\begin{align}\label{anotherintermedstepinprovingupperboundfornormalderivofeta}
    \partial_{n}\eta_{p}^{(\lambda_{0})}(q)
    &= \sum_{\{\theta\in\Theta_{k}:F_{\theta}(V_{0})\ni q\}} \sum_{\{p'\in V_{0}:F_{\theta}(p')\neq q\}} \eta_{p}^{(\lambda)}(F_{\theta}p')\,
        r^{-1}_{\theta} \bigl( \partial_{n} \eta_{p'}^{(r_{\theta}\mu_{\theta}\lambda)}\bigr) (F_{\theta}^{-1}(q)) \notag\\
    &\leq \tilde{c} \sum_{\{\theta\in\Theta_{k}:F_{\theta}(V_{0})\ni q\}} \sum_{p'\in V_{0}}
        r^{-1}_{\theta} \eta_{p}^{(\lambda)}(F_{\theta}p')
    \end{align}
and we have already shown that the values $r^{-1}_{\theta} \eta_{p}^{(\lambda)}(F_{\theta}p')>0$ when $F_{\theta}p'\not\in V_{0}$.

Finally the third statement follows from the second because $\partial_{n}\eta^{(\lambda)}_{p}(q)$ is continuous, non-positive and increasing in $\lambda$, and $r_{w}\mu_{w}\lambda$ is bounded above when $w\in \Theta_{k(\lambda)}$, with all of these depending only on the fractal and harmonic structure.
\end{proof}

\begin{lemma}
There is a constant $a>0$ such that for any point $x\in V(\Theta_{k(\lambda)})\setminus V_{0}$,
\begin{equation*}
    \eta^{(\lambda)}_{p}(x)
    \geq a \max\Bigl\{ \eta^{(\lambda)}_{p}(y): \exists\theta\in\Theta_{k}\text{ with }F_{\theta}(V_{0})\supset \{x,y\}\Bigr\}.
    \end{equation*}
These points $y$ are the neighbors of $x$ in $V(\Theta_{k(\lambda)})$.
\end{lemma}
\begin{proof}
Smoothness of $\eta^{(\lambda)}_{p}$ requires that the first step of~\eqref{intermedstepforupperboundofnormalderivative} holds, where we take $k=k(\lambda)$.  As we are not assuming $\eta^{(\lambda)}_{p}(x)=0$ we obtain instead of the second step of~\eqref{intermedstepforupperboundofnormalderivative}
\begin{align*}
    \lefteqn{ \eta_{p}^{(\lambda)}(x) \sum_{\{\theta\in\Theta_{k(\lambda)}:F_{\theta}(V_{0})\ni x\}} r_{\theta}^{-1} \partial_{n} \bigl( \eta_{F_{\theta}^{-1}(x)}^{(r_{\theta}\mu_{\theta}\lambda)} \bigr)( F_{\theta}^{-1}(x)) } \quad & \\
    &= - \sum_{\{\theta\in\Theta_{k(\lambda)}:F_{\theta}(V_{0})\ni x\}} \sum_{\{q\in V_{0}:F_{\theta}(q)\neq x\}} \eta_{p}^{(\lambda)}(F_{\theta}(q))
    r_{\theta}^{-1} \partial_{n} \bigl( \eta_{q}^{(r_{\theta}\mu_{\theta}\lambda)} \bigr)( F_{\theta}^{-1}(x))\\
    &\gtrsim \sum_{\{\theta\in\Theta_{k(\lambda)}:F_{\theta}(V_{0})\ni x\}} \sum_{\{q\in V_{0}:F_{\theta}(q)\neq x\}} r_{\theta}^{-1} \eta_{p}^{(\lambda)}(F_{\theta}(q))
    \end{align*}
where we used that $- \partial_{n} \bigl( \eta_{q}^{(r_{\theta}\mu_{\theta}\lambda)} \bigr)( F_{\theta}^{-1}(x))\gtrsim 1$ from the third part of Lemma~\ref{upperboundfornormalderivofnormalizedeta}. On the left of this inequality each of the normal derivatives $\partial_{n} \bigl( \eta_{F_{\theta}^{-1}(x)}^{(r_{\theta}\mu_{\theta}\lambda)} \bigr)( F_{\theta}^{-1}(x))$ is bounded above by a constant depending only on the harmonic structure, as shown in the first part of Lemma~\ref{normalderivativeestimate}.  Since the values $r_{\theta}^{-1}$ are comparable on both sides of the equation, the result follows.
\end{proof}

\begin{corollary}\label{lowerdecayestimateforeta}
Let $X_{k(\lambda),p}$ be the subset of $X$ obtained by deleting those cells $F_{\theta}(X)$, $\theta\in \Theta_{k(\lambda)}$ that intersect $V_{0}$ at points other than $p$.  There is $c>0$ such that for all $x\in X_{k(\lambda),p}$
\begin{equation*}
    \eta^{(\lambda)}_{p}(x)\gtrsim \exp\bigl(-c d_{k(\lambda)}(p,x) \bigr)
    \end{equation*}
\end{corollary}
\begin{proof}
If $x\in V(\Theta_{k(\lambda)})\setminus V_{0}$ and $d_{k(\lambda)}(p,x)=i$ then there is some $y\in V(\Theta_{k(\lambda)})$ with $d_{k(\lambda)}(x,y)=i-1$. By the previous result, $\eta^{(\lambda)}_{p}(x)\geq a \eta^{(\lambda)}_{p}(y)$ and by induction $\eta^{(\lambda)}_{p}(x)\geq a^{i} \eta^{(\lambda)}_{p}(p)=a^{i}$.  The result for points of $V(\Theta_{k(\lambda)})\setminus V_{0}$ follows by setting $c=-\log a$, and we note that the estimate $\eta^{(\lambda)}_{p}(x)\leq C\exp\bigl(- d_{k(\lambda)}(p,x) \bigr)$ from Corollary~\ref{etaexpdecayinchemicalmetric} ensures $c>0$.  To obtain the bound for a general point $x\in X_{k(\lambda),p}$, let $\theta\in\Theta_{k(\lambda)}$ be such that $F_{\theta}(X)\ni x$.  The value $d_{k(\lambda)}(p,x)$ is the distance from $p$ to the nearest point of $F_{\theta}(V_{0})$; all other points of this form are at most $d_{k(\lambda)}(p,x)+1$ from $p$, and none is in $V_{0}\setminus\{p\}$, so the result for $V(\Theta_{k(\lambda)})\setminus V_{0}$ applies to them.  It therefore suffices to know that the restriction of $\eta^{(\lambda)}_{p}$ to $F_{\theta}(X)$ is bounded below by a multiple of its boundary values.  Observe that  this function is a linear combination of the functions $\eta^{(r_{\theta}\mu_{\theta}\lambda)}_{q}\circ F_{\theta}^{-1}$ with the boundary values as coefficients.  Hence it is enough to know that each $\eta^{(r_{\theta}\mu_{\theta}\lambda)}_{q}(y)$ has a positive lower bound on those cells $F_{j}(X)$ that do not contain a point of $V_{0}\setminus\{q\}$.  This follows from the third part of Lemma~\ref{upperboundfornormalderivofnormalizedeta}, continuity of $\eta_{p}^{(\lambda)}$ and the fact that the values $r_{\theta}\mu_{\theta}\lambda$ lie in a bounded interval, and we see that the constant depends only on the harmonic structure of the fractal.
\end{proof}

\begin{corollary}\label{lowerdecayboundcorolforeta} For the constant $c$ of Corollary~\ref{lowerdecayestimateforeta}
\begin{equation*}
    \partial_{n}\eta_{p}^{(\lambda_{0})}(q)
    \lesssim -(1+\lambda)^{1/(S+1)} \exp\bigl(-c d_{k(\lambda)}(p,x) \bigr)
    \end{equation*}
\end{corollary}
\begin{proof}
If we rewrite~\eqref{anotherintermedstepinprovingupperboundfornormalderivofeta} with $k=k(\lambda)$ we have instead of $\tilde{c}$ a constant $\lesssim-1$ as determined in the third part of Lemma~\ref{upperboundfornormalderivofnormalizedeta}, so that
\begin{equation*}
    \partial_{n}\eta_{p}^{(\lambda_{0})}(q)
    \lesssim - \sum_{\{\theta\in\Theta_{k}:F_{\theta}(V_{0})\ni q\}} \sum_{\{p'\in V_{0}:F_{\theta}(p')\neq q\}}
        r^{-1}_{\theta} \eta_{p}^{(\lambda)}(F_{\theta}p')
    \end{equation*}
however on the right the values $r^{-1}_{\theta}\simeq(1+\lambda)^{1/(S+1)}$ and there is at least one point $F_{\theta}p'$ from  $V(\Theta_{k(\lambda)})\setminus V_{0}$ for which $d_{k(\lambda)}(p,F_{\theta}p')=d_{k(\lambda)}(p,q)-1$, so that $\eta_{p}^{(\lambda)}(F_{\theta}p')\geq C\exp\bigl(-c d_{k(\lambda)+c}(p,x) \bigr)$ as seen in Corollary~\ref{lowerdecayestimateforeta}.  Since all terms on the right have the same sign, this term gives an upper bound.
\end{proof}

%%%%%%%%%%%%%%%%%%%%%%%%%%%%%%%%%%%%%%%%%%%%%%%%%%%%%%%%%%%%%%%%%%%%%%%%%%%%%%%%%%%%%%%%%%%%%%%%%%%%%%%%%%%%%%%%%%%%%%%%%%%%%%%%%%%%%%%%%%

\section{Estimates of the resolvent kernel on the positive real axis}\label{resolventestimatesection}

In this section of the paper we use the estimates of Section~\ref{estimatesofpiecewiseefnssection} and the series expression~\eqref{maintheoremfromIPRRS} for the resolvent $G^{(\lambda)}(x,y)$ of the Laplacian to obtain estimates in the case $\lambda\in(0,\infty)$. The main result is as follows.

\begin{theorem}\label{pathdecayforGonXthm}
There are constants $\kappa_{1}$ and $\kappa_{2}$ depending only on the fractal, harmonic structure and measure, and such that
if $\lambda>0$,
\begin{equation}\label{pathdecayforGonXthmestimateofG}
    (1+\lambda)^{-1/(S+1)} \exp\Bigl( - \kappa_{1} d_{k(\lambda)}\bigl( x,y \bigr) \Bigr)
    \lesssim G^{(\lambda)}(x,y)
    \lesssim (1+\lambda)^{-1/(S+1)} \exp\Bigl( - \kappa_{2} d_{k(\lambda)}\bigl( x,y \bigr) \Bigr)
    \end{equation}
except that the lower bound is not valid if $x$ or $y$ is in a cell $F_{\theta}(X)$, $\theta\in\Theta_{k(\lambda)}$ such that $F_{\theta}(V_{0})$ contains a point of $V_{0}\setminus\{p\}$.  For these latter cells, the appropriate estimate is instead one on the normal derivative.  Specifically, if $p\in V_{0}$, $y$ is not in one of the above cells $F_{\theta}(X)$ and the normal derivative $\partial_{n}'$ is taken with respect to the first variable, then
\begin{equation}\label{pathdecayforGonXthmestimateofpartialnG}
    \exp\Bigl( - \kappa_{1} d_{k(\lambda)}\bigl( p,y \bigr) \Bigr)
    \lesssim -\partial_{n}'  G^{(\lambda)}(p,y)
    \lesssim \exp\Bigl( - \kappa_{2} d_{k(\lambda)}\bigl( p,y \bigr) \Bigr).
    \end{equation}
A symmetrical result holds for the normal derivative $-\partial_{n}''$ with respect to the second variable.  There are also bounds appropriate to the case where both points are near the boundary.  If $p$ and $q$ are points in $V_{0}$ then
\begin{equation}\label{pathdecayforGonXthmestimateofpartialnpartialnnGoffdiag}
    (1+\lambda)^{1/(S+1)} \exp\Bigl( - \kappa_{1} d_{k(\lambda)}\bigl( p,q \bigr) \Bigr)
    \lesssim \partial_{n}''\partial_{n}'  G^{(\lambda)}(p,q)
    \lesssim (1+\lambda)^{1/(S+1)} \exp\Bigl( - \kappa_{2} d_{k(\lambda)}\bigl( p,q \bigr) \Bigr).
    \end{equation}
It is easy to see these estimates are equivalent to the following global bounds.  Let $R(V_{0},x)$ denote the resistance distance from $x$ to $V_{0}$.  Then
\begin{equation*}
    \exp\Bigl( - \kappa_{1} d_{k(\lambda)}\bigl( p,y \bigr) \Bigr)
    \lesssim  \lambda^{-1/(S+1)} \Bigl( R(x,V_{0})^{-1}+\lambda^{1/(S+1)} \Bigr) \Bigl( R(y,V_{0})^{-1}+\lambda^{1/(S+1)} \Bigr) G^{(\lambda)}(x,y)
    \lesssim \exp\Bigl( - \kappa_{2} d_{k(\lambda)}\bigl( p,y \bigr) \Bigr).
\end{equation*}
\end{theorem}

With this in hand will be relatively easy to obtain a similar result for the Neumann resolvent $G_{N}^{(\lambda)}(x,y)$, which has vanishing normal derivatives rather than zero values at points of $V_{0}$.
\begin{corollary}\label{pathdecayforGNeumannonXthm}
For $\lambda>0$ the Neumann resolvent satisfies the upper and lower bounds of~\eqref{pathdecayforGonXthmestimateofG} everywhere on $X$.
\end{corollary}

The proofs of these results occupy the rest of this section.  We begin by expanding the expression~\eqref{maintheoremfromIPRRS} for the resolvent of the Laplacian and using~\eqref{PsiintermsofGpq} to obtain
\begin{equation}\label{resolventserieswithpsiterms}
    G^{(\lambda)}(x,y)
    = \sum_{w\in W_{\ast}}  \sum_{p,q\in V_{1}\setminus V_{0}} r_{w}
    G^{(r_{w}\mu_{w}\lambda)}_{pq}\:
    \psi^{(r_{w}\mu_{w}\lambda)}_{p}(F_{w}^{-1}x)\: \psi^{(r_{w}\mu_{w}\lambda)}_{q}(F_{w}^{-1}y),
    \end{equation}
where we recall that $G^{(\lambda)}_{pq}$ is the inverse of the matrix $B_{pq}^{(\lambda)}$ defined
in~\eqref{defnofBpq}.  In the next few lemmas we apply the estimates from Section~\ref{estimatesofpiecewiseefnssection} to the terms in this series, for which purpose we require the following assumption.

\begin{assumption}
For the remainder of this section we require that $\lambda>0$.
\end{assumption}

\begin{lemma}\label{Bpqclosetodiagonal}
Let $D^{(\lambda)}$ be the diagonal matrix with entries $D_{pp}^{(\lambda)}= (1+\lambda)^{1/(S+1)}\bigl(B_{pp}^{(\lambda)}\bigr)^{-1}$ and
$E^{(\lambda)}=\bigl(D^{(\lambda)}\bigr)^{-1}-(1+\lambda)^{-1/(S+1)}B^{(\lambda)}$.  Then for all $\lambda$ and $p$ we have $D^{(\lambda)}_{pp}\simeq1$,  $E^{(\lambda)}_{pp}=0$, and for $q\neq p$
\begin{equation*}
    \exp \Bigl( - c_{1}c d_{k(\lambda)}(p,q)  \Bigr)
        \lesssim E^{(\lambda)}_{pq}
        \lesssim \exp \Bigl( - c_{2} d_{k(\lambda)}(p,q)  \Bigr),
    \end{equation*}
where $c_{1}$ and $c_{2}$ depend only on the fractal and its harmonic structure and $c$ is the constant in Theorem~\ref{maintheoremoneta}.
\end{lemma}
\begin{proof}
Comparing~\eqref{defnofBpq} and~\eqref{definitionofpsi} we have
\begin{align*}
    B_{pq}^{(\lambda)}
    &=\sum_{j:F_{j}(V_{0})\ni q} \partial_{n}^{F_{j}(X)} \eta^{(r_{j}\mu_{j}\lambda)}_{F_{j}^{-1}p} \circ
    F_{j}^{-1}(q)\\
    &=\sum_{j:F_{j}(V_{0})\ni q} r_{j}^{-1} \partial_{n} \eta^{(r_{j}\mu_{j}\lambda)}_{F_{j}^{-1}p}
    \bigl( F_{j}^{-1}(q) \bigr).
    \end{align*}
However  Theorem~\ref{maintheoremoneta} then shows that for $p=q$ we have $B_{pp}^{(\lambda)}\simeq
(1+\lambda)^{1/(S+1)}$, so $D^{(\lambda)}_{pp}\simeq1$.  Also from Theorem~\ref{maintheoremoneta} we have for $p\neq q$
\begin{align*}
    -(1+\lambda)^{1/(S+1)} \exp \Bigl( - d_{k(r_{j}\mu_{j}\lambda)}(p,q) \Bigr)
        &\lesssim r_{j}^{-1} \partial_{n} \eta^{(r_{j}\mu_{j}\lambda)}_{F_{j}^{-1}p} \bigl( F_{j}^{-1}(q) \bigr) \\
        &\lesssim -(1+\lambda)^{1/(S+1)} \exp \Bigl( - c d_{k(r_{j}\mu_{j}\lambda)}(p,q) \Bigr),
    \end{align*}
and the observation that $d_{k(r_{j}\mu_{j}\lambda)}(p,q)\simeq d_{k(\lambda)}(p,q)$ lets us choose
appropriate constants $c_{1}$ and $c_{2}$.
\end{proof}

\begin{lemma}\label{Gpqestimate}
For any word $w\in W_{\ast}$ and all $\lambda>0$,
\begin{equation*}
     \exp \Bigl( - c_{1} c d_{k(r_{w}\mu_{w}\lambda)}(p,q)  \Bigr)
    \lesssim  \bigl((r_{w}\mu_{w})^{-1}+\lambda\bigr)^{\frac{1}{S+1}} r_{w} G_{pq}^{(r_{w}\mu_{w}\lambda)}
    \lesssim  \exp \Bigl( - c_{2} d_{k(r_{w}\mu_{w}\lambda)}(p,q)  \Bigr).
\end{equation*}
where $c_{1}$, $c_{2}$ and $c$ are as in Lemma~\ref{Bpqclosetodiagonal}.
\end{lemma}
\begin{proof}
Recall that $d_{k(\lambda)}(p,q)\gtrsim e^{k(\lambda)}\simeq (1+\lambda)^{1/(S+1)}$, so that $\exp \Bigl( - d_{k(\lambda)}(p,q) \Bigr)$ can be made
arbitrarily small by taking $\lambda$ large enough.  Lemma~\ref{Bpqclosetodiagonal} then implies $B^{(\lambda)}$ is close to diagonal, so we may find $G^{(\lambda)}=\bigl(B^{(\lambda)}\bigr)^{-1}$ via the Neumann series.  Specifically,
\begin{align*}
    (1+\lambda)^{1/(S+1)} G^{(\lambda)}
    = \bigl(I-D^{(\lambda)}E^{(\lambda)}\bigr)^{-1} D^{(\lambda)}
    = D^{(\lambda)} + \sum_{k=1}^{\infty} \bigl( D^{(\lambda)}E^{(\lambda)}\bigr)^{k}D^{(\lambda)}.
    \end{align*}
provided $\lambda\geq C_{1}$ where $C_{1}$ is chosen large enough that the series converges.  Observe that this $C_{1}$
depends only on the structure of the fractal, because $D^{(\lambda)}$ contains only values $\simeq1$ and the values in $E^{(\lambda)}$ satisfy the estimate in Lemma~\ref{Bpqclosetodiagonal}.  Notice also that all values in $D^{(\lambda)}$ and $E^{(\lambda)}$ are positive, hence the same is true of all terms in the series.  Making the obvious upper and lower estimates of the sum of the series we conclude that for all $p$ and $q$
\begin{equation*}
    \exp \Bigl( - c_{1}c d_{k(\lambda)}(p,q)  \Bigr)
        \lesssim (1+\lambda)^{1/(S+1)} G^{(\lambda)}_{pq}
        \lesssim \exp \Bigl( - c_{2} d_{k(\lambda)}(p,q)  \Bigr).
    \end{equation*}
Substituting $(r_{w}\mu_{w}\lambda)$ in place of $\lambda$ and using $r_{w}\mu_{w} = r_{w}^{S+1}$ gives
\begin{equation*}
     \exp \Bigl( - c_{1} c d_{k(r_{w}\mu_{w}\lambda)}(p,q)  \Bigr)
    \lesssim  \bigl((r_{w}\mu_{w})^{-1}+\lambda\bigr)^{\frac{1}{S+1}} r_{w} G_{pq}^{(r_{w}\mu_{w}\lambda)}
    \lesssim  \exp \Bigl( - c_{2} d_{k(r_{w}\mu_{w}\lambda)}(p,q)  \Bigr),
\end{equation*}
which proves the estimate for $\lambda\geq C_{1}(r_{w}\mu_{w})^{-1}$.  However if $0\leq\lambda\leq C_{1}(r_{w}\mu_{w})^{-1}$ then it is immediate that $\bigl((r_{w}\mu_{w})^{-1}+\lambda\bigr)^{\frac{1}{S+1}} r_{w}$ is bounded above and below, and we see $G_{pq}^{(r_{w}\mu_{w}\lambda)}$ is bounded above and below by continuity of its dependence on $r_{w}\mu_{w}\lambda$.  At the same time $d_{k(r_{w}\mu_{w}\lambda)}(p,q)$ is bounded because it is the distance on a cellular partition of bounded scale.  All constants depend only on the fractal and harmonic structure, so the result follows.
\end{proof}

\begin{lemma}\label{superexpdecayofpsi}
There are positive constants $c_{3}$ and $c_{4}$ depending only on the harmonic structure, such that if $c$ is the constant from Theorem~\ref{maintheoremoneta} and $p\in V_{1}\setminus V_{0}$, then for $x$ in the support of $\psi^{(r_{w}\mu_{w}\lambda)}_{p}(F_{w}^{-1}x)$, which is the union of all cells $F_{wj}(X)$ that contain $F_{w}(p)$,
\begin{equation}\label{superexpdecayofpsiboundforpsi}
    \exp \Bigl(- c_{3} c d_{k(\lambda)}(F_{w}p,x ) \Bigr)
    \lesssim \psi^{(r_{w}\mu_{w}\lambda)}_{p}(F_{w}^{-1}x)
    \lesssim \exp \Bigl(- c_{4} d_{k(\lambda)}(F_{w}p,x ) \Bigr),
    \end{equation}
except that the lower bound does not hold on the cells $F_{wj}\circ F_{\theta}(X)$ for those $\theta\in\Theta_{k(r_{wj}\mu_{wj}\lambda)}$ such that $F_{\theta}(X)$ intersects $V_{0}\setminus\{p\}$.  For these exceptional cells the correct estimate is that if $x=F_{w}(q)$ for some $q\in V_{0}$ then
\begin{equation}\label{superexpdecayofpsiboundforpartialnpsi}
    \begin{split}
    \exp \Bigl( - c_{3}cd_{k(\lambda)}(F_{w}p,F_{w}q )  \Bigr)
    &\lesssim -\bigl((r_{w}\mu_{w})^{-1}+\lambda\bigr)^{-1/(S+1)} \partial_{n}\bigl( \psi^{(r_{w}\mu_{w}\lambda)}_{p}\circ F_{w}^{-1}\bigr)(x)\\
    &\lesssim  \exp \Bigl( - c_{4} d_{k(\lambda)}(F_{w}p,F_{w}q )  \Bigr).
    \end{split}
    \end{equation}
\end{lemma}
\begin{proof}
From~\eqref{definitionofpsi} we see that
$\psi^{(r_{w}\mu_{w}\lambda)}_{p}(F_{w}^{-1}x)$ is a piecewise $\lambda$-eigenfunction on
$F_{w}(X)$ with value $1$ at $F_{w}(p)$ and zero at the other points of $F_{w}(V_{1})$. It is
non-zero precisely on the the cells $F_{wj}(X)$ such that $p\in F_{j}(V_{0})$.  On each such cell
it is equal to $\eta^{(r_{j}\mu_{j}r_{w}\mu_{w}\lambda)}_{F_{j}^{-1}p}\circ F_{wj}^{-1}$. According
to Theorem~\ref{maintheoremoneta} we have
\begin{equation}\label{subexplowerboundforetaeqn}
    \exp \Bigl(- c d_{k(r_{wj}\mu_{wj}\lambda)}({F_{j}^{-1}p},F_{wj}^{-1}(x) ) \Bigr)
    \lesssim \eta^{(r_{wj}\mu_{wj}\lambda)}_{F_{j}^{-1}p}\bigl( F_{wj}^{-1}(x) \bigr)
    \lesssim \exp \Bigl(- d_{k(r_{wj}\mu_{wj}\lambda)}({F_{j}^{-1}p},F_{wj}^{-1}(x) ) \Bigr),
    \end{equation}
except that the lower bound does not hold on the cells excluded in the statement of the lemma.

From the partition $\Theta_{k(r_{j}\mu_{j}r_{w}\mu_{w}\lambda)}$ of $X$, form the partition
$r_{w}\Theta_{k(r_{j}\mu_{j}r_{w}\mu_{w}\lambda)}$ of $F_{w}(X)$ and observe that the scale is comparable
to that of $\Theta_{k(\lambda)}$ restricted to $F_{w}(X)$.  It follows that there are positive $c_{3}$ and $c_{4}$
depending on the harmonic structure so that for $y,z\in X$,
\begin{equation}\label{distancecomparabilityremark}
    c_{3} d_{k(\lambda)}(F_{wj}y,F_{wj}z)
    \geq d_{k(r_{wj}\mu_{wj}\lambda)}(y,z)
    \geq c_{4} d_{k(\lambda)}(F_{wj}y,F_{wj}z)
    \end{equation}
provided we  round appropriately.
Substituting~\eqref{distancecomparabilityremark} into~\eqref{subexplowerboundforetaeqn} and eliminating the rounding by taking a
suitably large multiple of the exponential gives the desired estimate for  $\psi^{(r_{w}\mu_{w}\lambda)}_{p}(F_{w}^{-1}x)$ because it is a finite sum of such terms.

For the normal derivative we note that if $x=F_{w}q$ for $q\in V_{0}$ then $x=F_{wj}q'$ for some $q'\in V_{0}\setminus\{p\}$.  Using Theorem~\ref{maintheoremoneta} at $q'$ we obtain
\begin{equation}\label{subexplowerboundforetaeqnnormalderiv}
     \begin{split}
     \exp \Bigl( - c d_{k(r_{wj}\mu_{wj}\lambda)}({F_{j}^{-1}p},q' )  \Bigr)
     &\lesssim -(1+r_{wj}\mu_{wj}\lambda)^{-1/(S+1)} \partial_{n}\eta^{(r_{wj}\mu_{wj}\lambda)}_{F_{j}^{-1}p}\bigl( q' \bigr)\\
     &\lesssim \exp \Bigl( -  d_{k(r_{wj}\mu_{wj}\lambda)}({F_{j}^{-1}p}, q' )  \Bigr).
     \end{split}
     \end{equation}
Recalling that
\begin{equation*}
    \partial_{n}\bigl(\eta^{(r_{wj}\mu_{wj}\lambda)}_{F_{j}^{-1}p}\circ F_{wj}^{-1}\bigr)
    = r_{wj}^{-1} \bigl(\partial_{n}\eta^{(r_{wj}\mu_{wj}\lambda)}_{F_{j}^{-1}p}\bigr) \circ F_{wj}^{-1}
    \end{equation*}
we see that
\begin{equation*}
    -\bigl((r_{wj}\mu_{wj})^{-1}+\lambda\bigr)^{-1/(S+1)} \partial_{n}\bigl(\eta^{(r_{wj}\mu_{wj}\lambda)}_{F_{j}^{-1}p}\circ F_{wj}^{-1}\bigr)(x)
    \end{equation*}
also satisfies the estimate in~\eqref{subexplowerboundforetaeqnnormalderiv}.  As $\partial_{n}\psi^{(r_{w}\mu_{w}\lambda)}_{p}\circ F_{w}^{-1}$ is a finite sum of such terms, the bound~\eqref{superexpdecayofpsiboundforpartialnpsi} may be obtained by substituting~\eqref{distancecomparabilityremark} into this estimate.
\end{proof}

%%%%%%%%%%%%%%%%%%%%%%%

We divide the proof of Theorem~\ref{pathdecayforGonXthm} into two parts: the proof of the upper bounds and the proof of the lower bounds.

\begin{proof}[Proof of Theorem~\protect{\ref{pathdecayforGonXthm}}: Upper bounds]
Fix $x,y$ in $X$.  If $x\neq y$ let $w\in W_{\ast}$ be the longest word such that $x,y\in
F_{w}(X)$, and otherwise let $w$ be an infinite word such that $F_{w}(X)=\{x\}=\{y\}$ (in this
latter case the expression~\eqref{Gassumovercontractionsofaword} below may need an additional sum
over the possible choices of $w$, but we suppress this because it does not otherwise affect the
working). Then the series~\eqref{resolventserieswithpsiterms} terminates at scale $|w|$ (which is
$+\infty$ if $x=y$) and may be written
\begin{equation}\label{Gassumovercontractionsofaword}
    G^{(\lambda)}(x,y)
    = \sum_{i=0}^{|w|} \sum_{p,q\in V_{1}\setminus V_{0}} r_{[w]_{i}}
        G^{(r_{[w]_{i}}\mu_{[w]_{i}}\lambda)}_{pq}\:
        \psi^{(r_{[w]_{i}}\mu_{[w]_{i}}\lambda)}_{p}(F_{[w]_{i}}^{-1}x)\:
        \psi^{(r_{[w]_{i}}\mu_{[w]_{i}}\lambda)}_{q}(F_{[w]_{i}}^{-1}y).
    \end{equation}
It will be convenient to divide the sum into three pieces according to the size of $r_{w}\mu_{w}\lambda$.  Define
\begin{equation}\label{defnofinought}
    i_{0}=\begin{cases}
        \min\bigl\{i:r_{[w]_{i}}\mu_{[w]_{i}}\lambda\leq 1\bigr\}
            &\text{ if } r_{w}\mu_{w}\lambda\leq 1 \\
        +\infty&\text{ if } r_{w}\mu_{w}\lambda> 1,
        \end{cases}
    \end{equation}
so that if $i_{0}<\infty$ then $r_{[w]_{i_{0}}}\simeq \lambda^{-1/(S+1)}$. The three pieces of the
sum are as follows, where we note that it is possible for a piece to be empty.
\begin{align*}
    I_{1}
        &= \sum_{i_{0}}^{|w|} \sum_{p,q\in V_{1}\setminus V_{0}} r_{[w]_{i}}
        G^{(r_{[w]_{i}}\mu_{[w]_{i}}\lambda)}_{pq}\:
        \psi^{(r_{[w]_{i}}\mu_{[w]_{i}}\lambda)}_{p}(F_{[w]_{i}}^{-1}x)\:
        \psi^{(r_{[w]_{i}}\mu_{[w]_{i}}\lambda)}_{q}(F_{[w]_{i}}^{-1}y)\\
    I_{2}
        &= \sum_{0}^{\min\{i_{0}-1,|w|\}} \sum_{p\neq q\in V_{1}\setminus V_{0}} r_{[w]_{i}}
        G^{(r_{[w]_{i}}\mu_{[w]_{i}}\lambda)}_{pq}\:
        \psi^{(r_{[w]_{i}}\mu_{[w]_{i}}\lambda)}_{p}(F_{[w]_{i}}^{-1}x)\:
        \psi^{(r_{[w]_{i}}\mu_{[w]_{i}}\lambda)}_{q}(F_{[w]_{i}}^{-1}y)\\
    I_{3}
        &= \sum_{0}^{\min\{i_{0}-1,|w|\}} \sum_{p\in V_{1}\setminus V_{0}} r_{[w]_{i}}
        G^{(r_{[w]_{i}}\mu_{[w]_{i}}\lambda)}_{pp}\:
        \psi^{(r_{[w]_{i}}\mu_{[w]_{i}}\lambda)}_{p}(F_{[w]_{i}}^{-1}x)\:
        \psi^{(r_{[w]_{i}}\mu_{[w]_{i}}\lambda)}_{p}(F_{[w]_{i}}^{-1}y)
    \end{align*}

For both $I_{1}$ and $I_{2}$ we use the trivial estimate
\begin{equation*}
    0\leq  \psi^{(r_{[w]_{i}}\mu_{[w]_{i}}\lambda)}_{p}(F_{[w]_{i}}^{-1}x)\:
        \psi^{(r_{[w]_{i}}\mu_{[w]_{i}}\lambda)}_{q}(F_{[w]_{i}}^{-1}y)
    \leq 1
    \end{equation*}
for the $\psi$ factors.  Then for $I_{1}$ we have $r_{[w]_{i}}\mu_{[w]_{i}}\lambda\lesssim1$, from
which the values $\bigl|G^{(r_{[w]_{i}}\mu_{[w]_{i}}\lambda)}_{pq}\bigr| \lesssim 1$, and so
\begin{gather}\label{estimateforIone}
        0\leq  I_{1} \lesssim \sum_{i_{0}}^{|w|} r_{[w]_{i}} \lesssim r_{[w]_{i_{0}}} \lesssim (1+\lambda)^{-1/(S+1)}
        \text{ if } r_{w}\mu_{w}\lambda\leq 1\\
        I_{1}=0 \text{ if } r_{w}\mu_{w}\lambda> 1 \notag
    \end{gather}
where we used the geometric decay of the $r_{[w]_{i}}$ and the definition of $i_{0}$. Note that we have the upper bound $(1+\lambda)^{-1/(S+1)}$ rather than just $\lambda^{-1/(S+1)}$ because when $\lambda<1$ we have $i_{0}=0$ and $r_{[w]_{i_{0}}}=1$ corresponding to the empty word.

For $I_{2}$ we instead estimate the factors $r_{[w]_{i}}G^{(r_{[w]_{i}}\mu_{[w]_{i}}\lambda)}_{pq}$
using Lemma~\ref{Gpqestimate} to obtain
\begin{align}
    \bigl| I_{2} \bigr|
    &\leq
    \sum_{0}^{\min\{i_{0}-1,|w|\}} \sum_{p\neq q\in V_{1}\setminus V_{0}} r_{[w]_{i}}
        G^{(r_{[w]_{i}}\mu_{[w]_{i}}\lambda)}_{pq} \notag\\
    &\leq  \sum_{0}^{\min\{i_{0}-1,|w|\}} \bigl((r_{[w]_{i}}\mu_{[w]_{i}})^{-1} + \lambda\bigr)^{-1/(S+1)}
        \exp \Bigl( - c_{2} d_{k(r_{[w]_{i}}\mu_{[w]_{i}}\lambda)}(p,q) \Bigr) \notag\\
    &\lesssim \begin{cases}
        (1+\lambda)^{-1/(S+1)} &\text{ if }r_{w}\mu_{w}\lambda\leq 1\\
        (1+\lambda)^{-1/(S+1)} \exp \Bigl( - c_{2}  d_{k(r_{w}\mu_{w}\lambda)}(p,q) \Bigr)
            &\text{ if }r_{w}\mu_{w}\lambda> 1
        \end{cases}\label{estimateforItwo}
    \end{align}
because $d_{k}(p,q)$ is at least exponentially increasing in $k$, so the sum is dominated by a
constant multiple of the largest term, which is the one with the maximal value of $i$.  Moreover
the constant multiple depends only on the $r$ and $\mu$ values and the way path lengths on
$\Theta_{k}$ grow with $k$, all of which are properties only of the fractal, harmonic structure and
measure. In the case that the maximal value of $i$ is $i_{0}$ we have $r_{w}\mu_{w}\lambda\simeq1$
so $d_{k(r_{w}\mu_{w}\lambda)}(p,q)\leq1$ and the exponential term is trivial.  The factors $(1+\lambda)^{-1/(S+1)}$ can be used rather than the more obvious choice $\lambda^{-1/(S+1)}$ because if $\lambda<1$ the sum is empty.

For the $I_{3}$ term we must take a different approach, because in this case $i<i_{0}$ implies we
can apply Lemma~\ref{Gpqestimate} to obtain
\begin{equation*}
    r_{[w]_{i}}G^{(r_{[w]_{i}}\mu_{[w]_{i}}\lambda)}_{pq}
    \simeq \bigl((r_{[w]_{i}}\mu_{[w]_{i}})^{-1}+\lambda\bigr)^{-1/(S+1)} D^{(r_{[w]_{i}}\mu_{[w]_{i}}\lambda)}_{pq}
    \simeq \bigl((r_{[w]_{i}}\mu_{[w]_{i}})^{-1}+\lambda\bigr)^{-1/(S+1)}
    \end{equation*}
and therefore the only estimate we have is that from Lemma~\ref{superexpdecayofpsi}, which gives
\begin{align}
    I_{3}
    &\lesssim (1+\lambda)^{-1/(S+1)} \sum_{i=0}^{\min\{i_{0}-1,|w|\}} \sum_{p\in V_{1}\setminus V_{0}}
        \psi^{(r_{[w]_{i}}\mu_{[w]_{i}}\lambda)}_{p}(F_{[w]_{i}}^{-1}x)\:
        \psi^{(r_{[w]_{i}}\mu_{[w]_{i}}\lambda)}_{p}(F_{[w]_{i}}^{-1}y) \notag\\
    &\lesssim (1+\lambda)^{-1/(S+1)} \sum_{i=0}^{\min\{i_{0}-1,|w|\}} \sum_{p\in V_{1}\setminus V_{0}}
        \exp\Bigl( - c_{4} d_{k(\lambda)}\bigl(F_{[w]_{i}}p,x\bigr) - c_{4}d_{k(\lambda)}\bigl(F_{[w]_{i}}p,y\bigr)
        \Bigr) \label{Ithreeestimate}
    \end{align}
where we again have replaced the obvious $\lambda^{-1/(S+1)}$ factor with $(1+\lambda)^{-1/(S+1)}$ because the sum is empty if $\lambda<1$.

The above indicates that we need a lower bound on $d_{k(\lambda)}\bigl(F_{[w]_{i}}p,x\bigr) + d_{k(\lambda)}\bigl(F_{[w]_{i}}p,y\bigr)$ for the points
$F_{[w]_{i}}p$, $p\in V_{1}\setminus V_{0}$ and a mechanism for counting how many of them there are.  Let $L_{l}$ be the set of such points that are in the boundary of the $l$-cell containing $x$ but not in the boundary of the $(l+1)$-cell containing $x$, then the
number of points in $L_{l}$ is smaller than the number of points in $V_{1}\setminus V_{0}$, which is a constant
depending only on the structure of the fractal.  Moreover any point $z\in L_{l}$ is separated from
$x$ by an $(l+1)$-cell. By the same reasoning as in Lemma~\ref{superexpdecayofpsi}, the
$d_{k(\lambda)}$ diameter of any such $(l+1)$ cell $F_{[w]_{l}j}(X)$ is bounded below by
a constant multiple of the $d_{k(r_{[w]_{l}}\mu_{[w]_{l}}\lambda)}$ diameter of $X$, rounded down to the nearest integer.
Writing the diameter of $X$ with respect to $d_{k}$ as $\diam_{k}(X)$ we obtain:
\begin{equation*}
    \sum_{F_{[w]_{i}}p\in L_{l}} \exp\Bigl( - cd_{k(\lambda)}\bigl(F_{[w]_{i}}p,x\bigr) \Bigr)
    \lesssim \exp\Bigl( - c\diam_{k(r_{[w]_{l}}\mu_{[w]_{l}}\lambda)} (X) \Bigr).
    \end{equation*}
We use this estimate only for $l\leq|w|$.  Any points $F_{[w]_{i}}p$ which occur in the sum but are not in any of the $L_{l}$, $l\leq |w|$ must be in $F_{w}(V_{1})$, so the number of these is bounded by a constant depending only on the structure of the fractal.  For these points we use the triangle inequality to estimate the corresponding terms of~\eqref{Ithreeestimate}.
\begin{equation*}
    \exp\Bigl( - c_{4} d_{k(\lambda)}\bigl(F_{[w]_{i}}p,x\bigr) - c_{4} d_{k(\lambda)}\bigl(F_{[w]_{i}}p,y\bigr) \Bigr)
    \leq \exp\Bigl( - c_{4} d_{k(\lambda)}\bigl(x,y\bigr) \Bigr).
    \end{equation*}
Combining these estimates and substituting into~\eqref{Ithreeestimate} we have shown
\begin{align}
    I_{3}
    &\lesssim (1+\lambda)^{-1/(S+1)} \biggl(
    \exp\Bigl( - c_{4} d_{k(\lambda)}\bigl(x,y\bigr) \Bigr) +
    \sum_{l=0}^{\min\{i_{0},|w|\}} \exp\Bigl( - c\diam_{k(r_{[w]_{l}}\mu_{[w]_{l}}\lambda)} (X)
    \Bigr) \biggr) \notag\\
    &\lesssim \begin{cases}
        (1+\lambda)^{-1/(S+1)} &\text{ if }r_{w}\mu_{w}\lambda\leq 1\\
        (1+\lambda)^{-1/(S+1)} \exp\Bigl( - c_{4} d_{k(\lambda)}\bigl(x,y\bigr) \Bigr)
        &\text{ if }r_{w}\mu_{w}\lambda> 1
        \end{cases}\label{estimateforIthree}
    \end{align}
where we used the same reasoning about the exponential decay as was used for the $I_{2}$ term,
along with the fact that $d_{k(\lambda)}\bigl(x,y\bigr)\lesssim\diam_{k(r_{w}\mu_{w}\lambda)} (X)$
because $x$ and $y$ have resistance separation $\lesssim r_{w}$.

Comparing~\eqref{estimateforIone}, \eqref{estimateforItwo}
and~\eqref{estimateforIthree} we see that we always have a bound by $(1+\lambda)^{-1/(S+1)}$ if
$r_{w}\mu_{w}\lambda\leq 1$, and in this case the resistance separation of $x$ and $y$ is at most a
constant multiple of $r_{w}=(r_{w}\mu_{w})^{1/(S+1)}\lesssim\lambda^{-1/(S+1)}$, so that
$d_{k(\lambda)}\bigl(x,y\bigr)\lesssim 1$.  Setting $\kappa_{2}=\min\{c_{2},c_{4}\}$ we therefore obtain
\begin{equation*}
    G^{(\lambda)}(x,y)
    \lesssim (1+\lambda)^{-1/(S+1)}
    \lesssim (1+\lambda)^{-1/(S+1)} \exp\Bigl( - \kappa_{2} d_{k(\lambda)}\bigl(x,y\bigr) \Bigr).
    \end{equation*}
For the case $r_{w}\mu_{w}\lambda> 1$ we have $I_{1}=0$, and the estimate for $I_{2}$ is dominated
by a multiple of that for $I_{3}$, so again
\begin{equation*}
    G^{(\lambda)}(x,y)
    \lesssim (1+\lambda)^{-1/(S+1)} \exp\Bigl( - \kappa_{2} d_{k(\lambda)}\bigl(x,y\bigr) \Bigr)
    \end{equation*}
and we have established the upper bound of~\eqref{pathdecayforGonXthmestimateofG} stated in the theorem.

We now turn to the upper bounds for normal derivatives.  Recall that $\partial_{n}'$ denotes the normal derivative with respect to the first variable.
If $x\in V_{0}$ and $y\neq x$ then the series for $G^{(\lambda)}(x,y)$ is finite and we can compute the normal derivative term by term. The upper bound on $-\partial_{n}'G^{(\lambda)}(x,y)$ may then be obtained in almost the same way as the upper bound on $G^{(\lambda)}(x,y)$.  The reasoning for both $I_{2}$ and $I_{3}$ is unchanged except that instead of bounding $r_{[w]_{i}}G^{(r_{[w]_{i}}\mu_{[w]_{i}}\lambda)}_{pq}$ by $\lambda^{-1/(S+1)}$ we use the full upper bound $\bigl((r_{[w]_{i}}\mu_{[w]_{i}})^{-1}+\lambda\bigr)^{-1/(S+1)}$ from Lemma~\ref{Gpqestimate}.  The  upper bound for $-\partial_{n}'\psi^{(r_{[w]_{i}}\mu_{[w]_{i}}\lambda)}_{p}\circ F_{[w]_{i}}^{-1}(x)$ from~\eqref{superexpdecayofpsiboundforpartialnpsi} cancels this factor and leaves precisely the exponential decay terms seen in~\eqref{estimateforItwo} and~\eqref{Ithreeestimate}, so the only change to these estimates is that the $\lambda^{-1/(S+1)}$ factor is no longer present.  The $I_{1}$ term requires slightly more changes.  Again we use Lemma~\ref{Gpqestimate} and~\eqref{superexpdecayofpsiboundforpartialnpsi} to see that
\begin{equation}\label{boundfortermwithonenormalderivative}
    r_{[w]_{i}}G^{(r_{[w]_{i}}\mu_{[w]_{i}}\lambda)}_{pq} \Bigl( -\partial_{n}'\psi^{(r_{[w]_{i}}\mu_{[w]_{i}}\lambda)}_{p}\circ F_{[w]_{i}}^{-1}(x)  \Bigr)
    \lesssim \exp \Bigl(- c_{4} d_{k(\lambda)}(F_{[w]_{i}}p, x ) \Bigr)
    \end{equation}
is bounded, but then use it to conclude that
\begin{align*}
    \partial_{n}' I_{1}
    &\lesssim \sum_{i_{0}}^{|w|} \sum_{p,q\in V_{1}\setminus V_{0}}
        \exp \Bigl(- c_{4} d_{k(\lambda)}(F_{[w]_{i}}p, x ) \Bigr)
        \psi^{(r_{[w]_{i}}\mu_{[w]_{i}}\lambda)}_{q}(F_{[w]_{i}}^{-1}y)\\
    &\lesssim  \exp \Bigl(- c_{4} d_{k(\lambda)}(F_{[w]_{i}}p, y ) -c_{4} d_{k(\lambda)}(F_{[w]_{i}}p, x ) \Bigr),
    \end{align*}
where the last step is from~\eqref{superexpdecayofpsiboundforpsi}.  This sum can be bounded by the same argument as was used for $I_{3}$ in passing from~\eqref{Ithreeestimate} to~\eqref{estimateforIthree}, so we may sum the terms from $I_{1}$, $I_{2}$ and $I_{3}$ to obtain the desired bound $-\partial_{n}'G^{(\lambda)}(x,y)\lesssim \exp \Bigl(- \kappa_{2} d_{k(\lambda)}(x,y ) \Bigr)$ provided $y\neq x$.  This verifies the upper bound in~\eqref{pathdecayforGonXthmestimateofpartialnG}.  An easy argument then shows that for points $z$ within the $\Theta_{k(\lambda)}$ cell containing $x$, the peak size of $G^{(\lambda)}(z,\cdot)$ is comparable to $R(x,z)$ rather than $(1+\lambda)^{-1/(S+1)}$; we will later need the immediate consequence \begin{equation}\label{usefulboundaryLoneestforG}
    \|G^{(\lambda)}(z,\cdot) \|_{L^{1}} \lesssim (1+\lambda)^{-1/(S+1)}R(x,z) \text{ as } z\to x\in V_{0}.
    \end{equation}

To complete our upper bounds we must deal with what happens when $x$ and $y$ are both in $V_{0}$.  If this is true and they are distinct, then the sum~\eqref{Gassumovercontractionsofaword} reduces to the finite number of terms with $i=0$.  We can break each term into the product of a factor like that in~\eqref{boundfortermwithonenormalderivative} and a factor for the normal derivative in the second variable, which is of the form
\begin{equation*}
    -\partial_{n}''\bigl(\psi^{(r_{[w]_{i}}\mu_{[w]_{i}}\lambda)}_{q}\circ F_{[w]_{i}}^{-1}\bigr)(y)
    \lesssim \bigl((r_{[w]_{i}}\mu_{[w]_{i}})^{-1}+\lambda\bigr)^{1/(S+1)} \exp \Bigl(- c_{4} d_{k(\lambda)}(F_{[w]_{i}}p, y )\Bigr)
    \end{equation*}
by~\eqref{superexpdecayofpsiboundforpartialnpsi}, but with $i=0$.  However, the points $F_{[w]_{0}}p=p$ are then in $V_{1}\setminus V_{0}$, so both $d_{k(\lambda)}(F_{[w]_{i}}p, y )$ and $d_{k(\lambda)}(F_{[w]_{0}}p, x)$ are comparable to $d_{k(\lambda)}(x,y)$, and as $i=0$ the overall bound reduces to
\begin{equation*}
    \partial_{n}'\partial_{n}'' G^{(\lambda)}(x,y)
    \lesssim \bigl(1+\lambda\bigr)^{1/(S+1)}\exp \Bigl(- c_{4} d_{k(\lambda)}(x, y) \Bigr).
    \end{equation*}
Finally, we look at $\partial_{n}''\partial_{n}'G^{(\lambda)}(x,x)$ for $x\in V_{0}$.  Notice that in the above discussion of the normal derivative $\partial_{n}'G^{(\lambda)}(x,y)$ we established that (summing over those $w$ corresponding to $x$ if necessary)
\begin{equation}\label{expressionforpartialnGlambdaassumofpartialnterms}
    -\partial_{n}'G^{(\lambda)}(x,y)
    = \sum_{i} \sum_{p,q\in V_{1}\setminus V_{0}} r_{[w]_{i}}
        G^{(r_{[w]_{i}}\mu_{[w]_{i}}\lambda)}_{pq}\:
        \Bigl(-\partial_{n}'\psi^{(r_{[w]_{i}}\mu_{[w]_{i}}\lambda)}_{p}\circ F_{[w]_{i}}^{-1}(x)\Bigr)\:
        \psi^{(r_{[w]_{i}}\mu_{[w]_{i}}\lambda)}_{q}(F_{[w]_{i}}^{-1}y)
    \end{equation}
in which the sum is finite for any $y\neq x$.  However replacing $\psi^{(r_{[w]_{i}}\mu_{[w]_{i}}\lambda)}_{q}(F_{[w]_{i}}^{-1}y)$ with $-\partial_{n}''\psi^{(r_{[w]_{i}}\mu_{[w]_{i}}\lambda)}_{q}(F_{[w]_{i}}^{-1}y)$ throughout the series does not give a series for $\partial_{n}''\partial_{n}'G^{(\lambda)}(x,y)$.  The reason is that every finite truncation of the series is zero at $x$, whereas $-\partial_{n}'G^{(\lambda)}(x,x)$ is the continuous extension of $-\partial_{n}'G^{(\lambda)}(x,\cdot)$ from $X\setminus\{x\}$ to $x$, which exists because the series is a sum of continuous positive functions in which all partial sums are $\lesssim \exp \Bigl(- c_{4} d_{k(\lambda)}(x,y ) \Bigr)$ for $y\neq x$.  We therefore compute $\partial_{n}''\partial_{n}'G^{(\lambda)}(x,x)$ using
\begin{align}
    \partial_{n}''\partial_{n}'G^{(\lambda)}(x,x)
    &= \lim_{i\to\infty} E\bigl(-\partial_{n}'G^{(\lambda)}(x,\cdot),h_{i}^{x}(\cdot)\bigr) \notag\\
    &= \lim_{i\to\infty} \sum_{y\in F_{[w]_{i}}(V_{0})} \partial_{n}h_{i}^{x}(y) \bigl(-\partial_{n}'G^{(\lambda)}(x,y)\bigr) \notag\\
    &= \lim_{i\to\infty} \sum_{y\in F_{[w]_{i}}(V_{0})\setminus\{x\}} \partial_{n}h_{i}^{x}(y) \bigl(\partial_{n}'G^{(\lambda)}(x,x) -\partial_{n}'G^{(\lambda)}(x,y)\bigr) \label{computingthepartialonepartialtwoderivofresolventlinethree}\\
    &=\lim_{i\to\infty} \lim_{z\to x} \sum_{y\in F_{[w]_{i}}(V_{0})\setminus\{x\}} \partial_{n}h_{i}^{x}(y) \bigl(\partial_{n}'G^{(\lambda)}(x,z) -\partial_{n}'G^{(\lambda)}(x,y)\bigr) \label{computingthepartialonepartialtwoderivofresolventlinefour}
    \end{align}
where $h_{i}^{x}$ is harmonic on $F_{[w]_{i}}(X)$ and equal to $1$ at $x$ and zero at the other points of $F_{[w]_{i}}(V_{0})$.  We used that $\sum_{y\in F_{[w]_{i}}(V_{0})} \partial_{n}h_{i}^{x}(y)=0$ because $h_{i}^{x}$ is harmonic, and that $\partial_{n}'G^{(\lambda)}(x,x)=\lim_{z\to x}\partial_{n}'G^{(\lambda)}(x,z)$ by definition.  Now rewrite these terms using the Gauss-Green formula in the first variable. Since $(\lambda-\Delta')G^{(\lambda)}(s,y)=\delta_{y}(s)$ we have
\begin{equation*}
    \partial_{n}'G^{(\lambda)}(x,y)
    = \int_{X} \zeta_{x}(s)\Delta'G^{(\lambda)}(s,y)\, d\mu(s)
    = \lambda \int_{X} \zeta_{x}(s) G^{(\lambda)}(s,y)\, d\mu(s) - \zeta_{x}(y)
    \end{equation*}
from which~\eqref{computingthepartialonepartialtwoderivofresolventlinefour} may be written
\begin{align}
    \lefteqn{\partial_{n}''\partial_{n}'G^{(\lambda)}(x,x)}\quad& \notag\\
    &= \lim_{i\to\infty} \sum_{y\in F_{[w]_{i}}(V_{0})\setminus\{x\}} \partial_{n}h_{i}^{x}(y) \bigl(\zeta_{x}(y) - \lim_{z\to x} \zeta_{x}(z)\bigr) \notag\\
    & \quad - \lambda \lim_{i\to\infty} \int_{X}  \zeta_{x}(s)  \sum_{y\in F_{[w]_{i}}(V_{0})\setminus\{x\}} \partial_{n}h_{i}^{x}(y) G^{(\lambda)}(s,y)\, d\mu(s) \notag \\
    &\quad + \lambda \lim_{i\to\infty} \sum_{y\in F_{[w]_{i}}(V_{0})\setminus\{x\}} \partial_{n}h_{i}^{x}(y) \lim_{z\to x} \int_{X} \zeta_{x}(s) G^{(\lambda)}(s,z)\, d\mu(s). \notag\\
    &= \partial_{n}\zeta_{x}(x)  - \lambda \lim_{i\to\infty} \int_{X}  \zeta_{x}(s)  \sum_{y\in F_{[w]_{i}}(V_{0})\setminus\{x\}} \partial_{n}h_{i}^{x}(y) G^{(\lambda)}(s,y)\, d\mu(s) \label{computingthepartialonepartialtwoderivofresolventestimatetwo}
    \end{align}
In the last step of this computation we used that $\lim_{z\to x} \zeta_{x}(z)=\zeta_{x}(x)$, so the first term just produces $\partial_{n}\zeta_{x}(x)$ by the same argument that gave~\eqref{computingthepartialonepartialtwoderivofresolventlinethree}.  We also found that the third term was zero, because of~\eqref{usefulboundaryLoneestforG} and the boundedness of $h_{i}^{x}$.

To further simplify~\eqref{computingthepartialonepartialtwoderivofresolventestimatetwo} we again apply the Gauss-Green formula, this time in the second variable, to obtain
\begin{align*}
    \sum_{y\in F_{[w]_{i}}(V_{0})\setminus\{x\}} \partial_{n}h_{i}^{x}(y) G^{(\lambda)}(s,y)
    &= \partial_{n}'' G^{(\lambda)}(s,x) - \int_{F_{[w]_{j}}(X)} h_{i}^{x}(t) \Delta''G^{(\lambda)}(s,t) \, d\mu(t)\\
    &= \partial_{n}'' G^{(\lambda)}(s,x) + h_{i}^{x}(s) - \lambda \int_{F_{[w]_{j}}(X)} h_{i}^{x}(t) G^{(\lambda)}(s,t) \, d\mu(t)
    \end{align*}
Two pieces of this make no contribution in the limit as $i\to\infty$. One is the integral
\begin{equation*}
    \lambda \int_{X}  \zeta_{x}(s) h_{i}^{x}(s) \, d\mu(s)
    \lesssim  \mu\bigl(F_{[w]_{i}(X)}\bigr) \to 0 \text{ as }i\to\infty.
    \end{equation*}
For the other we use the upper bound of~\eqref{pathdecayforGonXthmestimateofG} to compute $\int_{X} G^{(\lambda)}(s,t) \, d\mu(s)\lesssim \lambda^{-1}$ and therefore
\begin{equation*}
    \lambda^{2} \int_{X} \zeta_{x}(s) \int_{F_{[w]_{j}}(X)} h_{i}^{x}(t) G^{(\lambda)}(s,t) \, d\mu(t)\, d\mu(s)
    \lesssim \lambda \int_{F_{[w]_{j}}(X)} h_{i}^{x}(t) \, d\mu(t)
    \leq \lambda \mu\bigl(F_{[w]_{i}(X)}\bigr)
    \end{equation*}
which also goes to $0$ as $i\to\infty$.  We used that $h_{i}^{x}$ is bounded by $1$.
On the remaining piece we use the previously established upper bound in~\eqref{pathdecayforGonXthmestimateofpartialnG} to see that that $-\partial_{n}'' G^{(\lambda)}(x,y)$ has faster than exponential decay with characteristic length $k(\lambda)$, so that its integral against a bounded function is estimated by the integral over a cell of diameter $k(\lambda)$ at $x$, which has measure comparable to $\lambda^{-S/(S+1)}$.  Thus
\begin{equation}\label{upperboundforterminvolvingpartialnofresolventagainstzeta}
    \lambda \int_{X}  \zeta_{x}(s) \bigl(-\partial_{n}'' G^{(\lambda)}(s,x)\bigr) \, d\mu(s)
    \lesssim  \lambda^{1/(S+1)}.
    \end{equation}
Substituting all of this into~\eqref{computingthepartialonepartialtwoderivofresolventestimatetwo} we have at last
\begin{equation}\label{usefulexpressionfortwopartialderivsofresolvent}
    \partial_{n}''\partial_{n}'G^{(\lambda)}(x,x)
    = \partial_{n}\zeta_{x}(x) +\lambda \int_{X}  \zeta_{x}(s) \bigl(-\partial_{n}'' G^{(\lambda)}(s,x)\bigr) \, d\mu(s)
    \lesssim (1+\lambda)^{1/(S+1)}.
    \end{equation}
This completes all the upper bounds in the theorem.
\end{proof}

%%%%%%%%%%%%%%%%%%%%%%%%%%%%%%%%%%%%%%%%%%%%%%%%%%%%%%%%%%%%%%%%%%%%%%%%

\begin{proof}[Proof of Theorem~\protect{\ref{pathdecayforGonXthm}}: Lower bounds]
Observe that all terms in the sum~\eqref{Gassumovercontractionsofaword} are positive, so to obtain the lower bound of~\eqref{pathdecayforGonXthmestimateofG} it suffices that we have the desired bound on a single term. Our upper estimates suggest we consider the dominant term for both the diagonal and off-diagonal series. To do so we will have to split into two different ranges of $\lambda$.  This splitting  depends on the lower estimate from~\eqref{superexpdecayofpsiboundforpsi} of Lemma~\ref{superexpdecayofpsi}, which is only valid off some exceptional cells.

The exceptional cells in the sense of Lemma~\ref{superexpdecayofpsi} are of the form $F_{[w]_{i}j}\circ F_{\theta}(X)$ for those $\theta\in\Theta_{k(r_{[w]_{i}j}\mu_{[w]_{i}j}\lambda)}$ such that $F_{\theta}(X)$ intersects $V_{0}\setminus\{p\}$.  Observe that if $r_{[w]_{i}j}\mu_{[w]_{i}j}\lambda$ is too small then these cells could cover $F_{[w]_{i}}(X)$, so to use this estimate at all we must assume that $r_{[w]_{i}}\mu_{[w]_{i}}\lambda>\tilde{c}$, where $\tilde{c}$ is chosen large enough this does not occur, and is a constant depending only on the fractal and harmonic structure.

Fix two points $x$ and $y$ and recall $w$ is the longest word so $F_{w}(X)$ contains both $x$ and $y$.  Our two cases are that where $r_{w}\mu_{w}\lambda>\tilde{c}$ and where $r_{w}\mu_{w}\lambda\leq\tilde{c}$, the latter including the case $x=y$.  In the former we prove the lower bound~\eqref{pathdecayforGonXthmestimateofG} including the exponential term, but in the latter the exponential term is bounded below and plays no role in the estimate.

First consider the case $r_{w}\mu_{w}\lambda>\tilde{c}$.  Then the lower estimate from~\eqref{superexpdecayofpsiboundforpsi} of Lemma~\ref{superexpdecayofpsi} is valid off the exceptional cells of $F_{w}(X)$, and these cells have size comparable to cells of the partition $\Theta_{k(\lambda)}$. Fix such a size and call the subcells of $F_{w}(X)$ that both have this size and contain a point of $F_{w}(V_{0})$ the boundary cells of $F_{w}(X)$.  We further arrange this size so that the exceptional cells for $\psi^{r_{w}\mu_{w}\lambda}_{F_{w}p}\circ F_{w}^{-1}$ are contained in those boundary cells that do not contain $F_{w}p$.

From Lemma~\ref{Gpqestimate} and Lemma~\ref{superexpdecayofpsi} we know that for any $i$ such that $r_{[w]_{i}}\mu_{[w]_{i}}\lambda>\tilde{c}$ and $x$ and $y$ not in the exceptional subcells for $\psi^{r_{[w]_{i}}\mu_{[w]_{i}}\lambda}_{F_{[w]_{i}p}}\circ F_{[w]_{i}}^{-1}$
\begin{align}
    \lefteqn{r_{[w]_{i}} G_{pq}^{(r_{[w]_{i}}\mu_{[w]_{i}}\lambda)} \psi^{(r_{[w]_{i}}\mu_{[w]_{i}}\lambda)}_{p}(F_{[w]_{i}}^{-1}x)\: \psi^{(r_{[w]_{i}}\mu_{[w]_{i}}\lambda)}_{q}(F_{[w]_{i}}^{-1}y)}\quad& \notag\\
    &\gtrsim \bigl((r_{[w]_{i}}\mu_{[w]_{i}}^{-1}+\lambda\bigr)^{-1/(S+1)} \exp \Bigl(- c_{1}d_{k(r_{[w]_{i}}\mu_{[w]_{i}}\lambda)}(x,y) - c_{3} c d_{k(\lambda)}(F_{[w]_{i}}p,x) - c_{3} c d_{k(\lambda)}(F_{[w]_{i}}q,y ) \Bigr) \notag\\
    &\geq \lambda^{-1/(S+1)} \exp \Bigl(- c_{1}c_{3} c d_{k(\lambda)}(F_{[w]_{i}}p,F_{[w]_{i}}q) - c_{3} c d_{k(\lambda)}(F_{[w]_{i}}p,x) - c_{3} c d_{k(\lambda)}(F_{[w]_{i}}q,y ) \Bigr) \label{lowerestimateoftermsinGseries}
    \end{align}
where the $d_{k(r_{[w]_{i}}\mu_{[w]_{i}}\lambda)}(x,y)$ term from Lemma~\ref{Gpqestimate} was converted to $c_{3}cd_{k(\lambda)}(F_{[w]_{i}}p,F_{[w]_{i}}q)$ using~\eqref{distancecomparabilityremark}.
To get our lower bound~\eqref{pathdecayforGonXthmestimateofG} on $G^{(\lambda)}(x,y)$ for $\lambda\gtrsim1$ it suffices that given $x$ and $y$, neither of which is in an exceptional cell containing a point of $V_{0}$,  we can find such $i$, $p$ and $q$  with $x$ and $y$  not in exceptional subcells of $F_{[w]_{i}}(X)$ and
\begin{equation}\label{desiredupperboundondistancefunctionstogetGlowerbd}
    d_{k(\lambda)}(F_{[w]_{i}}p,F_{[w]_{i}}q) + d_{k(\lambda)}(F_{[w]_{i}}p,x)+ d_{k(\lambda)}(F_{[w]_{i}}q,y )
    \leq d_{k(\lambda)}(x,y) + k'
    \end{equation}
for some constant $k'$ depending only on the fractal and harmonic structure, as the corresponding term from~\eqref{lowerestimateoftermsinGseries} will be a lower bound of the correct form for the non-negative series~\eqref{Gassumovercontractionsofaword} for $G^{(\lambda)}(x,y)$ once we set $\kappa_{1}=c_{1}c_{3}c$.

One possibility is that  neither $x$ nor $y$ is in a boundary cell of $F_{w}(X)$.  Our choice of $w$ ensures that $x$ and $y$ are in separate level-one subcells of $F_{w}(X)$ (i.e. cells of the form $F_{wj}(X)$), so the $d_{k(\lambda)}$ geodesic between them must pass through a vertex $F_{w}p$ for some $p\in V_{1}\setminus V_{0}$.  We conclude that $d_{k(\lambda)}(F_{w}p,x )+d_{k(\lambda)}(F_{w}p,y )=d_{k(\lambda)}(x,y)$, which is the desired bound~\eqref{desiredupperboundondistancefunctionstogetGlowerbd} in the case $q=p$ and $i=|w|$.

The alternative is that either $x$ or $y$ or both is in a boundary subcell of $F_{w}(X)$.  There is no loss of generality in assuming $x$ is in a boundary subcell, and the corresponding boundary point is $F_{w}p'$ for some $p'\in V_{0}$.  We will assume that $F_{w}p'\not\in V_{0}$, in which case there is a largest $i$ satisfying $F_{w}p'=F_{[w]_{i}}p$ for some $p\in V_{1}\setminus V_{0}$.  If $y$ is also in a boundary subcell then by taking $x$ to be the point with the larger value of $i$ of this type, we can assume that the boundary point corresponding to $y$ is of the form $F_{[w]_{i}}q$ for some $q\in V_{\ast}\setminus V_{0}$.

Since the boundary subcells have size comparable to cells from $\Theta_{k(\lambda)}$ they have bounded diameter in the $d_{k(\lambda)}$ distance. Therefore $d_{k(\lambda)}(F_{[w]_{i}}p,x)\leq k_{0}$ and the triangle inequality ensures $d_{k(\lambda)}(F_{[w]_{i}}p,y)\leq d_{k(\lambda)}(x,y) + k_{0}$. Examining the term of~\eqref{lowerestimateoftermsinGseries} with this $i$ and $p$ and taking $q=p$ we find that $x$ is too close to $F_{[w]_{i}}p$ to be in an exceptional subcell for $\psi^{r_{[w]_{i}}\mu_{[w]_{i}}\lambda}_{F_{[w]_{i}p}}\circ F_{[w]_{i}}^{-1}$.  If  $y$ is not in an exceptional subcell for this function then
\begin{equation*}
    d_{k(\lambda)}(F_{[w]_{i}}p,F_{[w]_{i}}p) + d_{k(\lambda)}(F_{[w]_{i}}p,x)+ d_{k(\lambda)}(F_{[w]_{i}}p,y )
    \leq k_{0} + d_{k(\lambda)}(x,y) + k_{0}
    \end{equation*}
which verifies the desired bound~\eqref{desiredupperboundondistancefunctionstogetGlowerbd}.  However for $y$ to be in an exceptional subcell for $\psi^{r_{[w]_{i}}\mu_{[w]_{i}}\lambda}_{F_{[w]_{i}p}}\circ F_{[w]_{i}}^{-1}$ it must be in a boundary subcell for $F_{[w]_{i}}(X)$ and thus for $F_{w}(X)$.  We conclude that $i+1=|w|$ and therefore $y$ is within $d_{k(\lambda)}$ distance $k_{0}$ of $F_{[w]_{i}}q$ for some $q\in V_{1}\setminus V_{0}$. Looking at the term from~\eqref{lowerestimateoftermsinGseries} with this $i$, $p$ and $q$ we use the triangle inequality $d_{k(\lambda)}(F_{[w]_{i}}p,F_{[w]_{i}}q)\leq d_{k(\lambda)}(x,y) + 2k_{0}$ to obtain
\begin{equation}\label{comparexydistancetotwodifferentjunctionpoints}
    d_{k(\lambda)}(F_{[w]_{i}}p,F_{[w]_{i}}q) + d_{k(\lambda)}(F_{[w]_{i}}p,x)+ d_{k(\lambda)}(F_{[w]_{i}}q,y )
    \leq d_{k(\lambda)}(x,y) + 4k_{0}
    \end{equation}
which verifies~\eqref{desiredupperboundondistancefunctionstogetGlowerbd}.  Then $x$ is not in an exceptional subcell for $\psi^{r_{[w]_{i}}\mu_{[w]_{i}}\lambda}_{F_{[w]_{i}p}}\circ F_{[w]_{i}}^{-1}$ and $y$ is not in an exceptional subcell for $\psi^{r_{[w]_{i}}\mu_{[w]_{i}}\lambda}_{F_{[w]_{i}q}}\circ F_{[w]_{i}}^{-1}$, so the lower bound follows from~\eqref{comparexydistancetotwodifferentjunctionpoints}.

Thus far we have verified the lower bound of~\eqref{pathdecayforGonXthmestimateofG} in the case that $r_{w}\mu_{w}\lambda>\tilde{c}$ and under the assumption that if either $x$ or $y$ is in a boundary subcell of $F_{w}(X)$, then this subcell does not contain a point of $V_{0}$.  However in the case $r_{w}\mu_{w}\lambda\leq \tilde{c}$ the lower bound amounts to saying that $G^{(\lambda)}(x,y)$ is bounded below if neither $x$ nor $y$ is within a $\Theta_{k(\lambda)}$ cell containing a point of $V_{0}$, which is vacuous when $\lambda$ is so small that these cells cover $X$, and evident from continuity of $G^{(\lambda)}(x,y)$ on the bounded interval in $\lambda$ otherwise.

Next we examine what happens on boundary subcells which contain points from $V_{0}$.  In this case we seek lower bounds for the negatives of partial derivatives as in~\eqref{pathdecayforGonXthmestimateofpartialnG} and~\eqref{pathdecayforGonXthmestimateofpartialnpartialnnGoffdiag}.  Again we assume $r_{w}\mu_{w}\lambda>\tilde{c}$ so as to use the lower bounds that go into~\eqref{lowerestimateoftermsinGseries}.

Recall that if $x\in V_{0}$ and $y\neq x$ we can write $-\partial_{n}'G^{(\lambda)}(x,y)$ as the series in~\eqref{expressionforpartialnGlambdaassumofpartialnterms}.  For this sum we can use termwise estimates like~\eqref{lowerestimateoftermsinGseries}, because replacing the lower bound on $\psi^{(r_{[w]_{i}}\mu_{[w]_{i}}\lambda)}_{p}(F_{[w]_{i}}^{-1}x)$ in~\eqref{lowerestimateoftermsinGseries}  by the lower bound on its normal derivative from~\eqref{superexpdecayofpsiboundforpartialnpsi} produces the same result except without the factor $\lambda^{-1/(S+1)}$.  As in the reasoning leading to~\eqref{desiredupperboundondistancefunctionstogetGlowerbd}, given $x\in V_{0}$ and $y$ not in an exceptional cell containing a point of $V_{0}$, our lower bound on $-\partial_{n}'G^{(\lambda)}(x,y)$ will follow if we can find $i$, $p$, and $q$ with $y$ not in an exceptional subcell of $F_{[w]_{i}}(X)$ and such that~\eqref{desiredupperboundondistancefunctionstogetGlowerbd} holds. Moreover the situation is simpler than it was there, because $x$ will always be in an exceptional subcell of $F_{[w]_{i}}(X)$ and we can take $p$ so $F_{[w]_{i}}p=x$.  For the largest number $j$ such that $y\in F_{[w]_{j}}(X)$, if $y$ is in an exceptional subcell of $F_{[w]_{j}}(X)$ then it is not in an exceptional subcell of $F_{[w]_{j-1}}(X)$, so we may take $i$ to be either $j$ or $j-1$ as appropriate.  In either case it is apparent that~\eqref{desiredupperboundondistancefunctionstogetGlowerbd} is valid, establishing the lower estimate of~\eqref{pathdecayforGonXthmestimateofpartialnG} when $r_{w}\mu_{w}\lambda>\tilde{c}$.  To deal with the case $r_{w}\mu_{w}\lambda\leq\tilde{c}$ we fix $y\neq x$ and observe that $G^{(\lambda)}(y,y)$ converges to a positive constant as $\lambda\to0$ and that  $\bigl(-\partial_{n}'G^{(\lambda)}(x,y)\bigr)\bigl(G^{(\lambda)}(y,y)\bigr)^{-1}$ is positive and decreasing in $\lambda$ by an argument similar to that in Corollary~\ref{signandorderingofnormalderivs}.  An easy computation using the series~\eqref{expressionforpartialnGlambdaassumofpartialnterms} with $\lambda=0$  shows $-\partial_{n}'G^{(0)}(x,y)$ is positive, and we conclude that $-\partial_{n}'G^{(\lambda)}(x,y)$ is bounded below by a constant depending only on the fractal and harmonic structure.  Together with the above, this shows that the lower estimate of~\eqref{pathdecayforGonXthmestimateofpartialnG} is valid for all positive $\lambda$.

We can now use this lower bound to improve the $\lesssim$ in~\eqref{upperboundforterminvolvingpartialnofresolventagainstzeta} to $\simeq$, because it guarantees that $-\partial_{n}''G^{(\lambda)}(x,y)$ achieves size comparable to a constant depending only on the fractal and harmonic structure.  The lower estimate on the integral then follows by an argument like that in Lemma~\ref{estimateofetaintegral}.  However this and~\eqref{usefulexpressionfortwopartialderivsofresolvent} immediately imply the lower bound in~\eqref{pathdecayforGonXthmestimateofpartialnpartialnnGoffdiag}, so all lower bounds in the theorem have been proved.
\end{proof}

\begin{proof}[Proof of Corollary~\ref{pathdecayforGNeumannonXthm}]
According to Theorem~4.2 of~\cite{IoPeRoHuSt2010TAMS}, the Neumann resolvent is obtained from the Dirichlet resolvent via
\begin{equation*}
    G_{N}^{(\lambda)}(x,y)
    = G^{(\lambda)}(x,y) + \sum_{p,q\in V_{0}} C_{pq}^{(\lambda)} \eta_{p}^{(\lambda)}(x)\eta_{q}^{(\lambda)}(y)
    \end{equation*}
where $C_{pq}^{(\lambda)}$ are the entries of the matrix inverse to that with entries $\partial_{n}\eta_{q}^{(\lambda)}(p)$.
From Theorem~\ref{maintheoremoneta} the values $C_{pq}^{(\lambda)}$ are all positive and the diagonal values $C_{pp}^{(\lambda)}$ are comparable to $(1+\lambda)^{-1/(S+1)}$. From the same theorem it is evident that $\eta_{p}^{(\lambda)}(x)\eta_{q}^{(\lambda)}(y)$ is bounded above by $\exp \Bigl(-\kappa_{2}d_{k(\lambda)}(x,y) \Bigr)$, so the terms added in passing from $G^{(\lambda)}(x,y)$ to $G_{N}^{(\lambda)}(x,y)$ also satisfy the upper bound~\eqref{pathdecayforGonXthmestimateofG}, and this estimate is therefore valid for $G_{N}^{(\lambda)}(x,y)$.  As we are adding positive terms, the lower bound  for $G^{(\lambda)}(x,y)$  in~\eqref{pathdecayforGonXthmestimateofG} applies also to $G_{N}^{(\lambda)}(x,y)$.  To complete the proof we need only check that this bound is still valid for $x$ and $y$ within a cell of scale $k(\lambda)$ at a boundary point.  At such points the exponential term is trivial, so the estimate is just that $G_{N}^{(\lambda)}(x,y)\gtrsim (1+\lambda)^{-1/(S+1)}$.  This is true because $\eta_{p}^{(\lambda)}$ is comparable to $1$ at both $x$ and $y$ from the estimates in Theorem~\ref{maintheoremoneta}, and $C_{pp}^{(\lambda)}\simeq(1+\lambda)^{-1/(S+1)}$ as previously noted.
\end{proof}

%%%%%%%%%%%%%%%%%%%%%%%%%%%%%%%%%%%%%%%%%%%%%%%%%%%%%%%%%%%%%%%%%%%%%%%%%%%%%%%%%%%%%%%%%%%%%%%%%%%%%%%%%%%%%%%%%%%%

\section{Resolvents on Blowups}\label{blowupsection}

\begin{definition}
For an infinite word $\w=\w_{1}\w_{2}\dotsm$ let
\begin{equation*}
    \Omega_{-n}=(F^{-1})_{[\w]_{n}}(X)= F_{\w_{1}}^{-1}\circ\dotsm\circ F_{\w_{n}}^{-1}(X).
    \end{equation*}
We call $\Omega_{-n}$ a finite blowup of $X$.  Observe that the sequence $\{\Omega_{-n}\}$ is
increasing because $F_{j}(X)\subset X$ for any $j$.  The blowup $\Omega$ of $X$ is the union
\begin{equation*}
    \Omega = \bigcup_{n\geq0} \Omega_{-n}.
    \end{equation*}
\end{definition}

In general there are uncountably many non-isometric blowups of a fractal.  If the word ends with
the infinite repetition of a single letter, then the blowup  will have non-empty boundary, but
otherwise there is no boundary. One reason to be interested in blowups is illustrated by the
familiar case of the unit interval $[0,1]$, which is a pcfss invariant set of the iterated function
system $F_{0}(x)= \frac{x}{2},\,F_{1}(x)=\frac{x+1}{2}$ on $\mathbb{R}$.  In this case any blowup
corresponding to a word that terminates with $000\dotsm$ is a ray to $+\infty$, while if the word
terminates with $111\dotsm$ it is a ray to $-\infty$; for all other words the blowup is all of
$\mathbb{R}$. Heuristically we think of the relationship between a fractal and its blowup as
somewhat akin to the relationship between $[0,1]$ and $\mathbb{R}$, inasmuch as they are compact
and non-compact examples of sets with similar local structure.  Some interesting results about
blowups of fractals may be found in~\cite{Strichartz1998CJM,Teplyaev1998JFA}.

Fix an infinite word $\w$.  It is straightforward to define a Laplacian on the blowup $\Omega$ that
is consistent with that on $X$. Recall that for a finite word $w$, we have $\Delta(u\circ F_{w}) =
r_{w}\mu_{w} (\Delta u)\circ F_{w}$.  Now if $u$ is defined on $\Omega_{-n}$ then $u\circ
(F^{-1})_{[\w]_{n}}$ is defined on $X$, so to maintain this composition rule for the Laplacian we
must set
\begin{align}\label{defnofLaplacianonfiniteblowup}
    \Delta u
    &= r_{[\w]_{n}}\mu_{[\w]_{n}} \Bigl( \Delta \bigl(  u\circ (F^{-1})_{[\w]_{n}} \bigr) \Bigr)
        \circ \Bigl( (F^{-1})_{[\w]_{n}} \Bigr)^{-1} \notag\\
    &= r_{[\w]_{n}}\mu_{[\w]_{n}} \Bigl( \Delta \bigl(  u\circ (F^{-1})_{[\w]_{n}} \bigr) \Bigr)
        \circ F_{\w_{n}}\circ F_{\w_{n-1}}\circ \dotsm\circ F_{\w_{1}}.
    \end{align}
It is easy to check that this ensures the definition of $\Delta$ on $\Omega_{-n}$ is consistent
with that on any other $\Omega_{-m}$, so suffices as a definition of $\Delta$ on $\Omega$ which is
consistent with that on $X=\Omega_{0}$.  We may similarly define $\DF$ and $\mu$ on $\Omega$ (note
that then $\mu$ is infinite but $\sigma$-finite) so as to be consistent with their definitions on
$X$.

Now the result~\eqref{maintheoremfromIPRRS} from~\cite{IoPeRoHuSt2010TAMS} may be transferred to each
$\Omega_{-n}$. For convenience of notation, for the remainder of this section the function
$G^{(\lambda)}(x,y)$ of~\eqref{maintheoremfromIPRRS} will be renamed $G_{0}^{(\lambda)}(x,y)$ to
emphasize that it corresponds to the zero blowup $X=\Omega_{0}$.

\begin{lemma}\label{greenonfiniteblowup}
Suppose $z\in\mathbb{C}$ is such that there is no finite word $w$ for which
$(r_{[\w]_{n}}\mu_{[\w]_{n}})^{-1}r_{w}\mu_{w}z$ is a Dirichlet eigenvalue of $\Delta$ on
$X$. Then defining
\begin{equation}\label{defnofGmlambda}
    G^{(z)}_{-n}(x,y)
    =  r_{[\w]_{n}}^{-1} G^{(r_{[\w]_{n}}^{-1}\mu_{[\w]_{n}}^{-1}z)}_{0}
        (F_{\w_{n}}\circ\dotsm\circ F_{\w_{1}}x,F_{\w_{n}}\circ\dotsm\circ F_{\w_{1}}y)
    \end{equation}
the function
\begin{equation*}
    u(x)=\int_{\Omega_{-n}} G^{(z)}_{-n}(x,y)f(y)\, d\mu(y)
    \end{equation*}
solves $(z\id-\Delta)u=f$ on $\Omega_{-n}$ with Dirichlet conditions at $\partial \bigl(
(F^{-1})_{[\w]_{n}}(X)\bigr)=(F^{-1})_{[\w]_{n}}(V_{0})$.
\end{lemma}

\begin{proof}
By direct computation setting $x'=F_{\w_{n}}\circ\dotsm\circ F_{\w_{1}}x$ and similarly for $y$,
\begin{align*}
    (z\id-\Delta_{x}) u(x)
    &= (z\id-\Delta_{x}) r_{[\w]_{n}}^{-1}
        \int_{\Omega_{-n}}  G_{0}^{(r_{[\w]_{n}}^{-1}\mu_{[\w]_{n}}^{-1}z)}
        (F_{\w_{n}}\circ\dotsm\circ F_{\w_{1}}x,F_{\w_{n}}\circ\dotsm\circ F_{\w_{1}}y)  f(y)\,
    d\mu(y)\\
    &= (z\id-\Delta_{x}) r_{[\w]_{n}}^{-1}\mu_{[\w]_{n}}^{-1}
        \int_{X} G_{0}^{(r_{[\w]_{n}}^{-1}\mu_{[\w]_{n}}^{-1}z)}(F_{\w_{n}}\circ\dotsm\circ F_{\w_{1}}x,y')
        f\bigl( (F^{-1})_{[\w]_{n}} y'\bigr)\, d\mu(y')\\
    &=(r_{[\w]_{n}}^{-1}\mu_{[\w]_{n}}^{-1}z - \Delta_{x'})
        \int_{X} G_{0}^{(r_{[\w]_{n}}^{-1}\mu_{[\w]_{n}}^{-1}z)}(x',y')
        f\bigl( (F^{-1})_{[\w]_{n}} y'\bigr)\, d\mu(y')\\
    &=f\bigl( (F^{-1})_{[\w]_{n}} x'\bigr)=f(x),
    \end{align*}
because $(r_{[\w]_{n}}^{-1}\mu_{[\w]_{n}}^{-1}z)$ satisfies the hypotheses of Proposition~\ref{propexistofeta}. The Dirichlet boundary conditions are immediate from the original result.
\end{proof}

\begin{theorem}\label{pathdecayforGonblowupthm}
For any $\lambda>0$ the sequence $G_{-n}^{(\lambda)}(x,y)$ is uniformly convergent on $\Omega\times\Omega$. The limit $G^{(\lambda)}_{\infty}(x,y)$ is such that if $f\in L^{1}(\Omega,d\mu)$ then
\begin{equation*}
    u(x)=\int_{\Omega} G^{(\lambda)}_{\infty}(x,y) f(y) \, d\mu(y)
    \end{equation*}
satisfies $(\lambda\id-\Delta)u=f$ on $\Omega$.  If $\Omega$ has boundary points then $u=0$ at these points.  Moreover there are positive constants $\kappa_{3}$ and $\kappa_{4}$ depending only on the fractal and harmonic structure and such that
\begin{equation}\label{pathdecayforGonXthmestimateofGonblowup}
    \lambda^{-1/(S+1)} \exp\Bigl( - \kappa_{3}d_{k(\lambda)}\bigl( x,y \bigr) \Bigr)
    \lesssim G^{(\lambda)}_{\infty}(x,y)
    \lesssim \lambda^{-1/(S+1)} \exp\Bigl( - \kappa_{4}d_{k(\lambda)}\bigl( x,y \bigr) \Bigr).
    \end{equation}
unless $x$ or $y$ is in a cell of $\Theta_{k(\lambda)}$ that contains a boundary point.  If there is a boundary point $p$ and $y$ is not in such a cell, then we have
\begin{equation}
    \exp\Bigl( - \kappa_{3} d_{k(\lambda)}\bigl( p,y \bigr) \Bigr)
    \lesssim -\partial_{n}'  G^{(\lambda)}_{\infty}(p,y)
    \lesssim \exp\Bigl( - \kappa_{4} d_{k(\lambda)}\bigl( p,y \bigr) \Bigr)
    \label{pathdecayforGonXthmestimateofpartialnGonblowup}
    \end{equation}
and if $p,q$ are boundary points then
\begin{equation}
    \lambda^{1/(S+1)}\exp\Bigl( - \kappa_{3} d_{k(\lambda)}\bigl( p,q \bigr) \Bigr)
    \lesssim \partial_{n}''\partial_{n}'  G^{(\lambda)}_{\infty}(p,q)
    \lesssim \lambda^{1/(S+1)} \exp\Bigl( - \kappa_{4} d_{k(\lambda)}\bigl( p,q \bigr) \Bigr).
    \label{pathdecayforGonXthmestimateofpartialnpartialnnGoffdiagonblowup}
    \end{equation}
\end{theorem}
\begin{proof}
We need only establish the convergence and the bounds; the other assertions follow immediately from
Lemma~\ref{greenonfiniteblowup}, which is applicable because $\lambda>0$ and the spectrum of
$\Delta$ contains only negative values.  Observe that instead of $x,y$ in a general compact set it
suffices to consider $x,y\in X=\Omega_{0}$, because the proof is the same for any $\Omega_{-n}$ and
any compact set is contained in some sufficiently large $\Omega_{-n}$.

For the moment we assume the convergence and examine the bounds.  These are readily deduced from the estimates in Theorem~\ref{pathdecayforGonXthm} applied to~\eqref{defnofGmlambda}. Direct substitution into~\eqref{pathdecayforGonXthmestimateofG} shows that both the upper and lower bound for $G^{(\lambda)}_{-n}(x,y)$ are of the form
\begin{equation*}
    r_{[\w]_{n}}^{-1} \bigl(1+r_{[\w]_{n}}^{-1}\mu_{[\w]_{n}}^{-1}\lambda\bigr)^{-1/(S+1)}
        \exp\Bigl( - \kappa d_{k(r_{[\w]_{n}}^{-1}\mu_{[\w]_{n}}^{-1}\lambda)}\bigl( F_{\w_{n}}\circ\dotsm\circ F_{\w_{1}}x,F_{\w_{n}}\circ\dotsm\circ F_{\w_{1}}y \bigr) \Bigr)
    \end{equation*}
however
\begin{equation*}
    r_{[\w]_{n}}^{-1} \bigl(1+r_{[\w]_{n}}^{-1}\mu_{[\w]_{n}}^{-1}\lambda\bigr)^{-1/(S+1)}
    = \bigl(r_{[\w]_{n}}\mu_{[\w]_{n}}+ \lambda\bigr)^{-1/(S+1)}
    \to \lambda^{-1/(S+1)}
    \end{equation*}
and
\begin{equation}\label{distancecomparabilityremarkrestated}
    d_{k(r_{[\w]_{n}}^{-1}\mu_{[\w]_{n}}^{-1}\lambda)}\bigl( F_{\w_{n}}\circ\dotsm\circ F_{\w_{1}}x,F_{\w_{n}}\circ\dotsm\circ F_{\w_{1}}y \bigr)
    \simeq d_{k(\lambda)}\bigl( x,y \bigr)
    \end{equation}
as described in the proof of~\eqref{distancecomparabilityremark}.  This demonstrates~\eqref{pathdecayforGonXthmestimateofGonblowup}.

Now using the fact that the normal derivative of a function composed with $F_{\w_{n}}\circ\dotsm\circ F_{\w_{1}}$ is the normal derivative of the function times the factor $r_{[\w]_{n}}$ we see from~\eqref{pathdecayforGonXthmestimateofpartialnG} that the upper and lower bound for $-\partial_{n}'G^{(\lambda)}_{-n}(p,y)$ are both of the form
\begin{equation*}
        \exp\Bigl( - \kappa d_{k(r_{[\w]_{n}}^{-1}\mu_{[\w]_{n}}^{-1}\lambda)}\bigl( F_{\w_{n}}\circ\dotsm\circ F_{\w_{1}}x,F_{\w_{n}}\circ\dotsm\circ F_{\w_{1}}y \bigr) \Bigr)
    \end{equation*}
so~\eqref{pathdecayforGonXthmestimateofpartialnGonblowup} follows using~\eqref{distancecomparabilityremarkrestated}.
Using the scaling of the normal derivative twice we have from~\eqref{pathdecayforGonXthmestimateofpartialnpartialnnGoffdiag} that the upper and lower bounds for $\partial_{n}''\partial_{n}'G^{(\lambda)}_{-n}(p,p)$ are both of the form
\begin{equation*}
    r_{[\w]_{n}} \bigl(1+r_{[\w]_{n}}^{-1}\mu_{[\w]_{n}}^{-1}\lambda\bigr)^{1/(S+1)}
    \exp\Bigl( - \kappa d_{k(r_{[\w]_{n}}^{-1}\mu_{[\w]_{n}}^{-1}\lambda)}\bigl( F_{\w_{n}}\circ\dotsm\circ F_{\w_{1}}x,F_{\w_{n}}\circ\dotsm\circ F_{\w_{1}}y \bigr) \Bigr)
    \end{equation*}
so using
\begin{equation*}
    r_{[\w]_{n}} \bigl(1+r_{[\w]_{n}}^{-1}\mu_{[\w]_{n}}^{-1}\lambda\bigr)^{1/(S+1)}
    = \bigl(r_{[\w]_{n}}\mu_{[\w]_{n}}+ \lambda\bigr)^{1/(S+1)}
    \to \lambda^{1/(S+1)}
    \end{equation*}
as $n\to\infty$ and~\eqref{distancecomparabilityremarkrestated} again we get~\eqref{pathdecayforGonXthmestimateofpartialnpartialnnGoffdiagonblowup}.  Note that these latter estimates are only of interest if there is at least one fixed point $p$ that is in all but finitely many $\Omega_{-n}$, from which it follows that $p$ is a boundary point of $\Omega$.

In order to establish the convergence it helps to reorganize the notation so that it is clear what we are summing.
From~\eqref{defnofGmlambda} and~\eqref{maintheoremfromIPRRS},
\begin{align*}
    G^{(\lambda)}_{-n}(x,y)
    &= r_{[\w]_{n}}^{-1} G^{(r_{[\w]_{n}}^{-1}\mu_{[\w]_{n}}^{-1}\lambda)}_{0}
        (F_{\w_{n}}\circ\dotsm\circ F_{\w_{1}}x,F_{\w_{n}}\circ\dotsm\circ F_{\w_{1}}y)\\
    &= r_{[\w]_{n}}^{-1} \sum_{w\in W_{\ast}} r_{w}
        \Psi^{(r_{w}\mu_{w}r_{[\w]_{n}}^{-1}\mu_{[\w]_{n}}^{-1}\lambda)}
        (F_{w}^{-1} \circ F_{\w_{n}}\circ\dotsm\circ F_{\w_{1}}x,F_{w}^{-1}\circ F_{\w_{n}}\circ\dotsm\circ F_{\w_{1}}y)
    \end{align*}
however if $w$ has length $m$ then $F_{w}^{-1}\circ F_{\w_{n}}\circ\dotsm\circ F_{\w_{1}} = F_{w_{m}}^{-1}\circ \dotsm\circ F_{w_{1}}^{-1} \circ F_{\w_{n}}\circ\dotsm\circ F_{\w_{1}}$,
so produces a non-trivial term in the sum only if  $w$ begins with $\w_{n}$, $\w_{n}\w_{n-1}$, and
so on. More precisely, if we think of $\w$ as infinite to the left, $\w=\dotsm\w_{3}\w_{2}\w_{1}$,
and write $[\w]_{-n}=\w_{n}\dotsm\w_{1}$ and $[\w]_{-n,m}=\w_{n}\dotsm\w_{n-m+1}$ then the words
producing a non-trivial term are of the form $[\w]_{-n,m}$ for $0\leq m<n$, or $[\w]_{-n}w$ for
some $w\in W_{\ast}$. Consequently
\begin{align*}
    \lefteqn{G^{(\lambda)}_{-n}(x,y)}\quad&\\
    &= r_{[\w]_{-n}}^{-1} \sum_{m=0}^{n-1} r_{[\w]_{-n,m}}
        \Psi^{(r_{[\w]_{-n,m}}\mu_{[\w]_{-n,m}}r_{[\w]_{-n}}^{-1}\mu_{[\w]_{-n}}^{-1}\lambda)}
        \bigl((F_{[\w]_{-n,m}})^{-1} \circ F_{[\w]_{-n}}x
            ,(F_{[\w]_{-n,m}})^{-1}\circ F_{[\w]_{-n}}y\bigr) \notag\\
    &\quad+  r_{[\w]_{-n}}^{-1} \sum_{w\in W_{\ast}} r_{[\w]_{-n}w}
        \Psi^{(r_{[\w]_{-n}w}\mu_{[\w]_{-n}w}r_{[\w]_{-n}}^{-1}\mu_{[\w]_{-n}}^{-1}\lambda)}
        \bigl((F_{[\w]_{-n}w})^{-1} \circ F_{[\w]_{-n}}x
            ,(F_{[\w]_{-n}w})^{-1}\circ F_{[\w]_{-n}}y\bigr) \notag\\
    &=\sum_{m=0}^{n-1} r_{[\w]_{m-n}}^{-1} \Psi^{(r_{[\w]_{m-n}}^{-1}\mu_{[\w]_{m-n}}^{-1}\lambda)}
        \bigl(F_{[\w]_{m-n}}x ,F_{[\w]_{m-n}}y\bigr)
        + \sum_{w\in W_{\ast}} r_{w}\Psi^{(r_{w}\mu_{w}\lambda)}
            (F_{w}^{-1} x, F_{w}^{-1}y) \notag\\
    &=\sum_{m=-n}^{-1} r_{[\w]_{m}}^{-1} \Psi^{(r_{[\w]_{m}}^{-1}\mu_{[\w]_{m}}^{-1}\lambda)}
        \bigl(F_{[\w]_{m}}x ,F_{[\w]_{m}}y\bigr)
        +  \sum_{w\in W_{\ast}} r_{w}\Psi^{(r_{w}\mu_{w}\lambda)}
            (F_{w}^{-1} x, F_{w}^{-1}y).
    \end{align*}
In particular, for $n'>n$ we have from~\eqref{PsiintermsofGpq}
\begin{align*}
    G^{(\lambda)}_{-n'}(x,y)-G^{(\lambda)}_{-n}(x,y)
    &= \sum_{m=-n'}^{-n-1} r_{[\w]_{m}}^{-1} \Psi^{(r_{[\w]_{m}}^{-1}\mu_{[\w]_{m}}^{-1}\lambda)}
        \bigl(F_{[\w]_{m}}x ,F_{[\w]_{m}}y\bigr)\\
    &= \sum_{m=-n'}^{-n-1} \sum_{p,q\in V_{1}\setminus V_{0}}  r_{[\w]_{m}}^{-1}
    G^{(r_{[\w]_{m}}^{-1}\mu_{[\w]_{m}}^{-1}\lambda)}_{pq}\:
    \psi^{(r_{[\w]_{m}}^{-1}\mu_{[\w]_{m}}^{-1}\lambda)}_{p}(F_{[\w]_{m}}x)\: \psi^{(r_{[\w]_{m}}^{-1}\mu_{[\w]_{m}}^{-1}\lambda)}_{q}(F_{[\w]_{m}}y).
    \end{align*}
From Lemma~\ref{Gpqestimate},
\begin{equation*}
    r_{[\w]_{m}}^{-1}
    G^{(r_{[\w]_{m}}^{-1}\mu_{[\w]_{m}}^{-1}\lambda)}_{pq}
    \lesssim  \bigl(r_{[\w]_{m}}^{-1}\mu_{[\w]_{m}}^{-1}+\lambda\bigr)^{\frac{-1}{S+1}}
        \exp \Bigl( - c_{2} d_{k(r_{[\w]_{m}}^{-1}\mu_{[\w]_{m}}^{-1}\lambda)}(p,q)  \Bigr)
    \lesssim  \lambda^{\frac{-1}{S+1}}
    \end{equation*}
and applying Lemma~\ref{superexpdecayofpsi} we find
\begin{align}
    \Bigl| G^{(\lambda)}_{-n'}(x,y)-G^{(\lambda)}_{-n}(x,y)\Bigr|
    &\lesssim \lambda^{-1/(S+1)} \sum_{m=-n'}^{-n-1} \sum_{p,q\in V_{1}\setminus V_{0}}
        \psi^{(r_{[\w]_{m}}^{-1}\mu_{[\w]_{m}}^{-1}\lambda)}_{p}(F_{[\w]_{m}}x)\: \psi^{(r_{[\w]_{m}}^{-1}\mu_{[\w]_{m}}^{-1}\lambda)}_{q}(F_{[\w]_{m}}y) \notag\\
    &\lesssim \lambda^{-1/(S+1)} \sum_{m=-n'}^{-n-1} \sum_{p,q\in V_{1}\setminus V_{0}}
        \exp\Bigl( - c_{4} d_{k(\lambda)}\bigl(F_{[\w]_{m}}^{-1}p,x\bigr) - c_{4}d_{k(\lambda)}\bigl(F_{[\w]_{m}}^{-1}q,y\bigr)
        \Bigr). \label{boundfordifferenceofGnGnprimeonblowup}
    \end{align}
Recall that we were able to assume $x,y\in \Omega_{0}$.  Now we wish to count how many $F_{[\w]_{m}}^{-1}p$ and $F_{[\w]_{m}}^{-1}q$ are within some distance of $\Omega_{0}$.  It is not easy to do this using the $d_{k(\lambda)}$ distance, but it is easy to make an estimate using the number of cells that have the same scale as $\Omega_{0}$, by which we mean the cells of the form $F_{[\w]_{m}}^{-1}\circ F_{w}(\Omega_{0})$ where $|w|=m$.  For cells of the same scale as $\Omega_{0}$, we think of those that intersect $\Omega_{0}$ as forming an annulus of size $1$, those intersecting this annulus (but not $\Omega_{0}$ itself) as forming an annulus size $2$, and so forth.  It is apparent that the number of points in the intersection of $\bigl\{F_{[\w]_{m}}^{-1}p: p\in V_{1}\setminus V_{0}, m\geq1\}$ with any such annulus is bounded by a constant depending only on the fractal.  We can also obtain a crude estimate on the growth of the $d_{k(\lambda)}$ distance from $\Omega_{0}$ to the $n^{\text{th}}$ annulus, just by noting that the resistance distance to this annulus must grow at least like $(\max_{j} r_{j})^{-n}$, and the $d_{k(\lambda)}(x,y)$ distance is bounded below by $(1+\lambda)^{1/(S+1)}$ times the resistance distance, so also grows geometrically.  It follows that~\eqref{boundfordifferenceofGnGnprimeonblowup} converges (pointwise) as $n'\to\infty$ and satisfies a bound
\begin{equation*}
    \Bigl| G^{(\lambda)}_{\infty}(x,y)-G^{(\lambda)}_{-n}(x,y)\Bigr|
    \lesssim C(\lambda)
        \exp\biggl( - c_{4} \min\Bigl\{d_{k(\lambda)}\bigl(F_{[\w]_{m}}^{-1}p,x\bigr): m\geq n,\, p\in V_{1}\setminus V_{0} \Bigr\} \biggr).
    \end{equation*}
However we established above that there are only finitely many points $F_{[\w]_{m}}^{-1}p$ within any of our annuli. Thus for any prescribed distance we can take the finite union of annuli covering that distance around $\Omega_{0}$, and by taking $n$ large enough we can be sure no $F_{[\w]_{m}}^{-1}p$ with $m\geq n$ lies in this union.  Hence $G^{(\lambda)}_{-n}(x,y)$ converges uniformly to $G^{(\lambda)}_{\infty}(x,y)$.
\end{proof}

%%%%%%%%%%%%%%%%%%%%%%%%%%%%%%%%%%%%%%%%%%%%%%%%%%%%%%%%%%%%%%%%%%%%%%%%%%%%%%%%%%%%%%%%%%%%%%%%%%%%%%%%%%%%%%%%%%%%%%%

\section{Phragmen-Lindel\"{o}f type theorems}\label{plsection}
In this section we prove some complex analytic estimates related to the Phragmen-Lindel\"{o}f theorem.  Later these will be used to prove bounds for the resolvent in a sector of $\mathbb{C}$ that omits the negative real axis.  Phragmen-Lindel\"{o}f theorems have previously been used to obtain off-diagonal decay estimates for heat kernels from Davies-Gaffney estimates~\cite{CoulhSikor2008PLMS}, however it seems that the techniques and results proved there are not applicable in our situation.

Fix an angle
$\alpha<\pi$ and consider the sector
\begin{equation*}
    \Sect = \bigl\{ z=|z|e^{i\beta}\in\mathbb{C}: -\alpha< \beta <\alpha \bigr\}.
    \end{equation*}
It will also be convenient to identify $\Sect^{+}=\Sect\cap\{\Im z>0\}$ and
$\Sect^{-}=\Sect\cap\{\Im z<0\}$.

One version of the classical Phragmen-Lindel\"{o}f theorem is as follows; a proof may be found in~\cite{Levin1996}.
\begin{theorem}[Phragmen--Lindel\"{o}f]
Suppose that $u(z)$ is a function analytic in an open sector of angular size $\alpha$ and continuous on the closure of the sector.  If $u$ is bounded by $M$ on the sides of the sector and satisfies $|u(z)|\leq C\exp (c|z|^{\alpha'})$ for some $\alpha'<\frac{\pi}{\alpha}$ and constants $C$ and $c$, then $|u|\leq M$ on the closed sector.
\end{theorem}
As a particular consequence we see that
\begin{corollary}\label{affinenestedversionofPLcorol}
Suppose $g(z)$ is analytic  on $\Sect^{+}$ and continuous and bounded by $1$ on the closure. If  $|g(z)|\leq \exp (-a_{1}|z|^{a_{2}})$ on the positive real axis for some $0<a_{2}\leq 1$ and $a_{1}>0$, then for $\beta\in[0,\alpha]$ we have
\begin{equation*}
    |g(|z|e^{i\beta})|
    \leq \exp \biggl( \frac{-a_{1}\sin\bigl( a_{2}(\alpha-\beta)\bigr) }{\sin (a_{2}\alpha) } |z|^{a_{2}} \biggr)
    \end{equation*}
\end{corollary}
\begin{proof}
With the convention that all power functions are defined by cutting the
$z$-plane along the negative real axis, let
\begin{equation*}
    v(z)
    =\exp\biggl( \frac{a_{1}e^{i(\pi/2-a_{2}\alpha)}}{\sin (a_{2}\alpha) } z^{a_{2}} \biggr)
    \end{equation*}
which is analytic in $\Sect^{+}$.  Note that the real part of the exponent vanishes on $\Arg(z)=\alpha$, and that on the positive real axis it is $a_{1}|z|^{a_{2}}$.  Thus $|v(z)g(z)|\leq1$ on $\partial\Sect^{+}$.  It is also apparent that
\begin{equation*}
    |v(z)|
    \leq \exp \biggl( \frac{a_{1}}{\sin (a_{2}\alpha) } |z|^{a_{2}} \biggr)
    \end{equation*}
within $\Sect^{+}$.  Since $a_{2}\leq1$ and the sector has angle $\alpha<\pi$, the Phragmen-Lindel\"{o}f theorem implies $|v(z)g(z)|\leq 1$ on the closure of $\Sect^{+}$.  Thus for $z=|z|e^{i\beta}$, $0\leq\beta\leq\alpha$ we have
\begin{equation*}
    |g(z)|
    \leq |v(z)|^{-1}
    = \exp \biggl( \frac{-a_{1}\sin\bigl( a_{2}(\alpha-\beta)\bigr) }{\sin (a_{2}\alpha) } |z|^{a_{2}} \biggr)
    \end{equation*}
\end{proof}
We shall later see that this could be applied to both piecewise eigenfunctions and the resolvent in the case of affine nested fractals, because these have decay $\exp\bigl(-a_{1}d_{k_(\lambda)}(x,y)\bigr)$ for $\lambda$ on the positive real axis and $d_{k(\lambda)}(x,y)\simeq (1+\lambda)^{\gamma/(S+1)} R(x,y)^{\gamma}$ by Proposition~\ref{affinenesteedversionofdk} and Definition~\ref{defnofklambda}.  In the more general pcfss case it cannot be used because we do not know whether $d_{k_(\lambda)}(x,y)$ grows like a power of $\lambda$.  In order to deal with the more general decay $\exp\bigl(-a_{1}d_{k_(\lambda)}(x,y)\bigr)$ we need a modified version of Corollary~\ref{affinenestedversionofPLcorol}.  The modification takes into account the fact that $d_{k_(\lambda)}(x,y)$ behaves like a power of $\lambda$ over exponential scales, but that the powers could perhaps be different for distinct scales.  In essence, what we must do is create functions analytic on a sector and growing according to different powers on different exponential scales.  This is done using Schwarz-Christoffel functions. The author anticipates that the existence of maps might be known, but does not know of a place where they are described.  Our replacement for Corollary~\ref{affinenestedversionofPLcorol} in the general case is as follows.
\begin{theorem}\label{generalpcfssversionofPLcorol}
Fix $\alpha<\pi$.  Suppose $f:[0,\infty)\to[0,\infty)$ has the property that there are  $0<\beta_{1}<\beta_{2}<1$ and $c_{1}>0$ and $c_{2}>0$ such that
\begin{equation}\label{hypothesisforSCfunctionconstruction}
    c_{1}\scale^{\beta_{1}}
    \leq \frac{f(\scale^{j+1})}{f(\scale^{j})}
    \leq c_{2} \scale^{\beta_{2}}.
    \end{equation}
for any sufficiently large $\scale\in[1,\infty)$ and all $j\in\mathbb{N}\cup\{0\}$.  If $g(z)$ is a function that is analytic on $\Sect^{+}$, continuous and bounded by $1$ on the closure of $\Sect^{+}$, and satisfies  $|g(z)|\leq \exp \bigl(-f(|z|^{\pi/\alpha})\bigr)$ on the positive real axis, then there are constants $c$ and $\tilde{c}$  such that for $\beta\in[0,\alpha]$
\begin{equation*}
    |g(|z|e^{i\beta})|
    \leq \exp \biggl(-c f(|z|^{\pi/\alpha}) \sin\tilde{c}\Bigl(1-\frac{\beta}{\alpha}\Bigr)  \biggr)
    \end{equation*}
\end{theorem}
\begin{proof}
From Lemma~\ref{SCfunctionforcomparisonwithdkgrowth} below we know that under these hypotheses for every $0<\epsilon\leq\frac{1}{3}\min\bigl\{(\beta_{1}), (1-\beta_{2}), (\beta_{2}-\beta_{1})\bigr\}$ there is a function  $F(z)$ analytic in the upper half plane and satisfying
\begin{equation*}
    c_{3}(\epsilon) f(|z|) \sin\bigl((\beta_{1}-\epsilon)(\pi-\beta)\bigr)
    \leq  \Re F(|z|e^{i\beta})
    \leq f(|z|)  \sin\bigl((\beta_{2}+\epsilon)(\pi-\beta)\bigr).
    \end{equation*}
Consider
\begin{equation*}
    v(z)=\exp\biggl( \frac{F(z^{\pi/\alpha})}{\sin\bigl(\pi(\beta_{2}+\epsilon)\bigr)} \biggr),
    \end{equation*}
which is analytic in $\Sect^{+}$.  Observe that $z^{\pi/\alpha}$ takes the ray with argument $\alpha$ to the negative real axis, so that $F(z^{\pi/\alpha})$ is imaginary and $|v(z)g(z)|\leq1$ on this ray.  Since the real part of $F(z^{\pi/\alpha})$ is bounded by $f(|z|^{\pi/\alpha})\sin\bigl(\pi(\beta_{2}+\epsilon)\bigr)$ on the positive real axis, we see that $|v(z)g(z)|\leq1$ on $\partial\Sect^{+}$.  We also have the bound
\begin{equation*}
    |v(z)| \leq \exp\biggl( \frac{f(|z|^{\pi/\alpha})}{\sin\bigl(\pi(\beta_{2}+\epsilon)\bigr)}\biggr)
    \end{equation*}
on $\Sect^{+}$. Since $f(|z|^{\pi/\alpha})\leq c_{2}|z|^{\beta_{2}\pi/\alpha}$ for all sufficiently large $|z|$ and $\frac{\beta_{2}\pi}{\alpha}<\frac{\pi}{\alpha}$, the Phragmen-Lindle\"{o}f theorem implies that $|v(z)g(z)|\leq1$ on $\Sect^{+}$.  We conclude that
\begin{equation*}
    |g(z)|
    \leq |v(z)|^{-1}
    \leq \exp\biggl(  \frac{-c_{3}(\epsilon) f(|z|^{\pi/\alpha}) \sin\bigl((\beta_{1}-\epsilon)\pi(1-\frac{\beta}{\alpha})\bigr)}{\sin\bigl(\pi(\beta_{2}+\epsilon)\bigr)} \biggr)
    \end{equation*}
and note that the constants are non-zero.
\end{proof}

The construction and properties of the Schwarz-Christoffel function used in the proof are contained in the following lemma.
\begin{lemma}\label{SCfunctionhasrightests}
Given $0<\alpha_{1}<\alpha_{2}<1$ and numbers $\tau_{j}$ such that that $\alpha_{1}<\sum_{j=0}^{k}\tau_{j}<\alpha_{2}$ for all $k$, and $\sup_{j}|\tau_{j}|=C<1$, the Schwarz-Christoffel function
\begin{equation*}
    H(z) = \int_{0}^{z} w^{-\tau_{0}} \prod_{j=1}^{\infty} \Bigl( 1-\frac{w}{\scale^{j}}\Bigr)^{-\tau_{j}} \, dw
    \end{equation*}
is holomorphic in the upper half plane and satisfies
\begin{align*}
    \lefteqn{|z|^{1-\sum_{j=0}^{K}(1-j/K)\tau_{j}}
        \min\Bigl\{\sin\bigl(\alpha_{1}\pi+(1-\alpha_{1})\beta \bigr) ,\sin\bigl(\alpha_{2}\pi+(1-\alpha_{2})\beta \bigr)  \Bigr\}}
    \quad& \\
    &\lesssim \Im H(|z|e^{i\beta})\\
    &\lesssim |z|^{1-\sum_{j=0}^{K}(1-j/K)\tau_{j}}
        \max\Bigl\{\sin\bigl(\alpha_{1}\pi+(1-\alpha_{1})\beta \bigr) ,\sin\bigl(\alpha_{2}\pi+(1-\alpha_{2})\beta \bigr)  \Bigr\}.
    \end{align*}
for every $\beta\in[0,\pi]$, with constants depending only upon the assumed bounds and on $\scale$.
\end{lemma}
\begin{proof}
It is well known that $H(z)$ maps the upper half plane conformally onto the region bounded by a polygonal path with angles $\pi\tau_{j}$ at the vertices and extends continuously to the boundary, so that the negative real axis is mapped to itself and at $w\in(\scale^{k},\scale^{k+1})\subset\mathbb{R}$ the direction of the boundary path is $\pi\sum_{0}^{k}\tau_{j}$, which by hypothesis is in $[\alpha_{1}\pi,\alpha_{2}\pi]$.  We will need to know an estimate for the direction of the image of a ray from $0$ at angle $\beta\in(0,\pi)$.  At a point $w$ on this ray, the angles $\beta_{j}$ from $w$ to $\scale^{j}\in\mathbb{R}$ form a strictly increasing sequence beginning at $\beta_{0}=-\pi+\beta$ and converging to $0$.  The direction tangent to the image curve at $H(w)$ is then $\beta+\sum_{j}(-\tau_{j})\beta_{j}$.  Now
\begin{equation*}
    \beta+\sum_{j=0}^{k}(-\tau_{j})\beta_{j}
    =\beta + (\pi-\beta)\sum_{j=0}^{k}\tau_{j} + \sum_{j=0}^{k} \tau_{j}(\beta-\pi-\beta_{j})
    \end{equation*}
and $\beta-\pi-\beta_{j}\leq0$ for all $j$, so the second sum has the opposite sign to $\sum_{j}\tau_{j}>0$.  It follows that
\begin{equation*}
    \beta+\sum_{j=0}^{k}(-\tau_{j})\beta_{j}
    \leq \beta + (\pi-\beta)\sum_{j=0}^{k}\tau_{j}
    \leq \beta + (\pi-\beta)\alpha_{2}
    =\pi - (1-\alpha_{2})(\pi-\beta).
    \end{equation*}
To see a lower bound, write
\begin{equation*}
    \beta+\sum_{j=0}^{k}(-\tau_{j})\beta_{j}
    =\beta + (\pi-\beta)\tau_{0} + (\pi-\beta)\sum_{j=1}^{k}\tau_{j} + \sum_{j=0}^{k} \tau_{j}(\beta-\pi-\beta_{j}).
    \end{equation*}
If $\sum_{1}^{k}\tau_{j}\geq0$, then so is $\sum_{1}^{k}-\tau_{j}\beta_{j}$, so
\begin{equation*}
    \beta+\sum_{j=0}^{k}(-\tau_{j})\beta_{j}
    \geq\beta+(\pi-\beta)\tau_{0}
    \geq\beta+(\pi-\beta)\alpha_{1}.
    \end{equation*}
However, at any $k$ such that $\sum_{1}^{k}\tau_{j}<0$ we have $\sum_{j=0}^{k} \tau_{j}(\beta-\pi-\beta_{j})>0$, so that the same lower bound
\begin{equation*}
    \beta+\sum_{j=0}^{k}(-\tau_{j})\beta_{j}
    \geq \beta+(\pi-\beta) \sum_{j=0}^{k}\tau_{j}
    \geq \beta+(\pi-\beta)\alpha_{1}
    = \pi - (1-\alpha_{1})(\pi-\beta)
    \end{equation*}
holds.  Thus the direction of the image of the ray from $0$ at angle $\beta$ is between $\pi - (1-\alpha_{1})(\pi-\beta)$ and $\pi - (1-\alpha_{2})(\pi-\beta)$.  In particular we may relate the integral along this ray to the integral of the magnitude of the integrand along the ray.  Writing $h(w)= w^{-\tau_{0}} \prod_{j=1}^{\infty} \bigl( 1-w\scale^{-j}\bigr)^{-\tau_{j}}$ and $z=|z|e^{i\beta}$ we have
\begin{align}\label{heatfromsector_compareIMofSCtomagnitude}
    \lefteqn{\min\Bigl\{\sin\bigl(\pi - (1-\alpha_{1})(\pi-\beta) \bigr) ,\sin\bigl(\pi - (1-\alpha_{2})(\pi-\beta) \bigr)  \Bigr\}
        \int_{0}^{|z|} \Bigl| h(te^{i\beta} ) \Bigr|\, dt}  \quad & \notag \\
    &\leq \Im \int_{0}^{|z|} h(te^{i\beta}) \, d(te^{i\beta}) = \Im H(z) \notag\\
    &\leq  \max\Bigl\{\sin\bigl(\pi - (1-\alpha_{1})(\pi-\beta) \bigr) ,\sin\bigl(\pi - (1-\alpha_{2})(\pi-\beta) \bigr)  \Bigr\}
        \int_{0}^{|z|} \Bigl| h(te^{i\beta} ) \Bigr|\, dt.
    \end{align}

To proceed we need an estimate on the magnitude of the integrand $h(w)$.  For fixed $z$ let $k$ be such that $\scale^{k-1/2}\leq|z|\leq\scale^{k+1/2}$ and write
\begin{align*}
    \bigl|h(w) \bigr|
    &= |w|^{-\tau_{0}}\Bigl| 1-\frac{w}{\scale^{k}}\Bigr|^{-\tau_{k}} \prod_{j=1}^{k-1} \Bigl|\frac{w}{\scale^{j}} \Bigr|^{-\tau_{j}}
        \prod_{j=1}^{k-1} \Bigl| 1-\frac{\scale^{j}}{w}\Bigr|^{-\tau_{j}}
        \prod_{j=k+1}^{\infty} \Bigl| 1-\frac{w}{\scale^{j}}\Bigr|^{-\tau_{j}}\\
    &= |w|^{-\tau_{0}}\Bigl| 1-\frac{w}{\scale^{k}}\Bigr|^{-\tau_{k}} \prod_{j=1}^{k-1} \Bigl|\frac{w}{\scale^{j}} \Bigr|^{-\tau_{j}}
        \exp \biggl( \sum_{j=1}^{k-1} -\tau_{j} \log \Bigl| 1-\frac{\scale^{j}}{w} \Bigr|
            + \sum_{j=k+1}^{\infty} -\tau_{j} \log \Bigl| 1-\frac{w}{\scale^{j}} \Bigr| \biggr).
    \end{align*}
Now it is easy to see that
\begin{align*}
    \biggl| \sum_{j=1}^{k-1} -\tau_{j} \log \Bigl| 1-\frac{\scale^{j}}{w} \Bigr| \biggr|
    &\leq  \frac{C}{|w|} \sum_{j=1}^{k-1} \scale^{j} \leq \tilde{C}\\
    \biggl| \sum_{j=k+1}^{\infty} -\tau_{j} \log \Bigl| 1-\frac{w}{\scale^{j}} \Bigr| \biggr|
    &\leq  C|w| \sum_{j=k+1}^{\infty} \scale^{-j} \leq \tilde{C}
    \end{align*}
so that
\begin{equation*}
    \bigl|h(w) \bigr|
    \simeq \Bigl| 1-\frac{w}{\scale^{k}}\Bigr|^{-\tau_{k}} |w|^{-\tau_{0}} \prod_{j=1}^{k-1} \Bigl|\frac{w}{\scale^{j}} \Bigr|^{-\tau_{j}}
    \simeq \Bigl| 1-\frac{w}{\scale^{k}}\Bigr|^{-\tau_{k}} \prod_{j=0}^{k-1} \scale^{(j-k)\tau_{j}}
    \end{equation*}
with constants independent of $k$ (here we used boundedness of $\sum_{0}^{k}\tau_{j}$).
We may use this to estimate the integral in~\eqref{heatfromsector_compareIMofSCtomagnitude} from $t=A^{k-1/2}$ to $t=A^{k+1/2}$.  The fact that $|\tau_{k}|<1$ ensures integrability of the first term, and that the integration introduces a bounded multiple of $\scale^{k}$ into the product.  Letting $K$ be the integer part of $\log|z|/\log\scale$ we find that
\begin{equation*}
    \int_{\scale^{-1/2}}^{|z|} \Bigl| h(te^{i\beta} ) \Bigr|\, dt
    \simeq \sum_{k=0}^{K}  \scale^{k}\prod_{j=0}^{k-1} \scale^{(j-k)\tau_{j}}
    \geq \scale^{K}\prod_{j=0}^{K-1} \scale^{(j-K)\tau_{j}}.
    \end{equation*}
However the ratio of the $(k+1)$-th term to the $k$-th term is
\begin{equation*}
    \left( \scale^{k+1}\prod_{j=0}^{k} \scale^{(j-k-1)\tau_{j}}\right)\left( \scale^{-k}\prod_{j=0}^{k-1} \scale^{(k-j)\tau_{j}}\right)
    = \scale^{1-\sum_{0}^{k}\tau_{j}}
    \geq  \scale^{1-\alpha_{2}}
    \end{equation*}
by hypothesis.  Thus we have the upper bound
\begin{equation*}
    \sum_{k=0}^{K}  \scale^{k}\prod_{j=0}^{k-1} \scale^{(j-k)\tau_{j}}
    \leq \left(\scale^{K}\prod_{j=0}^{K-1} \scale^{(j-K)\tau_{j}} \right)
        \sum_{k=0}^{K} \scale^{-k(1-\alpha_{2})}
    \lesssim \scale^{K}\prod_{j=0}^{K-1} \scale^{(j-K)\tau_{j}}
    \end{equation*}
and finally may conclude (again using boundedness of $\sum_{0}^{k}\tau_{j}$)
\begin{equation*}
    \int_{\scale^{0}}^{|z|} \Bigl| h(te^{i\beta} ) \Bigr|\, dt
    \simeq |z|^{1-\sum_{j=0}^{K}(1-j/K)\tau_{j}}.
    \end{equation*}
from which the result follows by~\eqref{heatfromsector_compareIMofSCtomagnitude}.
\end{proof}

The construction in the preceding lemma is applicable under the hypotheses of the theorem, as shown in the following.

\begin{lemma}\label{SCfunctionforcomparisonwithdkgrowth}
Suppose $f$ satisfies the hypotheses of Theorem~\ref{generalpcfssversionofPLcorol}.  For any $0<\epsilon\leq\frac{1}{3}\min\bigl\{(\beta_{1}), (1-\beta_{2}), (\beta_{2}-\beta_{1})\bigr\}$ there is a constant $c_{3}=c_{3}(\epsilon)$ and a function $F(z)$, holomorphic on the upper half plane, such that for all $\beta\in[0,\pi]$
\begin{equation*}
    c_{3}f(|z|) \sin\bigl((\beta_{1}+\epsilon)(\pi-\beta)\bigr)
    \leq  \Re F(|z|e^{i\beta})
    \leq f(|z|)  \sin\bigl((\beta_{2}-\epsilon)(\pi-\beta)\bigr).
    \end{equation*}
\end{lemma}
\begin{proof}
By hypothesis we have
\begin{equation*}
    \beta_{1} - \frac{c_{1}}{\log\scale}
    \leq \frac{ \log f(\scale^{i}) - \log f(\scale^{i-1})} {\log \scale }
    \leq \beta_{2} + \frac{c_{2}}{\log\scale}.
    \end{equation*}
Choose $\scale$ sufficiently large that both $\frac{c_{1}}{\log\scale}$ and $\frac{c_{2}}{\log\scale}$ are smaller than $\epsilon$.  Then
\begin{equation}\label{neededestimateforsumoftau}
    0<\beta_{1}-\epsilon
    \leq \frac{ \log f(\scale^{i}) - \log f(\scale^{i-1})} {\log \scale }
    \leq \beta_{2}+\epsilon
    <1.
    \end{equation}
By making $\scale$ larger if necessary we also may require $\beta_{1}-\epsilon\leq  \frac{\log f(\scale) } {\log \scale}\leq \beta_{2}+\epsilon$, because $\log (f(1))/\log \scale\to 0$ as $\scale\to\infty$.

Now let
\begin{equation*}
    l_{j}=\frac{\log f(\scale^{j+1})}{(j+1)\log\scale},\quad j\geq0.
    \end{equation*}
and
\begin{equation*}
    \tau_{j}
    = \begin{cases}
        1-l_{j_{0}} &\text{ if $j=0$}\\
        2l_{0}-2l_{1}  &\text{ if  $j=1$}\\
        \tau_{j}=j\bigl(2l_{j-1}-l_{j}-l_{j-2}\bigr)+ l_{j-2}-l_{j} &\text{ if $j\geq 2$}
        \end{cases}
    \end{equation*}
We claim that
\begin{equation}
    1 - \sum_{j=0}^{k}\bigl(1-j/k \bigr)\tau_{j}
    =l_{k-1} \text{ for $k\geq 1$}.
    \label{inductionontheljvalues}
    \end{equation}
This is easily verified when $k=1,2$.  For $k>2$ it follows by induction, because if it is true up to $k$ then the the following are true, and we must prove $X=(k+1)l_{k}$.
\begin{gather*}
    k - \sum_{j=0}^{k} (k-j)\tau_{j}
        = kl_{k-1}\\
    k-1 - \sum_{j=0}^{k-1} (k-1-j)\tau_{j}
        = (k-1)l_{k-2}\\
    k+1 - \sum_{j=0}^{k+1} (k+1-j)\tau_{j}
        = X.
    \end{gather*}
If we subtract the latter two equations from twice the first we obtain
\begin{equation*}
    \tau_{k}
    =2kl_{k-1} - (k-1)l_{k-2} -X
    \end{equation*}
from which we have the desired equality:
\begin{align*}
    X
    &= 2kl_{k-1} - (k-1)l_{k-2} -\tau_{k}\\
    &= 2kl_{k-1} - (k-1)l_{k-2} - k\bigl(2l_{k-1}-l_{k}-l_{k-2}\bigr)- l_{k-2}+l_{k}\\
    &= (k+1)l_{k}.
    \end{align*}

Our goal is to apply Lemma~\ref{SCfunctionhasrightests}, for which we need estimates on  $\sum_{0}^{i}\tau_{j}$.  Define $\alpha_{1}=1-\beta_{2}-\epsilon$ and $\alpha_{2}=1-\beta_{1}+\epsilon$; it is easy to check that $0<\alpha_{1}<\alpha_{2}<1$.  For $i=0$,  $\tau_{0}=1-l_{0}$ is in $[\alpha_{1},\alpha_{2}]$ by the estimate following~\eqref{neededestimateforsumoftau}.  For $i\geq1$ we may use~\eqref{inductionontheljvalues} to obtain
\begin{align*}
    \sum_{j=0}^{i} \tau_{j}
    &= \sum_{j=0}^{i+1} (i+1-j)\tau_{j} - \sum_{0}^{i}(i-j)\tau_{j}\\
    &= 1 - (i+1)l_{i} + il_{i-1}\\
    &=1-\frac{\log f(\scale^{i+1}) - \log f(\scale^{i})}{\log\scale}\\
    &\in[\alpha_{1},\alpha_{2}]
    \end{align*}
directly from~\eqref{neededestimateforsumoftau}.  We also need that $\sup_{j}|\tau_{j}|<1$.  As already noted, $\tau_{0}\in[\alpha_{1},\alpha_{2}]$.  For $j=1$ we write
\begin{equation*}
    \tau_{1}
    = l_{0} - (2l_{1}-l_{0})
    = l_{0} - \frac{\log f(\scale^{2}) - \log f(\scale)}{\log\scale}
    \end{equation*}
which is a difference of two values from $\bigl[\beta_{1}-\epsilon,\beta_{2}+\epsilon\bigr]$, so has magnitude bounded by $2\epsilon+\beta_{2}-\beta_{1}$. Similarly, for $j\geq2$ rewrite $\tau_{j}$ as
\begin{equation*}
    \Bigl( jl_{j-1} - (j-1)l_{j-2}\Bigr) - \Bigl( (j+1)l_{j}-jl_{j-1}\Bigr)
    =\frac{\log f(\scale^{j}) - \log f(\scale^{j-1})}{\log\scale}
    - \frac{\log f(\scale^{j+1}) - \log f(\scale^{j})}{\log\scale}
    \end{equation*}
which is a difference of the same type.  Since we chose that $3\epsilon<\beta_{2}-\beta_{1}$ we conclude that all $|\tau_{j}|$, $j\geq1$ are bounded by $\beta_{2}-\beta_{1}<1$.

We have verified that Lemma~\ref{SCfunctionhasrightests} may be applied to the sequence $\tau_{j}$ to produce a function $H(z)$. Since $\alpha_{2}>\alpha_{1}=1-\beta_{2}-\epsilon>0$ we may find a constant $c'$ such that $\sin\bigl(\beta+\alpha_{2}(\pi-\beta)\bigr)\leq c'\sin\bigl(\beta+\alpha_{1}(\pi-\beta)\bigr)$ for all $\beta\in[0,\pi]$; inserting this into the estimate of $\Im H(z)$ in Lemma~\ref{SCfunctionhasrightests} yields
\begin{equation*}
    |z|^{1-\sum_{j=0}^{K}(1-j/K)\tau_{j}} \sin\bigl(\beta+\alpha_{2}(\pi-\beta)\bigr)
    \lesssim  \Im H(|z|e^{i\beta})
    \lesssim |z|^{1-\sum_{j=0}^{K}(1-j/K)\tau_{j}} \sin\bigl(\beta+\alpha_{1}(\pi-\beta)\bigr)
    \end{equation*}
with $K$ the integer part of $\log|z|/\log\scale$.  However we determined in~\eqref{inductionontheljvalues} that
\begin{equation*}
    |z|^{1-\sum_{j=0}^{K}(1-j/K)\tau_{j}}
    =|z|^{l_{K-1}}
    =|z|^{\frac{\log f(\scale^{K})}{(K+1)\log\scale}}
    = \Bigl(f(\scale^{K})\Bigr)^{\frac{\log|z|}{(K+1)\log\scale}}.
    \end{equation*}
Moreover $\frac{\log|z|}{(K+1)\log\scale}\in\bigl[1-\frac{1}{K+1},1\bigr]$ and $f(\scale^{K})^{1/(K+1)}$ is bounded, as is the ratio of $f(\scale^{K})/f(|z|)$, so we have in fact
\begin{equation*}
    f(|z|) \sin\bigl(\beta+\alpha_{2}(\pi-\beta)\bigr)
    \lesssim  \Im H(|z|e^{i\beta})
    \lesssim f(|z|) \sin\bigl(\beta+\alpha_{1}(\pi-\beta)\bigr)
    \end{equation*}
which may be rewritten as
\begin{equation*}
    f(|z|) \sin\bigl((1-\alpha_{2})(\pi-\beta)\bigr)
    \lesssim  \Im H(|z|e^{i\beta})
    \lesssim f(|z|)  \sin\bigl((1-\alpha_{1})(\pi-\beta)\bigr).
    \end{equation*}
Hence there is some constant $c$ such that $F(z)=-iH(z)/c$ satisfies the conclusion of the lemma.
\end{proof}

%%%%%%%%%%%%%%%%%%%%%%%%%%%%%%%%%%%%%%%%%%%%%%%%%%%%%%%%%%%%%%%%%%%%%%%%%%%%%%%%%%%%%%%%%%%%%%%%%%%%%%%%%%%%%%

\section{Estimates away from the negative real axis.}\label{offnegaxissection}

In previous sections we have estimated piecewise eigenfunctions and the resolvent for real positive values of $\lambda$.  In this section we
combine them with a weak decay estimate for the resolvent and use the Phragmen-Lindel\"{o}f theorems to obtain decay estimates in a sector in $\mathbb{C}$. This gives bounds on the resolvent everywhere away from the negative real axis.

Our weak decay estimate for the resolvent relies on knowing $L^{\infty}$ bounds for projection kernels onto eigenspaces.  To obtain them we first prove a modified version of Theorem 4.5.4 of~\cite{Kigami2001}, which relates $L^{\infty}$ and $L^{2}$ norms for eigenfunctions.
\begin{theorem}\label{LinftyfromLtwoforeigenfns}
If $u$ is in the span of the eigenfunctions with eigenvalues not exceeding $\Lambda$ then $\|u\|_{\infty}\lesssim \Lambda^{S/2(S+1)} \|u\|_{2}$.
\end{theorem}
\begin{proof}
Let $u=\sum a_{j}\phi_{j}(x)$, where the $\phi_{j}$ are an orthonormal set of eigenfunctions with eigenvalues $-\lambda_{j}$, and $0\leq\lambda_{j}\leq\Lambda$ for all $j$.  The main point of our argument is that $u$ does not vary very much on any cell $F_{w}(X)$ for which $r_{w}\mu_{w}\Lambda\leq\frac{1}{2}$. Specifically, if $F_{w}(X)$ is a cell and $h_{w}$ is the harmonic function on $X$ with $h_{w}(p)=u(F_{w}(p))$ for all $p\in V_{0}$, then the difference between $u$ and $h_{w}$ on $F_{w}(X)$ can be obtained by integrating against the Green kernel $G^{(0)}(x,y)$.  Since we know $\|G^{(0)}(x,y)\|_{\infty}\lesssim1$ we have
\begin{align}
    \bigl\| u\circ F_{w} - h_{w}\circ F_{w} \bigr\|_{\infty}
    &= \Bigl\| \int_{X} G^{(0)}(x,y)\, \bigl(-\Delta (u\circ F_{w})(y)\bigr) \, d\mu(y) \Bigr\|_{\infty} \notag\\
    &\lesssim \bigl\| \Delta (u \circ F_{w}) \bigr\|_{1} \notag\\
    &= r_{w}\mu_{w} \bigl\| (\Delta u ) \circ F_{w} \bigr\|_{1} \notag\\
    &= r_{w} \bigl\| \Delta u  \bigr\|_{L^{1}(F_{w}(X))} \notag\\
    &\leq r_{w}\mu_{w}^{1/2} \bigl\| \Delta u  \bigr\|_{L^{2}(F_{w}(X))}. \label{eprojestone}
    \end{align}
If $\Theta$ is a partition of $X$ into cells with $r_{w}\mu_{w}\Lambda\leq\frac{1}{2}$ then
\begin{equation}
    \max_{w} r_{w}\mu_{w} \bigl\| \Delta u  \bigr\|_{2}
    = \max_{w} r_{w}\mu_{w} \Bigl(\sum_{j} a_{j}^{2}\lambda_{j}^{2} \Bigr)^{1/2}
    \leq \max_{w} r_{w}\mu_{w} \Lambda \|u\|_{2}
    \leq \frac{1}{2}\|u\|_{2}. \label{eprojesttwo}
    \end{equation}
Now if we make the above decomposition on each cell of $\Theta$ we obtain a piecewise harmonic function $h$ by setting $h\circ F_{w}=h_{w}\circ F_{w}$ and by~\eqref{eprojestone} and~\eqref{eprojesttwo}
\begin{equation}
    \|u-h\|_{\infty}
    \lesssim \max_{w} r_{w}\mu_{w}^{1/2} \bigl\| \Delta u  \bigr\|_{L^{2}(F_{w}(X))}
    \leq \max_{w} r_{w}\mu_{w}^{1/2} \bigl\| \Delta u  \bigr\|_{2}
    \leq \frac{1}{2} \max \mu_{w}^{-1/2} \|u\|_{2} \label{eprojestthree}
    \end{equation}
However~\eqref{eprojestone} and~\eqref{eprojesttwo} also yield
\begin{align}
    \|u-h\|_{2}^{2}
    =\sum_{w} \|u - h \|_{L^{2}(F_{w}(X))}^{2}
    &=\sum_{w} \mu_{w} \|u\circ F_{w} - h_{w}\circ F_{w} \|_{L^{2}(X)}^{2} \notag\\
    &\leq \sum_{w} \mu_{w} \bigl\| u\circ F_{w} - h_{w}\circ F_{w} \bigr\|_{\infty}^{2} \notag\\
    &\leq \sum_{w} r_{w}^{2}\mu_{w}^{2} \bigl\| \Delta u  \bigr\|_{L^{2}(F_{w}(X))}^{2} \notag\\
    &\leq \max_{w} r_{w}^{2}\mu_{w}^{2}  \bigl\| \Delta u  \bigr\|_{2}^{2}
    \leq \frac{1}{4}\|u\|_{2}^{2}. \label{eprojestfour}
    \end{align}
Moreover the $L^{\infty}$ and $L^{2}$ norms of $h_{w}\circ F_{w}$ are comparable because the space of harmonic functions is finite dimensional.  As in Lemma~4.5.5 of~\cite{Kigami2001} we conclude
\begin{equation*}
    \|h\|_{\infty}
    \leq \sum_{w} \|h_{w}\circ F_{w}\|_{\infty}
    \lesssim \sum_{w} \|h_{w}\circ F_{w}\|_{2}
    = \sum_{w} \mu_{w}^{-1/2} \|h_{w}\|_{L^{2}(F_{w}(X))}
    = \max_{w} \mu_{w}^{-1/2} \|h\|_{2}
    \end{equation*}
and by~\eqref{eprojestfour}, $\|h\|_{\infty}\leq \frac{3}{2}\max_{w} \mu_{w}^{-1/2} \|u\|_{2}$.  Combining this with~\eqref{eprojestthree} gives
\begin{equation*}
    \|u\|_{\infty}
    \leq \|h\|_{\infty} + \frac{1}{2} \max \mu_{w}^{-1/2} \|u\|_{2}
    \leq 2\max_{w} \mu_{w}^{-1/2} \|u\|_{2}.
    \end{equation*}
We required only $r_{w}\mu_{w}\Lambda\leq\frac{1}{2}$, so may choose it such that $\mu_{w}^{(S+1)/S}=r_{w}\mu_{w}\gtrsim\Lambda^{-1}$, from which $\mu_{w}^{-1/2}\lesssim \Lambda^{S/2(S+1)}$.
\end{proof}

\begin{theorem}\label{multiplierestimates}
If $\{\phi_{j}\}$ is a finite orthonormal set of eigenfunctions with eigenvalues not exceeding $\Lambda$ and $\{b_{j}\}$ is a set of complex numbers with all $|b_{j}|\leq B$ then
\begin{equation*}
    \Bigl\|\sum_{j} b_{j}\phi_{j}(x)\phi_{j}(y)\Bigr\|_{L^{\infty}(X\times X)}
    \lesssim B \Lambda^{S/(S+1)}.
    \end{equation*}
\end{theorem}
\begin{proof}
Fix $y\in X$ and take the $L^{\infty}$ norm with respect to $x$.  We see that $\sum_{j} \phi_{j}(x)\phi_{j}(y)$ is in the span of the eigenfunctions with eigenvalues at most $\Lambda$, so by Theorem~\ref{LinftyfromLtwoforeigenfns} it suffices to compute the $L^{2}$ norm with respect to $x$.  Since the $\phi_{j}(x)$ are orthonormal and appear with coefficients $\phi_{j}(y)$, we have
\begin{equation*}
    \Lambda^{-S/(S+1)}\Bigl\| \sum_{j} \phi_{j}(x)\phi_{j}(y)  \Bigr\|_{L^{\infty}(x)}^{2}
    \lesssim  \sum_{j} \phi_{j}(y)^{2}
    \leq  \Bigl\| \sum_{j}  \phi_{j}(y)^{2}  \Bigr\|_{\infty}
    \leq \Bigl\| \sum_{j}  \phi_{j}(x)\phi_{j}(y) \Bigr\|_{L^{\infty}(X\times X)}
    \end{equation*}
Taking the supremum over $y$ gives
\begin{equation*}
    \Bigl\| \sum_{j} \phi_{j}(x)\phi_{j}(y) \Bigr\|_{L^{\infty}(X\times X)}
    \lesssim \Lambda^{S/(S+1)}
    \end{equation*}
because the functions are continuous and $X$ is compact.  However the the fact that the $\phi_{j}$ are real-valued then ensures
\begin{equation*}
    \Bigl|\sum_{j} b_{j}\phi_{j}(x)\phi_{j}(y)\Bigr|
    \leq B \sum_{j} |\phi_{j}(x)\phi_{j}(y)|
    \leq B\Bigl( \sum_{j} \phi_{j}(x)^{2} \Bigr)^{1/2} \Bigl(\sum_{j}\phi_{j}(y)^{2} \Bigr)^{1/2}
    \lesssim B \Lambda^{S/(S+1)}.
    \end{equation*}
\end{proof}

\begin{remark}
The estimates in Theorems~\ref{LinftyfromLtwoforeigenfns} and~\ref{multiplierestimates} are sharp in the sense that there are fractals on which they are achieved. Indeed, on certain fractals with sufficient symmetry it is also known that there are individual eigenfunctions of eigenvalue $\lambda$ with support in cells of size $\lambda^{-S/(S+1)}$ when $\lambda$ is large.  For these to have $L^{2}$ norm $1$ they must then have $L^{\infty}$ size $\lambda^{S/2(S+1)}$, and the product $\phi_{j}(x)\phi_{j}(y)$ then has size $\lambda^{S/(S+1)}$ (see Theorem~4.5.4 of~\cite{Kigami2001}).  We do not know whether they are sharp on all pcfss sets, though it is not difficult to get a lower bound from the Weyl law in Proposition~\ref{KigamiLapidusWeylest}.  According to that result, the number of eigenfunctions with eigenvalue at most $\Lambda$ is bounded below by a multiple of $\Lambda^{S/(S+1)}$ if $\Lambda$ is sufficiently large.  Calling them $\phi_{j}$ we see $\sum_{j}\phi_{j}(x)\phi_{j}(y)$ has $L^{2}(X\times X)$ norm bounded below by a multiple of $\Lambda^{S/2(S+1)}$, and hence $L^{\infty}(X\times X)$ norm that is at least this great.  Hence as a converse to Theorem~\ref{multiplierestimates} we have the existence of a set of eigenfunctions satisfying the assumptions and for which there is $x$ at which $\sum_{j}\phi_{j}(x)\phi_{j}(x)\gtrsim\Lambda^{S/2(S+1)}$.  The consequence for Theorem~\ref{LinftyfromLtwoforeigenfns} is that setting $a_{j}=\phi_{j}(x)$ we have the $L^{\infty}$ norm of $u=\sum_{j}a_{j}\phi_{j}(x)$ is equal $\sum_{j}a_{j}^{2}\gtrsim\Lambda^{S/2(S+1)}$, while its $L^{2}$ norm is $\bigl(\sum_{j}a_{j}^{2}\bigr)^{1/2}$, so $u$ satisfies the hypotheses and has $\|u\|_{\infty}\geq\Lambda^{S/4(S+1)}\|u\|_{2}$.
\end{remark}

\begin{lemma}\label{etaradialgrowthest}
For $z\in\Sect$, the piecewise eigenfunction $\eta_{p}^{(z)}$ satisfies
\begin{equation*}
    \bigl\| \eta_{p}^{(z)}(x) \bigr\|_{L^{\infty}(X)}
    \lesssim \Bigl( 1+  \tan \bigl(\frac{\alpha}{2}\bigr) \Bigr).
    \end{equation*}
\end{lemma}
\begin{proof}
Recall $\zeta_{p}$ is the harmonic function on $X$ with value $1$ at $p\in V_{0}$ and $0$ at the other points of $V_{0}$. Let  $\{\phi_{j}\}$ be an orthonormal eigenfunction basis of $L^{2}$ with corresponding eigenvalues $\{-\lambda_{j}\}$ repeated according to multiplicity. It is easy to verify (see Lemma~3.3 of~\cite{IoPeRoHuSt2010TAMS}) that if the $L^{2}$ expansion of $\zeta_{p}$ in this basis is $\zeta_{p}=\sum_{j}b_{j}\phi_{j}$ then we have $\zeta_{p}-\eta_{p}^{(z)}=z\sum_{\lambda}\frac{b_{j}}{z+\lambda_{j}}\phi_{j}$, with $L^{2}$ convergence of the function and its Laplacian.  We also have the $L^{2}$ expansion $\eta_{p}^{(z)}=\sum_{\lambda}\frac{\lambda_{j}b_{j}}{z+\lambda_{j}}\phi_{j}$.
Now we compute at $z\in\Sect$
\begin{equation}\label{etazminusetamodz}
    \eta_{p}^{(z)}(x) - \eta_{p}^{(|z|)}(x)
    = \bigl(|z|-z\bigr) \sum_{j} \frac{\lambda_{j}b_{j}}{(z+\lambda_{j})(|z|+\lambda_{j})}\phi_{j}(x).
    \end{equation}
We break up the sum into pieces of the form $\lambda_{j}\in[2^{k},2^{k+1})$ and estimate each piece using Theorem~\ref{LinftyfromLtwoforeigenfns} to obtain
\begin{align*}
    \bigl\|\eta_{p}^{(z)}(x) - \eta_{p}^{(|z|)}(x)\|_{\infty}
    &\lesssim \bigl||z|-z\bigr| \sum_{k} \Bigl\| \sum_{2^{k}\leq\lambda_{j}<2^{k+1}} \frac{\lambda_{j}b_{j}}{(z+\lambda_{j})(|z|+\lambda_{j})}\phi_{j}(x) \Bigr\|_{\infty}\\
    &\leq \bigl||z|-z\bigr| \sum_{k} 2^{(k+1)S/2(S+1)} \Bigl\| \sum_{2^{k}\leq\lambda_{j}<2^{k+1}} \frac{\lambda_{j}b_{j}}{(z+\lambda_{j})(|z|+\lambda_{j})}\phi_{j}(x) \Bigr\|_{2}\\
    &= \bigl||z|-z\bigr| \sum_{k} 2^{(k+1)S/2(S+1)} \biggl( \sum_{2^{k}\leq\lambda_{j}<2^{k+1}} \frac{|\lambda_{j}b_{j}|^{2}}{\bigl|z+\lambda_{j}\bigr|^{2} \bigl||z|+\lambda_{j}\bigr|^{2}} \biggr)^{1/2}\\
    &\leq \sum_{k} \biggl( 2^{kS/2(S+1)} \max_{2^{k}\leq\lambda_{j}<2^{k+1}} \frac{\bigl||z|-z\bigr|}{|z+\lambda_{j}|} \biggr)
        \biggl( \sum_{2^{k}\leq\lambda_{j}<2^{k+1}} \frac{|\lambda_{j}b_{j}|^{2}}{\bigl||z|+\lambda_{j}\bigr|^{2}} \biggr)^{1/2}.
    \end{align*}
Now for $k\geq k_{0}=\log_{2}|z|+1$ we have $\lambda_{j}\geq2|z|$ so $|z+\lambda_{j}|\geq\frac{\lambda_{j}}{2}\geq2^{k-1}$ and the maximum in the above is at most $8|z|2^{-k}$.  For $k\leq k_{0}$ we instead bound the maximum using the fact that $z\in\Sect$. Clearly we have $|z+\lambda|\geq |z|\sin\alpha$, but this is a poor estimate when $\alpha<\frac{\pi}{2}$, where $|z+\lambda|\geq |z|$.  Instead observe that $\sin\alpha>\cos\frac{\alpha}{2}$ for $\alpha\in[\frac{\pi}{2},\pi]$ and so $|z+\lambda_{j}|\geq |z|\cos\frac{\alpha}{2}$ for all $\alpha\in[0,\pi]$.  Also $||z|-z|\leq2|z|\sin\frac{\alpha}{2}$, so the maximum is bounded by $2\tan \frac{\alpha}{2}$.  Inserting these bounds for the specified $k$ and using the Cauchy-Schwarz inequality gives
\begin{align*}
    \lefteqn{ \bigl\|\eta_{p}^{(z)}(x) - \eta_{p}^{(|z|)}(x)\|_{\infty}}\quad &\\
    &\leq 2\tan \Bigl( \frac{\alpha}{2}\Bigr) \sum_{k\leq k_{0}} 2^{kS/2(S+1)} \biggl( \sum_{2^{k}\leq\lambda_{j}<2^{k+1}}
        \frac{|\lambda_{j}b_{j}|^{2}}{\bigl||z|+\lambda_{j}\bigr|^{2}} \biggr)^{1/2} + 8|z| \sum_{k\geq k_{0}} 2^{kS/2(S+1)}2^{-k} \biggl( \sum_{2^{k}\leq\lambda_{j}<2^{k+1}} \frac{|\lambda_{j}b_{j}|^{2}}{\bigl||z|+\lambda_{j}\bigr|^{2}} \biggr)^{1/2}\\
    &\leq 2\tan \Bigl( \frac{\alpha}{2}\Bigr) \Bigl( \sum_{k\leq k_{0}}2^{kS/(S+1)}\Bigr)^{1/2} \biggl(
        \sum_{\lambda_{j}\leq 2^{k_{0}+1}} \frac{|\lambda_{j}b_{j}|^{2}}{\bigl||z|+\lambda_{j}\bigr|^{2}} \biggr)^{1/2}
        + 8|z| \Bigl( \sum_{k\geq k_{0}} 2^{kS/(S+1)}2^{-2k}\Bigr)^{1/2} \biggl( \sum_{\lambda_{j}\geq2^{k_{0}}} \frac{|\lambda_{j}b_{j}|^{2}}{\bigl||z|+\lambda_{j}\bigr|^{2}} \biggr)^{1/2}\\
    &\lesssim \Bigl( 2\tan \Bigl( \frac{\alpha}{2}\Bigr) 2^{k_{0}S/2(S+1)} + 8|z|2^{k_{0}S/2(S+1)}2^{-k_{0}} \Bigr) \Bigl( \sum_{j} \frac{|\lambda_{j}b_{j}|^{2}}{\bigl||z|+\lambda_{j}\bigr|^{2}} \Bigr)^{1/2}\\
    &\leq \Bigl( 4\tan \Bigl( \frac{\alpha}{2}\Bigr) |z|^{S/2(S+1)} + 16|z|^{S/2(S+1)} \Bigr)\bigl\| \eta_{p}^{(|z|)} \bigr\|_{2}
    \end{align*}
where in the final step we used that $k_{0}=\log_{2}|z|+1$ and the known $L^{2}$ expansion of $\eta_{p}^{(|z|)}$.  The latter can be estimated using~\eqref{finalboundsforeta} from Theorem~\ref{maintheoremoneta}.  This shows $\eta_{p}^{(|z|)}$ has exponential decay with scale $k(|z|)$, so that its $L^{2}$ norm is dominated by integrating over the cell containing $p$ and of measure $|z|^{-S/(S+1)}$.  Since the function is bounded by $1$ we find $\bigl\| \eta_{p}^{(|z|)} \bigr\|_{2}\lesssim |z|^{-S/2(S+1)}$.  Again using~\eqref{finalboundsforeta}, this time as an $L^{\infty}$ bound on $\eta_{p}^{(|z|)}$, completes the proof.
\end{proof}

\begin{lemma}\label{Gradialgrowthest} For $z\in\Sect$, the Green kernel for a pcfss fractal with regular harmonic structure satisfies
\begin{equation*}
    \bigl\| G^{(z)}(x,y) \bigr\|_{L^{\infty}(X\times X)}
    \lesssim (1+ |z|)^{-1/(S+1)} \Bigl( 1+  \tan \bigl(\frac{\alpha}{2}\bigr) \Bigr).
    \end{equation*}
\end{lemma}
\begin{proof}
Expanding the resolvent kernel at both $z\in\Sect$ and at $|z|\in\mathbb{R}$ with respect to the basis $\{\phi_{j}\}$ of eigenfunctions from the previous theorem we find
\begin{equation*}
    G^{(z)}(x,y) - G^{(|z|)}(x,y)
    =\sum_{j} \Bigl( \frac{1}{z+\lambda_{j}} - \frac{1}{|z|+\lambda_{j}} \Bigr)\phi_{j}(x)\phi_{j}(y)
    =\sum_{j} \frac{|z|-z}{(z+\lambda_{j})(|z|+\lambda_{j})} \phi_{j}(x)\phi_{j}(y).
    \end{equation*}
As before we break up the sum, taking one piece where $\lambda_{j}\leq4|z|$, and the rest to be of the form $\lambda_{j}\in[2^{k},2^{k+1})$, $k\geq k_{0}=\log_{2}|z|+1$.  For each of the the latter we have $|z|+\lambda_{j}\geq|z+\lambda_{j}|\geq \frac{\lambda_j}{2}\geq 2^{k-1}$ so the coefficients multiplying $\phi_{j}(x)\phi_{j}(y)$ are bounded by $8|z|2^{-2k}$.  Applying Theorem~\ref{multiplierestimates} we find
\begin{equation}\label{largelambdapieceofresolventweakest}
    \Bigl\| \sum_{2^{k}|z|\leq\lambda_{j}\leq 2^{k+1}|z|} \frac{|z|-z}{(z+\lambda_{j})(|z|+\lambda_{j})} \phi_{j}(x)\phi_{j}(y) \Bigr\|_{\infty}
    \leq 8|z| 2^{-2k} 2^{(k+1)S/(S+1)}
    \end{equation}
For the piece where $\lambda_{j}\leq 4|z|$ we use $|z+\lambda_{j}|\geq |z|\cos\frac{\alpha}{2}$ and $||z|-z|\leq2|z|\sin\frac{\alpha}{2}$ for all $\alpha\in[0,\pi]$, because $z\in\Sect$. The coefficients are then bounded by $2|z|^{-1}\tan \frac{\alpha}{2}$, and applying Theorem~\ref{multiplierestimates} gives
\begin{equation*}
    \Bigl\| \sum_{\lambda_{j}\leq 4|z|} \frac{|z|-z}{(z+\lambda_{j})(|z|+\lambda_{j})} \phi_{j}(x)\phi_{j}(y) \Bigr\|_{\infty}
    \lesssim 2|z|^{-1}\tan \bigl(\frac{\alpha}{2}\bigr)  |z|^{S/(S+1)}
    =  2 |z|^{-1/(S+1)} \tan \bigl(\frac{\alpha}{2}\bigr).
    \end{equation*}
Combining this with~\eqref{largelambdapieceofresolventweakest} we have
\begin{align*}
    \bigl\|G^{(z)}(x,y) - G^{(|z|)}(x,y)\bigr\|_{\infty}
    &\leq 2 |z|^{-1/(S+1)} \tan \bigl(\frac{\alpha}{2}\bigr)  + \sum_{k\geq k_{0}} \Bigl\| \sum_{2^{k}\leq\lambda_{j}\leq 2^{k+1}} \frac{|z|-z}{(z+\lambda_{j})(|z|+\lambda_{j})} \phi_{j}(x)\phi_{j}(y) \Bigr\|_{\infty}\\
    &\leq  2  |z|^{-1/(S+1)} \tan \bigl(\frac{\alpha}{2}\bigr)  +  8|z| \sum_{k\geq k_{0}} 2^{-2k} 2^{(k+1)S/(S+1)}\\
    &\lesssim  |z|^{-1/(S+1)} \tan \bigl(\frac{\alpha}{2}\bigr) + 8|z|  2^{-2k_{0}} 2^{(k_{0}+1)S/(S+1)}\\
    &\lesssim  |z|^{-1/(S+1)}\Bigl( 1+  \tan \bigl(\frac{\alpha}{2}\bigr) \Bigr)
    \end{align*}
However the bound can be improved for small $z$, because if $2|z|<\lambda_{0}$, the smallest Dirichlet eigenvalue, then the first sum is empty and the second begins at $\log\lambda_{0}$ rather than $\log|z|$, giving a bound independent of $z$.  We may therefore replace $|z|^{-1/(S+1)}$ with $(1+|z|)^{-1/(S+1)}$ in the estimate.  The result now follows by using~\eqref{pathdecayforGonXthmestimateofG} from Theorem~\ref{pathdecayforGonXthm} to bound $G^{(|z|)}(x,y)$.
\end{proof}

\begin{remark}
A small modification of the method in the previous proof provides bounds off suitably small neighborhoods of the eigenvalues.    If $D_{j}$ is the disc radius  $\delta \lambda_{j}^{1/(S+1)}$ centered at $-\lambda_{j}$ then outside $\cup_{j}D_{j}$ we have $|z+\lambda_{j}|\geq \delta^{-1}\lambda_{j}^{1/(S+1)}$ and using this for those $\lambda_{j}\leq 4|z|$ rather than the estimate in the above proof we see that
\begin{equation*}
    \bigl\|G^{(z)}(x,y)\bigr\|_{\infty}
    \lesssim \delta^{-1} |z|^{(S-1)/(S+1)} \text{ off }\cup_{j}D_{j}
\end{equation*}
Moreover the Weyl estimate in Proposition~\ref{KigamiLapidusWeylest} guarantees that the number of eigenvalues of size at most $\Lambda$ is bounded by a constant multiple of $\Lambda^{S/(S+1)}$ and therefore a suitably small $\delta$ depending only on the fractal and harmonic structure ensures $\cup_{j}D_{j}$ contains at most half of any interval of the form $[-2^{k+1},-2^{k}]$.  It is easy to check that this estimate and that in the lemma are comparable at points of the form $-\lambda_{j}+i\delta\lambda_{j}^{1/(S+1)}$.
\end{remark}

We complete our estimates in $\Sect$ with some for the normal derivatives.
\begin{lemma}
Under the assumptions of Lemmas~\ref{etaradialgrowthest} and~\ref{Gradialgrowthest},
\begin{align}
    \bigl| \partial_{n}\eta_{p}^{(z)}(q)  \bigr|
    &\lesssim (1+ |z|)^{1/(S+1)} \Bigl( 1+  \tan \bigl(\frac{\alpha}{2}\bigr) \Bigr) \label{partialnetaweakbounds}\\
    \bigl\| \partial_{n}'G^{(z)}(p,x) \bigr\|_{L^{\infty}(X)}
    &\lesssim \Bigl( 1+  \tan \bigl(\frac{\alpha}{2}\bigr) \Bigr) \label{partialnGweakbounds}\\
    \bigl| \partial_{n}''\partial_{n}'G^{(z)}(p,q) \bigr|
    &\lesssim (1+ |z|)^{1/(S+1)} \Bigl( 1+  \tan \bigl(\frac{\alpha}{2}\bigr) \Bigr). \label{partialnpartialnGweakbounds}
    \end{align}
\end{lemma}
\begin{proof}
Applying the Laplacian to
\begin{equation*}
    G^{(z)}(x,y) - G^{(|z|)}(x,y)
    =\sum_{j} \frac{|z|-z}{(z+\lambda_{j})(|z|+\lambda_{j})} \phi_{j}(x)\phi_{j}(y).
    \end{equation*}
gives an $L^{2}$ convergent series.  We may therefore expand in the Gauss-Green formula
\begin{align*}
    \lefteqn{\partial_{n}''G^{(z)}(x,p) - \partial_{n}''G^{(|z|)}(x,p)}\quad&\\
    &= \sum_{q\in V_{0}} \partial_{n}''\bigl(G^{(z)}(x,q) - G^{(|z|)}(x,q)\bigr) \zeta_{p}(q) - \bigl(G^{(z)}(x,q) - G^{(|z|)}(x,q)\bigr)\partial_{n}\zeta_{p}(q)\\
    &= \int \Delta\bigl(G^{(z)}(x,y) - G^{(|z|)}(x,y)\bigr) \zeta_{p}(y)\, d\mu(y)\\
    &= \int \sum_{j} \frac{(|z|-z)\lambda_{j}}{(z+\lambda_{j})(|z|+\lambda_{j})} \phi_{j}(x)\phi_{j}(y)\zeta_{p}(y)\, d\mu(y)\\
    &= \bigl(|z|-z\bigr)\sum_{j} \frac{\lambda_{j}b_{j}}{(z+\lambda_{j})(|z|+\lambda_{j})} \phi_{j}(x)\\
    &= \eta_{p}^{(z)}(x) - \eta_{p}^{(|z|)}(x)
    \end{align*}
as we saw in~\eqref{etazminusetamodz}.  The justification for exchanging the integral and sum is that  $\|\phi_{j}\|_{\infty}\leq \lambda_{j}^{S/2(S+1)}$ by Theorem~\ref{LinftyfromLtwoforeigenfns} and $\zeta_{p}$ is in $L^{2}$, so the series form of the integrand is $L^{1}$ convergent uniformly in $x$.  The same is true of the series expansion of $\eta_{p}^{(z)}(x) - \eta_{p}^{(|z|)}(x)$ by Lemma~\ref{etaradialgrowthest}, and combining the estimate of that lemma with~\eqref{pathdecayforGonXthmestimateofpartialnG} from Theorem~\ref{pathdecayforGonXthm} gives~\eqref{partialnGweakbounds}.

Similarly, the fact that applying the Laplacian to $\eta_{p}^{(z)}(x) - \eta_{p}^{(|z|)}(x)$ gives an $L^{2}$ convergent series allows us to compute
\begin{align*}
    \partial_{n}\eta_{p}^{(z)}(p) - \partial_{n}\eta_{p}^{(|z|)}(p)
    &= \sum_{q\in V_{0}} \partial_{n}\bigl(\eta_{p}^{(z)} - \eta_{p}^{(|z|)}\bigr)(q) \zeta_{p}(q) - \bigl(\eta_{p}^{(z)} - \eta_{p}^{(|z|)}\bigr)(q)\partial_{n}\zeta_{p}(q)\\
    &= \int \Delta \bigl(\eta_{p}^{(z)}(x) - \eta_{p}^{(|z|)}(x)\bigr) \zeta_{p}(x) \, d\mu(x)\\
    &= \int  \bigl(|z|-z\bigr)\sum_{j} \frac{\lambda_{j}^{2} b_{j}}{(z+\lambda_{j})(|z|+\lambda_{j})} \phi_{j}(x) \zeta_{p}(x)\, d\mu(x)\\
    &=  \bigl(|z|-z\bigr)\sum_{j} \frac{\lambda_{j}^{2} b_{j}^{2}}{(z+\lambda_{j})(|z|+\lambda_{j})}\\
    &\leq \bigl(|z|-z\bigr)\sum_{j} \frac{|\lambda_{j} b_{j}|^{2}}{(|z|+\lambda_{j})^{2}}\\
    &\leq 2|z| \bigl\| \eta_{p}^{(|z|)}\bigr\|_{2}^{2}
    \lesssim |z|^{1/(S+1)}
    \end{align*}
where we used the decay estimate~\eqref{finalboundsforeta} as we did in Lemma~\ref{etaradialgrowthest}. This and the normal derivative estimates from Theorem~\ref{maintheoremoneta} give~\eqref{partialnetaweakbounds} and together with~\eqref{pathdecayforGonXthmestimateofpartialnpartialnnGoffdiag} from Theorem~\ref{pathdecayforGonXthm} establishes~\eqref{partialnpartialnGweakbounds}.
\end{proof}

%%%%%%%%%%%%%%%%%%%%%%%%%%%%%

With the preceding weak decay estimates in hand we may prove the main result of this section, which gives bounds on both the functions $\eta_{p}^{(z)}$ and the resolvent kernel $G^{(z)}(x,y)$, as well as their normal derivatives, at points $z\in\mathbb{C}$ that are not on the negative real axis.
\begin{theorem}\label{Allestimatesoffnegaxis}
Suppose $X$ is a post-critically finite fractal with regular harmonic structure.  Then there is a constant $\kappa_{5}>0$ depending only on the fractal and harmonic structure such that for $z=|z|e^{i\beta}$, $|\beta|<\pi$ we have
\begin{align*}
    \bigl| \eta_{p}^{(z)}(x) \bigr|
    &\lesssim \Phi(\beta,|z|,p,x)\\
    \bigl| \partial_{n} \eta_{p}^{(z)}(q) \bigr|
    &\lesssim |z+1|^{1/(S+1)} \Phi(\beta,|z|,p,q)\\
    \bigl| G^{(z)}(x,y) \bigr|
    &\lesssim |z+1|^{-1/(S+1)} \Phi(\beta,|z|,x,y) \\
    \bigl| \partial_{n}'G^{(z)}(p,y) \bigr|
    &\lesssim \Phi(\beta,|z|,p,y) \\
    \bigl| \partial_{n}''\partial_{n}'G^{(z)}(p,q) \bigr|
    &\lesssim |z+1|^{1/(S+1)} \Phi(\beta,|z|,p,q)
    \end{align*}
where
\begin{equation*}
    \Phi(\beta,|z|,x,y)
    = \Bigl( 1+  \tan \bigl(\frac{\pi+|\beta|}{4}\bigr) \Bigr) \exp \biggl(- \kappa_{5} d_{k(|z|)}(x,y) \sin\Bigl(\frac{\tilde{c}(\pi-|\beta|)}{\pi+|\beta|}\Bigr)\biggr)
    \end{equation*}
\end{theorem}
\begin{proof}
We give only the proof for $G^{(z)}(x,y)$ because the others are entirely analogous.  Fix $x$ and $y$ in $X$.  Recall from~\eqref{pathdecayforGonXthmestimateofG} that the resolvent satisfies
\begin{equation}\label{sectorestimates_Gupperbound}
    G^{(\lambda)}(x,y)
    \lesssim (1+\lambda)^{-1/(S+1)} \exp\Bigl( - \kappa_{2} d_{k(\lambda)}\bigl( x,y \bigr) \Bigr).
    \end{equation}
for $\lambda\in(0,\infty)$.

Let $f(\lambda)= \kappa_{2} d_{k(\lambda^{\alpha/\pi})}(x,y)$.  Recalling Definition~\ref{defnofklambda} we may substitute into Lemma~\ref{boundonratioofchemicaldistancesfordifferentscales} with $k'=\frac{\alpha}{\pi}\log\scale$ and $k=jk'$ to find that there are constants $c_{1}$ and $c_{2}$, depending only on the harmonic structure, for which
\begin{equation*}
    c_{1} \scale^{\frac{\alpha}{\pi(S+1)}} \leq \frac{d_{k(\scale^{ (j+1)\alpha/\pi})}(x,y)}{ d_{k(\scale^{j\alpha/\pi})}(x,y)}
    \leq c_{2} \scale^{\frac{\alpha}{2\pi}}
    \end{equation*}
provided $M$ is large enough. (Note that the latter restriction is so $k(\scale^{(j+1)\alpha/\pi})$ is given by the formula in Definition~\ref{defnofklambda} and is not zero.) Thus~\eqref{hypothesisforSCfunctionconstruction} holds for $f(\lambda)$ with $0<\beta_{1}=\frac{\alpha}{\pi(S+1)}<\beta_{2}=\frac{\alpha}{2\pi}<1$.  Consider the product
\begin{equation*}
    g(z)
    = \frac{(z+1)^{1/(S+1)}}{C \bigl( 1+\tan\frac{\alpha}{2}\bigr)}
    G^{(z)}(x,y)
    \end{equation*}
on the sector $\Sect^{+}$.  Lemma~\ref{Gradialgrowthest} implies $|g|$ is bounded by $1$ on $\Sect^{+}$ if $C$ is large enough, while~\eqref{sectorestimates_Gupperbound} implies $|g(z)|\leq \exp\bigl(-f(|z|^{\pi/\alpha})\bigr)$ on the positive real axis.  Applying Theorem~\ref{generalpcfssversionofPLcorol} and multiplying out to retrieve $G^{(z)}$ we have for $\beta\in[0,\alpha]$ that
\begin{equation*}
    \bigl|G^{(|z|e^{i\beta})}(x,y) \bigr|
    \lesssim |z+1|^{-1/(S+1)} \Bigl( 1+  \tan \bigl(\frac{\alpha}{2}\bigr) \Bigr) \exp \biggl(-c \kappa_{2} d_{k(|z|)}(x,y) \sin\tilde{c}\Bigl(1-\frac{\beta}{\alpha}\Bigr)  \biggr).
    \end{equation*}
In order to obtain an estimate on a general ray $z=|z|e^{i\beta}$, $\beta\neq\pi$ it then suffices to take $\alpha=\frac{1}{2}(\pi+\beta)$, so
\begin{equation*}
    \bigl|G^{(|z|e^{i\beta})}(x,y) \bigr|
    \lesssim |z+1|^{-1/(S+1)}\Bigl( 1+  \tan \bigl(\frac{\pi+\beta}{4}\bigr) \Bigr) \exp \biggl(-c \kappa_{2} d_{k(|z|)}(x,y) \sin\Bigl(\frac{\tilde{c}(\pi-\beta)}{\pi+\beta}\Bigr)\biggr)
    \end{equation*}
and by symmetry corresponding estimates are valid on the sector $\Sect^{-}$ obtained by reflection in the real axis, simply by replacing $\beta$ with $|\beta|$.
\end{proof}

Note that in the case of an affine nested fractal we could have proved this result using the  classical Phragmen-Lindel\"{o}f theorem, because $d_{k(\lambda)}(x,y)\simeq (1+\lambda)^{\gamma/(S+1)} R(x,y)^{\gamma}$ by Proposition~\ref{affinenesteedversionofdk} and Definition~\ref{defnofklambda}.  In this case one also obtains better constants.

\begin{theorem}\label{mainblowuptheorem}
The sequence $G_{-n}^{(z)}(x,y)$ defined in Theorem~\ref{pathdecayforGonblowupthm} for a blow-up $\Omega$ of $X$ converges uniformly on compact subsets of $\Omega\times\Omega\times\bigl(\mathbb{C}\setminus(-\infty,0]\bigr)$ to the Laplacian resolvent $G^{(\lambda)}_{\infty}(x,y)$, and there is $\kappa_{6}>0$ such that
\begin{equation*}
    \tilde{\Phi}(\beta,|z|,x,y)
    = \Bigl( 1+  \tan \bigl(\frac{\pi+|\beta|}{4}\bigr) \Bigr) \exp \biggl(- \kappa_{6} d_{k(|z|)}(x,y) \sin\Bigl(\frac{\tilde{c}(\pi-|\beta|)}{\pi+|\beta|}\Bigr)\biggr).
    \end{equation*}
\end{theorem}
\begin{proof}
Using the definition of $G^{(z)}_{-n}(x,y)$ and the bound from Theorem~\ref{Allestimatesoffnegaxis} we have
\begin{align*}
    \lefteqn{\bigl| G^{(z)}_{-n}(x,y) \bigr|}\quad&\\
    &=  r_{[\w]_{n}}^{-1} \bigl| G^{(r_{[\w]_{n}}^{-1}\mu_{[\w]_{n}}^{-1}z)}_{0}
        (F_{\w_{n}}\circ\dotsm\circ F_{\w_{1}}x,F_{\w_{n}}\circ\dotsm\circ F_{\w_{1}}y) \bigr| \\
    &\lesssim r_{[\w]_{n}}^{-1} \bigl| r_{[\w]_{n}}^{-1}\mu_{[\w]_{n}}^{-1}z +1 \bigr|^{-1/(S+1)}
        \Phi\bigl(\beta,r_{[\w]_{n}}^{-1}\mu_{[\w]_{n}}^{-1}|z|,F_{\w_{n}}\circ\dotsm\circ F_{\w_{1}}x,F_{\w_{n}}\circ\dotsm\circ F_{\w_{1}}y \bigr)\\
    &= \bigl| z + r_{[\w]_{n}}\mu_{[\w]_{n}}\bigr|^{-1/(S+1)} \Phi\bigl(\beta,r_{[\w]_{n}}^{-1}\mu_{[\w]_{n}}^{-1}|z|,F_{\w_{n}}\circ\dotsm\circ F_{\w_{1}}x,F_{\w_{n}}\circ\dotsm\circ F_{\w_{1}}y \bigr)
    \end{align*}
but we recall from~\eqref{distancecomparabilityremarkrestated} that
\begin{equation*}
    d_{k(r_{[\w]_{n}}^{-1}\mu_{[\w]_{n}}^{-1}\lambda)}\bigl( F_{\w_{n}}\circ\dotsm\circ F_{\w_{1}}x,F_{\w_{n}}\circ\dotsm\circ F_{\w_{1}}y \bigr)
    \simeq d_{k(\lambda)}\bigl( x,y \bigr)
    \end{equation*}
from which there is $\kappa_{6}>0$ such that
\begin{equation*}
    \bigl| G^{(z)}_{-n}(x,y)\bigr|
    \lesssim \bigl| z + r_{[\w]_{n}}\mu_{[\w]_{n}}\bigr|^{-1/(S+1)} \tilde{\Phi}(\beta,|z|,x,y)
    \end{equation*}
independent of $n$.  Thus the sequence is uniformly bounded on any compact set of the type in the hypotheses, and since it is analytic in $z$ it is a normal family by Montel's theorem.  The function $G^{(\lambda)}_{\infty}(x,y)$ is the unique limit point because the sequence converges on the positive real axis by Theorem~\ref{pathdecayforGonblowupthm}.
\end{proof}

%%%%%%%%%%%%%%%%%%%%%%%%%%%%%%%%%%%%

\section{Estimates for other kernels}\label{otherkernelsection}

One of the main purposes for proving estimates on the resolvent kernel is that we can obtain estimates of the kernels of other operators from it using functional calculus.  Specifically, if $\Gamma$ is a contour in $\mathbb{C}$ that surrounds the spectrum $\{-\lambda_{j}\}$ of $\Delta$ and $h(z)$ is analytic in a neighborhood of $\Gamma$ and its interior then we define
\begin{equation}\label{defnofH}
    H(x,y)
    = \frac{1}{2\pi i} \int_{\Gamma} G^{(z)} (x,y) h(z) \, dz
    \end{equation}
provided $h$ is such that the integral converges. Under reasonable assumptions this should be the kernel of $h(\Delta)$.  At times it is useful to think of the kernel of $h(\Delta)$ as $\tilde{H}(x,y)=\sum_{j}h(-\lambda_{j})\phi_{j}(x)\phi_{j}(y)$, so we begin by describing sufficient conditions for these to coincide.  It will be convenient for us to work with contours that lie in a sector $\Sect$ as in Section~\ref{plsection}.

\begin{lemma}
Suppose $\Gamma$ is a contour in $\Sect$ that surrounds the spectrum of $\Delta$, that $h(z)$ is analytic in a neighborhood of $\Gamma$ and its interior, and that there is some $a>\frac{S}{2(S+1)}$ such that
\begin{equation}\label{conditforhintegralsokernelssame}
    \int_{\Gamma} |h(z)| |z|^{a-1} \, |dz| <\infty.
    \end{equation}
If in addition $\sum_{j} |h(-\lambda_{j})|^{2}<\infty$ then both $H(x,y)$ and $\tilde{H}(x,y)$ converge in $L^{2}(X\times X)$ to a common limit.
If the stronger condition
\begin{equation}\label{sumabilityconditonhlambdasokernelssame}
    \sum_{k} 2^{kS/(S+1)}\sup_{2^{k}\leq\lambda_{j}<2^{k+1}} \bigl|h(-\lambda_{j})\bigr| <\infty
    \end{equation}
holds then both $H(x,y)$ and $\tilde{H}(x,y)$ converge in $L^{\infty}(X\times X)$ to a common limit.
\end{lemma}
\begin{proof}
The Cauchy formula guarantees that
\begin{equation*}
    \tilde{H}(x,y)
    = \sum_{j}h(-\lambda_{j})\phi_{j}(x)\phi_{j}(y)
    = \sum_{j} \frac{1}{2\pi i}\int_{\Gamma} \frac{h(z)}{z+\lambda_{j}}\, dz\,  \phi_{j}(x)\phi_{j}(y)
    \end{equation*}
so we must prove that we can exchange sum and integral to obtain
\begin{equation*}
    \frac{1}{2\pi i}\int_{\Gamma} \sum_{j} \frac{h(z)}{z+\lambda_{j}} \phi_{j}(x)\phi_{j}(y) \, dz
    =H(x,y)
    \end{equation*}
where the convergence is in $L^{p}(X\times X)$ for $p=2$ or $p=\infty$.  So it suffices to prove convergence of
\begin{gather}
    \Bigl\| \sum_{j} h(-\lambda_{j})\phi_{j}(x)\phi_{j}(y) \Bigr\|_{L^{p}(X\times X)} \text{ and,} \label{firstsummabilityforopkernel}\\
    \Bigl\| \int_{\Gamma} \sum_{j} \frac{h(z)}{z+\lambda_{j}} \phi_{j}(x)\phi_{j}(y) \, dz \Bigr\|_{L^{p}(X\times X)}. \label{secondsummabilityforopkernel}
    \end{gather}
For~\eqref{firstsummabilityforopkernel} we have
\begin{align*}
    \Bigl\| \sum_{j} h(-\lambda_{j})\phi_{j}(x)\phi_{j}(y) \Bigr\|_{L^{2}(X\times X)}
    &= \Bigl(\sum_{j} \bigl|h(-\lambda_{j})\bigr|^{2}\Bigr)^{1/2}\\
    \Bigl\| \sum_{j} h(-\lambda_{j})\phi_{j}(x)\phi_{j}(y) \Bigr\|_{L^{\infty}(X\times X)}
    &\lesssim \sum_{k} 2^{kS/(S+1)} \sup_{2^{k}\leq\lambda_{j}<2^{k+1}} \bigl|h(-\lambda_{j})\bigr|
    \end{align*}
where the latter is from Theorem~\ref{multiplierestimates}.  Summability of these are the conditions given in the hypotheses for $p=2$ and $p=\infty$ respectively.

For~\eqref{secondsummabilityforopkernel} we use the integral Minkowski inequality and Theorem~\ref{multiplierestimates} to obtain
\begin{align*}
    \lefteqn{\Bigl\| \int_{\Gamma} \sum_{j} \frac{h(z)}{z+\lambda_{j}} \phi_{j}(x)\phi_{j}(y) \, dz \Bigr\|_{L^{p}(X\times X)}}\quad&\\
    &\leq  \int_{\Gamma} |h(z)| \Bigl\| \sum_{j} \frac{1}{z+\lambda_{j}} \phi_{j}(x)\phi_{j}(y) \Bigr\|_{L^{p}(X\times X)} \, |dz|\\
    &\lesssim \begin{cases}
        \displaystyle \int_{\Gamma} |h(z)| \Bigl( \sum_{j} |z+\lambda_{j}|^{-2}\Bigr)^{1/2} \, |dz| &\quad p=2\\
        \displaystyle\int_{\Gamma} |h(z)| \sum_{k} 2^{kS/(S+1)} \sup_{2^{k}\leq\lambda_{j}<2^{k+1}} |z+\lambda_{j}|^{-1} \, |dz| &\quad p=\infty
        \end{cases}
    \end{align*}
But now $|z+\lambda_{j}|^{-1}\leq \tan\bigl(\frac{\alpha}{2}\bigr)\min\bigl\{|z|^{-1},\lambda_{j}^{-1}\bigr\}$.
For $a>\frac{S}{2(S+1)}$ and using the Weyl law of Proposition~\ref{KigamiLapidusWeylest}
\begin{equation*}
    \Bigl( \sum_{j} |z+\lambda_{j}|^{-2}\Bigr)^{1/2}
    \lesssim |z|^{a-1}  \tan\Bigl(\frac{\alpha}{2}\Bigr)  \Bigl(\sum_{k} 2^{kS/(S+1)} 2^{-2ak} \Bigr)^{1/2}
    \leq c_{a} |z|^{a-1}  \tan\Bigl(\frac{\alpha}{2}\Bigr)
    \end{equation*}
and also
\begin{equation*}
    \sum_{k} 2^{kS/(S+1)} \sup_{2^{k}\leq\lambda_{j}<2^{k+1}} |z+\lambda_{j}|^{-1}
    \leq |z|^{a-1} \tan\Bigl(\frac{\alpha}{2}\Bigr)  \Bigl(\sum_{k} 2^{kS/(S+1)} 2^{-2ak} \Bigr)^{1/2}
    \leq c_{a} |z|^{a-1}  \tan\Bigl(\frac{\alpha}{2}\Bigr)
    \end{equation*}
so that we have for both $p=2$ and $p=\infty$ and any $a>\frac{S}{2(S+1)}$ there is $c_{a}$ so that
\begin{equation*}
    \Bigl\| \int_{\Gamma} \sum_{j} \frac{h(z)}{z+\lambda_{j}} \phi_{j}(x)\phi_{j}(y) \, dz \Bigr\|_{L^{p}(X\times X)}
    \lesssim c_{a} \int_{\Gamma} |h(z)|  |z|^{a-1}\, |dz|
    \end{equation*}
and finiteness of this integral is enough to ensure~\eqref{secondsummabilityforopkernel}.
\end{proof}

We can then obtain $L^{\infty}$ bounds for the kernel $H(x,y)$ by inserting our estimates on $G^{(z)}(x,y)$ into~\eqref{defnofH} and integrating over $\Gamma$.  How best to do this depends very much on the function $h(z)$, so rather than attempt to formulate a general theorem we simply give some examples of how this can be used.  There are many other examples that could be treated in a similar manner to those given below, such as exponentials of complex powers of $\Delta$, some of which we may return to in future work.

\begin{example}[Heat Kernel]
The heat kernel is defined for $t>0$ by
\begin{equation*}
    p_{t}(x,y)= \sum_{j} e^{-\lambda_{j}t} \phi_{j}(x)\phi_{j}(y)
    \end{equation*}
so in the above notation it is $H(x,y)$ for $h(z)=e^{tz}$.  To estimate it we introduce a suitable family of contours.

Fix $t>0$ and $\alpha\in(0,\pi)$, and let $\Gamma_{\alpha,t}$ be as in Figure~\ref{Gammafigure}. It consists of the arc of the circle $t|z|=1$ in the sector $\Sect$, as well as the rays $z=|z|e^{\pm i\alpha}$, $|z|\in[t^{-1},\infty)$, traversed so as to wind once around any point on the
negative real axis.
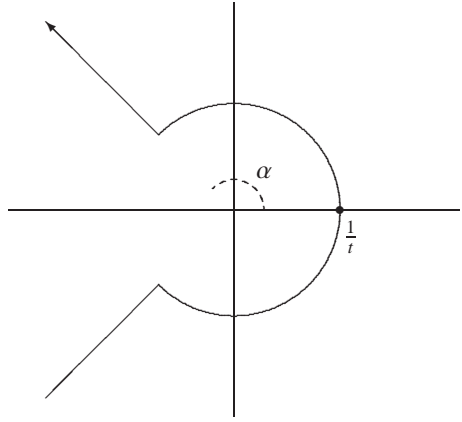
\begin{figure}
\begin{center}
\setlength{\unitlength}{1mm}
\begin{picture}(60,60)(-30,-30)
\put(0,-27.5){\line(0,1){55}}%
\put(-30,0){\line(1,0){60}}%
\put(0,0){\arc(-10,-10){270}}%
\put(-10,-10){\line(-1,-1){15}}%
\put(-10,10){\vector(-1,1){15}}%
\put(14.05,0){\circle*{1}}%
\put(14.5,-4.5){\mbox{$\frac{1}{t}$}}%
\curvedashes{}%
\put(0,0){\arc[-1](4,0){135}}%
\setlength{\overhang}{0.5\curvelength}
\curvedashes[.09\curvelength]{1,1}%
\put(0,0){\arc(4,0){135}}%
\put(3,4.2){\mbox{$\alpha$}}
\end{picture}
\end{center}
\caption{The contour $\Gamma_{\alpha,t}$}\label{Gammafigure}
\end{figure}

Observe that on the circular arc of $\Gamma_{3\pi/4,t}$ we have $|e^{tz}|\leq e$, while on the rays $|e^{tz}|\leq e^{-t|z|/2}$ where $|z|\geq t^{-1}$.  It is then easy to verify~\eqref{conditforhintegralsokernelssame}, for example with $a=1$.  We also note that $h(-\lambda_{j})=e^{-t\lambda_{j}}$, so that~\eqref{sumabilityconditonhlambdasokernelssame} is just finiteness of $\sum_{k}2^{kS/(S+1)}e^{-t2^{k}}$.  We conclude that
\begin{equation*}
    p_{t}(x,y)=\frac{1}{2\pi i} \int_{\Gamma_{3\pi/4,t}} G^{(z)}(x,y)e^{zt} \,dz
    \end{equation*}
and can apply our estimates of $G^{(z)}$ to estimate $p_{t}$.
\begin{theorem}\label{heatkernelestimates}
On a pcfss set with regular harmonic structure there is $\kappa_{7}>0$ such that
\begin{equation*}
    |p_{t}(x,y)|
    \lesssim t^{-S/(S+1)} \exp\Bigl(-\kappa_{7} d_{k(t^{-1})}(x,y) \Bigr).
    \end{equation*}
\end{theorem}
\begin{proof}
We make the trivial computation
\begin{equation*}
    \left|\int_{\Gamma_{\alpha,t}}  G^{(z)}(x,y)e^{zt}\, dz \right|
    \leq \sup_{z\in\Gamma_{\alpha,t}}|G^{(z)}(x,y)| \int_{\Gamma_{\alpha,t}} |e^{zt}|\, |dz|
    \lesssim \frac{1}{t} \sup_{z\in\Gamma_{\alpha,t}}|G^{(z)}(x,y)|
    \end{equation*}
and use the bound $|G^{(z)}(x,y)|\lesssim |z+1|^{-1/(S+1)} \Phi(\frac{\pi}{4},|z|,x,y)$ in the sector
$A_{\pi/4}$ from Theorem~\ref{Allestimatesoffnegaxis}, and the fact that $|z|\geq t^{-1}$ on the contour.
\end{proof}
This is the same upper bound as in Theorem~1.1 of~\cite{HamblKumag1999PLMS}, though their result is stronger because it also includes a lower estimate from which it follows that this upper bound is sharp.

Note that in the case of an affine nested fractal, where by Proposition~\ref{affinenesteedversionofdk} and Definition~\ref{defnofklambda} we have $d_{k(\lambda)}\simeq (1+\lambda)^{\gamma/(S+1)} R(x,y)^{\gamma}$ this becomes
\begin{equation*}
    |p_{t}(x,y)|
    \lesssim t^{-S/(S+1)} \exp\Bigl(-\kappa_{7} R(x,y)^{\gamma} (1+t^{-1})^{\gamma/(S+1)} \Bigr)
    \end{equation*}
which should be compared to the upper bound in Theorem~1.1 of~\cite{FitzHamKum1994CMP} or Theorem~6.1 of~\cite{HamblKumag1999PLMS}.  Our $\gamma$ corresponds to $\frac{\gamma'(S+1)}{S+1-\gamma'}$ in the latter reference.

It is also worth mentioning that the above method gives estimates on the derivatives of the heat kernel, simply because differentiation in $t$ introduces a polynomial power of $z$ into the integrand.
\end{example}

\begin{example}[$e^{w\Delta},w\in\mathbb{C}$]
Consider $h_{w}(\Delta)$ for $h_{w}(z)=e^{wz}$ and $w\in\mathbb{C}$.  In order that that~\eqref{sumabilityconditonhlambdasokernelssame} holds it is necessary and sufficient that $w=|w|e^{i\beta}$ with $|\beta|<\frac{\pi}{2}$.  For ~\eqref{conditforhintegralsokernelssame} we must be careful to select a suitable contour.  If we use $\Gamma_{\alpha,t}$ then on the radial parts of the contour we have $|h(z)|=e^{|zw|e^{i(\beta\pm\alpha)}}$, with the $\pm$ sign determining whether we are on the upper or lower ray.  Then~\eqref{conditforhintegralsokernelssame} holds for some (indeed any) $a>\frac{S}{S+1}$ iff $\beta\pm\alpha\in(\frac{\pi}{2},\pi)$, meaning that rotating the boundary rays of the sector by $\beta$ still leaves them in the left half plane.  This is true iff $|\pi-\alpha|<\beta$.
With these constraints, following the proof of the previous theorem with $t=|w|$ we have a bound of the form
\begin{equation*}
    \bigl| H_{w}(x,y) \bigr|
    \lesssim |w|^{-S/(S+1)}  \bigl( \cot |\beta|\bigr) \exp \biggl(- \kappa_{5} d_{k(|w|^{-1})}(x,y) \sin (c|\beta|)\biggr).
    \end{equation*}
\end{example}

\begin{example}[Complex powers]
If we take the usual logarithm with branch cut on the negative real axis then $h(z)=\exp(w\log (-z))=(-z)^{w}$ is analytic in a neighborhood of the spectrum of $\Delta$.  We cannot use one of the $\Gamma_{\alpha,t}$ contours, but can take the rays at angle $\pm\alpha$ from radial distance $\frac{\lambda_{0}}{2|\cos\alpha|}$, where $-\lambda_{0}$ is the smallest Dirichlet eigenvalue, and connect these with a vertical line through $-\frac{\lambda_{0}}{2}$ on the real axis.  This change does not affect the convergence of~\eqref{conditforhintegralsokernelssame}, which depends only on the behavior of $h(z)$ on the rays.  On these $|h(z)|$ is dominated by $e^{|\im(w)|(\pi-|\alpha|)}|z|^{\re(w)}$, so~\eqref{conditforhintegralsokernelssame} is satisfied if $a+\re(w)<0$ for some $a>\frac{S}{2(S+1)}$, hence we need only $\re(w)<-\frac{S}{2(S+1)}$.  Applying the same estimate to $h(-\lambda_{j})$ we see from the Weyl estimate that
\begin{equation*}
    \sum_{j}|h(-\lambda_{j})|^{2}
    \leq C \sum_{j} \lambda_{j}^{2\re(w)}
    \leq C\sum_{k} 2^{kS/(S+1)} 2^{2\re(w)(k+1)}
    <\infty
    \end{equation*}
when $\re(w)<-\frac{S}{2(S+1)}$.  This is sufficient to conclude that $H(x,y)$ and $\tilde{H}(x,y)$ converge to the same limit in $L^{2}(X\times X)$.  If we want convergence of $\tilde{H}(x,y)$ also in $L^{\infty}(X\times X)$ then we instead need $\re(w)<-\frac{S}{S+1}$.

In estimating $|H(x,y)|$ along the rays using Theorem~\ref{Allestimatesoffnegaxis} we do not need any assumptions on $w$.  We have
\begin{align*}
    \lefteqn{\int_{\text{Ray}} |h(z)G^{(z)}(x,y)|\, |dz| }\quad&\\
    &\lesssim e^{|\im(w)|(\pi-|\alpha|)} \biggl( 1+\tan\Bigl(\frac{\pi+|\alpha|}{4} \Bigr) \biggr) \int_{\Gamma} |z|^{\re(w)-1/(S+1)} \exp \biggl(- \kappa_{5} d_{k(|z|)}(x,y) \sin \Bigl( \tilde{c}\frac{\pi-|\alpha|}{\pi+|\alpha|} \Bigr)\biggr)\, |dz|
    \end{align*}
however $d_{k(|z|)}(x,y)\gtrsim (1+|z|)^{1/(S+1)} R(x,y)$ from the discussion preceding Lemma~\ref{boundonratioofchemicaldistancesfordifferentscales} and Definition~\ref{defnofklambda}, so this integral is bounded by a Gamma function.  No matter how large a negative value of $\re(w)$ we have, we cannot obtain anything better than $|H(x,y)|\simeq1$ because the contour integral for $H(x,y)$ includes integration along the vertical line on which the most we can say is that $h(z)G^{(z)}$ is bounded by a constant depending on the fractal and harmonic structure and the length is bounded in the same way.  Abusing notation to use $\Gamma(\cdot)$ for the Gamma function as well as the contour we find that our bound is
\begin{equation*}
    \bigl|H(x,y)\bigr|
    \lesssim 1+ e^{|\im(w)|(\pi-|\alpha|)} \biggl( 1+\tan\Bigl(\frac{\pi+|\alpha|}{4} \Bigr) \biggr) \biggl( \csc\Bigl( \tilde{c}\frac{\pi-|\alpha|}{\pi+|\alpha|} \Bigr) \frac{1}{R(x,y)} \biggr)^{S+(S+1)\re(w)}  \Gamma\Bigl(S+(S+1)\re(w)\Bigr).
    \end{equation*}
If $|\im(w)|$ is larger than $S+(S+1)\re(w)$ then taking $\alpha$ so $|\im(w)|(\pi-|\alpha|)=S+(S+1)\re(w)$ we find that the terms involving $\alpha$ cancel with the Gamma factor to leave only $|\im(w)|^{(S+1)(\re(w)+1)}$.  If $\im(w)$ is smaller than this then we may take $\alpha=\pi-1$ and bound by a multiple of $e^{\im(w)}\Gamma\Bigl(S+(S+1)\re(w)\Bigr)$.  In either case we get no worse than
\begin{equation*}
    \bigl|H(x,y)\bigr|
    \lesssim 1 + (S+(S+1)|w|)^{(S+1)(\re(w)+1) } R(x,y)^{1-(S+1)(\re(w)+1)}.
    \end{equation*}

We further remark that if the fractal is affine nested, so that $d_{k}(x,y)\simeq e^{k\gamma}R(x,y)^{\gamma}$ for some $\gamma\in \Bigl[\frac{1}{S+1}, \frac{1}{2}\Bigr]$ as in Proposition~\ref{affinenesteedversionofdk}, then we can improve the bound to
\begin{equation*}
    \bigl|H(x,y)\bigr|
    \lesssim 1 + (S+(S+1)|w|)^{((S+1)(\re(w)+1)+\gamma-1)/\gamma } R(x,y)^{(1-(S+1)(\re(w)+1))/\gamma}.
    \end{equation*}
Some results related to these appear in~\cite{IonesRoger}, where they were obtained for affine nested fractals by recognizing the kernels as being of Calder\'{o}n-Zygmund type.
\end{example}

%%%%%%%%%%%%%%%%%%%%

\providecommand{\bysame}{\leavevmode\hbox to3em{\hrulefill}\thinspace}
\providecommand{\MR}{\relax\ifhmode\unskip\space\fi MR }
\providecommand{\MRhref}[2]{%
  \href{http://www.ams.org/mathscinet-getitem?mr=#1}{#2}
}
\providecommand{\href}[2]{#2}

\end{document}